\documentclass[10pt,leqno]{article}

\usepackage{amsmath,amsfonts,amscd,amssymb,theorem}
\usepackage{xypic}
\xyoption{all}

\long\def\comment#1\endcomment{}

\makeatletter
\begingroup
\gdef\th@dotted{\normalfont\itshape
  \def\@begintheorem##1##2{%
        \item[\hskip\labelsep \theorem@headerfont ##1\ ##2.]}%
\def\@opargbegintheorem##1##2##3{%
   \item[\hskip\labelsep \theorem@headerfont ##1\ ##2\ (##3).]}}
\endgroup
\makeatother

\theoremstyle{dotted}

\newtheorem{theorem}{Theorem}[section]
\newtheorem{lemma}[theorem]{Lemma}
\newtheorem{conjecture}[theorem]{Conjecture}
\newtheorem{prop}[theorem]{Proposition}
\newtheorem{corr}[theorem]{Corollary}

\newtheorem{prob}[theorem]{Problem}
\newcommand{\dis}{{\displaystyle}}

\makeatletter
\begingroup
\gdef\th@upshape{\normalfont
  \def\@begintheorem##1##2{%
        \item[\hskip\labelsep \theorem@headerfont ##1\ ##2.]}%
\def\@opargbegintheorem##1##2##3{%
   \item[\hskip\labelsep \theorem@headerfont ##1\ ##2\ (##3).]}}
\endgroup
\makeatother

\theoremstyle{upshape}

\newtheorem{defn}[theorem]{Definition}
\newtheorem{remark}[theorem]{Remark}

\makeatletter
\renewcommand{\subsection}{\@startsection{subsection}{2}{0pt}{-3ex
plus -1ex minus -0.2ex}{-2mm plus -0pt minus
-2pt}{\normalfont\bfseries}} \makeatother

\makeatletter
\@addtoreset{equation}{section}
\makeatother

\newcommand{\cntrct}                
{\hspace{2pt}\raisebox{1pt}{\text{$\lrcorner$}}\hspace{2pt}}

\newcommand{\proof}[1][Proof.]{\smallskip\noindent{\em #1}}
\def\endproof{\hfill\ensuremath{\square}\par\medskip}

\def\eqref#1{\thetag{\ref{#1}}}

\let\latexref=\ref
\def\ref#1{{\normalfont{\latexref{#1}}}}

\newcommand{\wt}{\widetilde}
\newcommand{\wh}{\widehat}
\newcommand{\op}{\operatorname}
\newcommand{\ratto}{\dasharrow}     

\newcommand{\idot}{{\:\raisebox{1pt}{\text{\circle*{1.5}}}}}
\newcommand{\hdot}{{\:\raisebox{3pt}{\text{\circle*{1.5}}}}}
\newcommand{\C}{{\mathbb C}}
\newcommand{\R}{{\mathbb R}}
\newcommand{\Q}{{\mathbb Q}}
\newcommand{\Z}{{\mathbb Z}}
\newcommand{\A}{{\cal A}}
\newcommand{\B}{{\cal B}}
\newcommand{\F}{{\cal F}}
\newcommand{\N}{{\cal N}}
\newcommand{\M}{{\cal M}}
\newcommand{\LL}{{\cal L}}
\newcommand{\T}{{\cal T}}

\newcommand{\X}{{\mathfrak X}}
\newcommand{\Y}{{\cal Y}}
\newcommand{\V}{{\cal V}}
\newcommand{\Pp}{{\sf P}}
\newcommand{\Hf}{{\sf H}}
\newcommand{\ZZ}{{\cal Z}}
\newcommand{\Zf}{{\sf Z}}

\newcommand{\calo}{{\cal O}}

\newcommand{\m}{{\mathfrak m}}
\newcommand{\g}{{\mathfrak g}}
\newcommand{\p}{{\mathfrak p}}
\newcommand{\h}{{\mathfrak h}}

\newcommand{\eps}{\varepsilon}
\renewcommand{\phi}{\varphi}
\newcommand{\QQ}{{\mathsf{Quant}}}

\newcommand{\qis}{\cong}
\newcommand{\smsh}{\operatorname{{\bf \#}}}

\renewcommand{\dim}{\operatorname{\sf dim}}
\newcommand{\codim}{\operatorname{\sf codim}}
\newcommand{\gr}{\operatorname{\sf gr}}
\newcommand{\id}{\operatorname{\sf id}}
\newcommand{\rk}{\operatorname{\sf rk}}
\newcommand{\sh}{\operatorname{\sf sh}}
\newcommand{\sgn}{\operatorname{\sf sign}}
\newcommand{\mmod}{\operatorname{\sf{-mod}}}
\newcommand{\Hom}{\operatorname{Hom}}
\newcommand{\Ext}{\operatorname{Ext}}

\newcommand{\RHom}{\operatorname{{\bf R}{\cal H}{\it om}}}

\newcommand{\RG}{\operatorname{{\bf R}\Gamma}}
\newcommand{\HH}{\operatorname{{\mathbb H}}}
\newcommand{\hhom}{\operatorname{{\cal H}{\it om}}}
\newcommand{\Der}{\operatorname{Der}}
\newcommand{\Def}{\operatorname{Def}}
\newcommand{\Ker}{\operatorname{{\sf Ker}}}
\newcommand{\Spec}{\operatorname{Spec}}
\newcommand{\Proj}{\operatorname{Proj}}
\newcommand{\Spf}{\operatorname{Spf}}

\newcommand{\Har}{\operatorname{{\sf Har}}}

\newcommand{\HP}{\operatorname{{\cal H}{\cal P}}}

\newcommand{\aff}{^{{\tt{aff}}}}
\newcommand{\into}{\,\,\hookrightarrow\,\,}
\newcommand{\onto}{\,\,\twoheadrightarrow\,\,}
\newcommand{\bn}{{\mathbf{n}}}
\newcommand{\HHH}{H\!H}
\newcommand{\om}{\omega}
\newcommand{\yvg}{{V/G}}

\newcommand{\base}{B}
\newcommand{\scx}{$\text{symplectic convex }$}
\newcommand{\iso}{{\;\stackrel{_\sim}{\longrightarrow}\;}}

\title{Poisson deformations of symplectic quotient singularities}

\author{Victor Ginzburg and Dmitry Kaledin}

\begin{document}

\maketitle

\begin{abstract}{\footnotesize  We  establish a connection
between smooth symplectic  resolutions and  symplectic deformations
of a (possibly singular)  affine Poisson variety.

 In particular, let
 $V$ be a finite-dimensional complex symplectic
vector space and $G\subset Sp(V)$ a finite subgroup.
Our main result says that the so-called
{\it Calogero-Moser deformation} of
the orbifold   $V/G$ is,  in an appropriate sense, a {\it versal} 
Poisson  deformation.
That enables us to
determine the algebra
structure on the  cohomology  $H^\hdot(X,\C)$
of any smooth  symplectic  resolution $X\onto V/G$
(multiplicative McKay correspondence).
 We prove further that if $G\subset GL(\h)$ is an irreducible   Weyl
group and $V=\h\oplus\h^*$, then 
no smooth  symplectic  resolution of $V/G$
exists unless $G$ is of types
$\mathbf{A},\mathbf{B},\mathbf{C}.$}
\end{abstract}
{\footnotesize \tableofcontents}

\section{Main results.}

\subsection{Introduction.}
Let $V$ be a finite-dimensional symplectic vector space over $\C$,
and $G \subset Sp(V)$ a finite subgroup. The quotient ${\yvg}$ has a
natural structure of an irreducible affine algebraic variety with
coordinate ring $\C[V/G]:=\C[V]^G$, the subalgebra of $G$-invariant
polynomials on $V$. The algebraic variety $\yvg$ is normal and
Gorenstein. However, it is always {\em singular} whenever $G\neq
\{1\}$.

The symplectic structure on $V$ gives rise to a
Poisson algebra structure on $\C[V]$ which induces, by restriction,
a Poisson algebra structure on $\C[V]^G$. Thus, the space ${\yvg}$
becomes a Poisson variety. 

Recall that 
a resolution of singularities $X\to Y$ of
 an irreducible Gorenstein variety $Y$ is called {\it crepant}
provided the smooth manifold $X$ has trivial canonical class.
Assume now $Y$ has a Poisson structure which is generically non-degenerate
(i.e., symplectic). Then a resolution of singularities $X\to Y$ is
called {\it symplectic} provided the pull-back of the symplectic 2-form
on the generic locus of $Y$ extends to a (non-degenerate)
 symplectic 2-form on the whole of $X$. Any symplectic resolution is 
crepant. Conversely, it was shown in \cite{K1} that
 any crepant  resolution $X \to Y$ is 
necessarily symplectic, i.e.,  the pull-back of  symplectic 2-form
on the generic locus of $Y$ automatically extends to a non-degenerate
 2-form on $X$.  
\bigskip

In this paper we study crepant (equivalently, symplectic) resolutions $X
\to \yvg$ of a symplectic quotient singularity ${\yvg}$. Our first
result concerns the existence of such resolutions. It has been shown
by M. Verbitsky  \cite{V} that if the quotient $V/G$ admits a
crepant  resolution, then the group $G$ must be generated
by the so-called {\em symplectic reflections} (see \cite{V} or
Definition~\ref{symplectic_refl} below). 
\medskip

\noindent
{\bf Example.\;} If $\dim V =2$ then any finite subgroup
$G\subset SL(V)$ is  generated by  symplectic reflections,
and  ${\yvg}$ is the so-called du Val surface
singularity. It is well-known that there exists a canonical minimal
resolution of singularities $X \to \yvg$,
which is at the same time a symplectic resolution.
\medskip

An important series of groups  generated
by symplectic reflections are provided by Coxeter groups.
Specifically, let
$\h$ be the complexification of a euclidean real
vector space with root system of some Dynkin type, and let $G \subset
GL(\h)$ be the corresponding Weyl group. Further let $\C^2$ be the
standard space with a fixed volume  $2$-form, and set
$V:=\h \otimes \C^2$.  The tensor product of the symmetric bilinear
form (coming from the euclidean inner product) on $\h$ and the
2-form on $\C^2$ gives a nondegenerate 2-form on $V$, hence, makes
$V$ a symplectic vector space.  The group $G$ acts on $V=\h \otimes
\C^2$, via the action on the first factor, by symplectic
automorphisms.

For types $A_n$, $B_n$, $C_n$ the quotient $V/G$ is known to admit a
crepant resolution  of singularities (of Hilbert scheme type,
see e.g. \cite{Ku}).  Using the results
of I. Gordon [Go], we resolve the existence question
for other Dynkin graphs to the negative:

\begin{theorem}\label{g2}
For a root system of type $D_n$, $E_n$, $F_4$ or $G_2$, the quotient
$V/G$ where $V=\h \otimes
\C^2$, does not admit a  resolution of singularities
$X \to V/G$ with trivial canonical bundle.
\end{theorem}

Following \cite{arg},
 define an increasing filtration
$F_\idot(\C[G])$ on the group algebra $\C[G]$
by letting $F_k(\C[G]),$ $k\geq 0,$ be the
$\C$-linear span of the elements $g\in{G}$ such that
$\rk(\id_{_V}-g)\leq k$. (Thus, $1\in F_0(\C[G])$ and
symplectic reflections belong to $F_2(\C[G])$). This filtration is
clearly compatible with the algebra structure on $\C[G]$. Let
$F_\idot(\Zf{G})$ denote the induced filtration on $\Zf{G}$, the
center of $\C[G]$.  Write $\gr^F_\idot(\Zf{G})$ for the
corresponding associated graded algebra.

Let  $X \to \yvg$ be a
symplectic resolution. Our second result describes the
algebra structure in the cohomology $H^\hdot(X,\C)$ of the manifold
$X$.

\begin{theorem}[Multiplicative McKay correspondence]\label{vasserot}
Let $X \onto V/G$ be a resolution  of singularities 
with trivial canonical bundle.  Then there is a canonical graded algebra isomorphism:
$H^\hdot(X,\C)\cong \gr^F_\idot(\Zf{G})$. In particular, $X$ has no
odd (rational) cohomology.
\end{theorem}

We note that the dimension equality: $\dim H^i(X,\C)=
\dim \gr_i (\Zf{G})$ has been known for some time,
see \cite{batyrev}. Later on, an explicit basis
in the Borel-Moore homology $H_\hdot^{\tt{BM}}(X,\C)$
parametrized by conjugacy classes in $G$
was constructed in \cite{K}.
This  constitues the so-called {\em generalized McKay
correspondence} and amounts, in our notation,
to a linear bijection $\gamma_*: H_{d-\idot}^{\tt{BM}}(X,\C)\iso
\gr_\hdot (\Zf{G})$,
where $d=\dim_{\mathbb{R}} X$.

The multiplication structure in the cohomology was first computed independently
by Vasserot \cite{Vas} and Lehn-Sorger \cite{LS} (cf. also
\cite{wang}) in the case when $X$ is the Hilbert scheme of $n$ points
on $\C^2$. At the same time,
Vasserot (and the first author) conjectured 
the existence of a natural  algebra isomorphism:
$\gamma^*: \gr^F_\idot(\Zf{G})\iso H^\hdot(X,\C)$
 for arbitrary symplectic quotient
singularities; our result proves this conjecture.

The natural $\C^*$-action on $V$ by dilations induces a $\C^*$-action on
$V/G$. The latter may be canonically lifted to the symplectic resolution
$X$, cf. \cite{K1} and Proposition \ref{omega.0} below.
Recall further that $\C^*$-equavariant cohomology
of a $\C^*$-variety is an algebra over $\C[u]$, the cohomology
of the classifying space.  Further, using the filtration
$F_\idot(\Zf{G})$ one forms the corresponding graded
{\it Rees algebra} $\text{\sl{Rees}}_\idot(\Zf{G}):=
\sum_i\, F_i(\Zf{G})\cdot u^i\subset \Zf{G}[u]$.
\begin{conjecture}[$\C^*$-equavariant cohomology] 
There is a canonical graded $\C[u]$-algebra isomorphism:
$H^\hdot_{\C^*}(X,\C)\,\cong\, \text{\sl{Rees}}_\idot(\Zf{G}).$
\end{conjecture}
In the special case of the Hilbert scheme of $n$ points in $\C^2$ 
a proof of this conjecture is implicitly contained in \cite{Vas}.

Theorem \ref{vasserot}
leads further to an interesting question of
computing the Poin-car\'e duality isomorphism for $X$ in terms of 
the group $G$; in other words we propose 

\begin{prob}\label{prob1}
Compute the composite map
$$
\xymatrix{
\gr^F_\idot(\Zf{G})\;\ar[r]_<>(.5){\sim}^<>(.5){\gamma^*}&
\;H^\hdot(X,\C)\;\ar[rr]_{\sim}^<>(.5){\text{Poincar\'e duality}}&&
\;H_{d-\hdot}^{\tt{BM}}(X,\C)\;\ar[r]_{\sim}^<>(.5){\gamma_*}&
\;\gr_\hdot (\Zf{G}).
}
$$
\end{prob}
This map seems to be closely related to the  character table of the group $G$.

\subsection{Orbifold cohomology and Quantization.}\label{remark_orbi}
Let $M$ be an arbitrary smo\-oth symplectic algebraic variety,
and
$G$ a finite group of  symplectic automorphisms of $M$.
There are many examples, cf. \cite{Ba}, \cite{BKR}, in which, 
given a crepant resolution $X\to M/G$ one has
a canonical equivalence of triangulated categories:
$D^b(Coh(X))\cong D^b(Coh_G(M))$. In such a case, taking
the Grothendick groups of both categories, one obtaines
an isomorphism: $K(X)\cong K_G(M)$ of the algebraic K-groups.

In the special case where $V=M$ is a symplectic vector space
and $G\subset Sp(V),$ the Borel-Moore homology group
$H^{\tt{BM}}_\idot(X,\C)$ is known, by \cite{K}, to be
spanned by the algebraic cycles.
Hence the Chern character map gives a ring isomorphism
$ch: \C\otimes K(X)\iso H^\hdot(X,\C)$.
One also has the Thom isomorphism: $K_G(V)\cong R(G)$,
where $R(G)$ stands for the representation ring of $G$.
Thus, in addition to Problem \ref{prob1} we arrive at
the following

\begin{prob}\label{prob2}
Compute the composite map (cf. also \cite[Problem 17.11]{EG}):
$$
\xymatrix{
\Zf{G}\stackrel{\sim}{\to} \C\otimes R(G)\stackrel{\sim}{\to}\C\otimes
K_G(V)\stackrel{\sim}{\to}
\C\otimes K(X)\ar[r]^<>(.5){\sim}_<>(.5){ch}&H^\hdot(X,\C)
\ar[r]^<>(.5){\sim}_<>(.5){\ref{prob1}}&\gr (\Zf{G}).
}
$$
\end {prob}
Note that the group $\Zf{G}$ on the left is viewed as the
algebra of class-functions on $G$ with pointwise multiplication.
Note further, that the isomorphism $K(X)\cong K_G(V)$
is {\it not} compatible with the ring structures.
\smallskip

Associated with a finite
group action on a Calabi-Yau manifold $M$,
one can introduce an  orbifold (= `stringy')
cohomology  $H^\hdot_{\tt{orb}}(M;G)$, see \cite{Ba}.
We consider the special case where $M$ is a holomorphic
symplectic manifold and the group $G$ acts by symplectic automorphisms.
In such a case the definition reads:
\begin{equation}\label{Horb}
H^\hdot_{\tt{orb}}(M;G)\,:=\,\Bigl(\bigoplus\nolimits_{g\in G}\,\,
H^{\hdot-\dim M^g}(M^g)\Bigr)^G,
\end{equation}
where $M^g$ denotes the fixed point set of 
$g\in G$, and where $H^{\hdot-\dim M^g}(M^g)$
is a shorthand notation for 
$\bigoplus_\alpha\, H^{\hdot-\dim M_\alpha}(M_\alpha)$, a
direct sum
ranging over the set of connected components, $M_\alpha$,
of the manifold $M^g$. Further, for any $g,h\in G$, one has a
cup-product pairing
$\cup: H^\hdot(M^g)\times H^\hdot(M^h)\to H^\hdot(M^g\cap M^h)$,
and  the Gysin map
$\imath_*: H^\hdot(M^g\cap M^h)
\to H^\hdot(M^{gh}),$
induced by the imbedding $\imath: M^g\cap M^h\into M^{gh}.$ Following
Ruan, cf. \cite{R1}, \cite{R2},
Fantechi and G\"ottsche  \cite{FG} have introduced
a certain cohomology class $c(g,h) \in H^\hdot(M^g\cap M^h)$,
and showed that the assignment
$$H^\hdot(M^g) \times H^\hdot(M^h)\longrightarrow
H^\hdot(M^{gh}),\quad
a,b\longmapsto \imath_*(a\cup b \cup c(g,h))$$
gives rise to an associative product on
the direct sum in the RHS of \eqref{Horb}.
This product puts
a  graded algebra structure
on orbifold cohomology.

It is known further, see \cite{Ba}, \cite{DL} and also \cite{Ba},
that, given  an
arbitrary  Calabi-Yau orbifold
$M/G$ and a smooth  crepant  resolution
 $X\to M/G$, there is a graded space isomorphism
$H^\hdot_{\tt{orb}}(M;G)\cong H^\hdot(X)$.
Moreover, it has been conjectured in \cite{CR},\cite{R1} and \cite{FG},
that  there is a
  graded  {\em algebra} isomorphism
$H^\hdot_{\tt{orb}}(M;G)\cong H^\hdot(X)$,
provided 
$M$ is a
holomorphic symplectic manifold
with symplectic $G$-action and  $X\to M/G$ is
 a symplectic resolution.

The above conjecture is supported by
the special case where  $M=V$ is
a symplectic vector space and  $G\subset Sp(V)$.
Then each fixed point set $V^g$ is contractible,
and formula \eqref{Horb} reduces
to $H^\hdot_{\tt{orb}}(V;G)=\gr_\idot
(\Zf{G})$
(lemma \ref{easy} from \S6 below
 insures compatibility of the gradings on both sides).
Thus, the conjectured algebra 
isomorphism $H^\hdot_{\tt{orb}}(M;G)\cong H^\hdot(X)$
becomes nothing but our Theorem \ref{vasserot} above.

For a general algebraic symplectic
manifold $M$, the  above
conjecture may be approached as follows.
The symplectic form on $M$ makes
the structure sheaf into a sheaf,
$(\calo_M,\,\{-,-\}_M),$
of Poisson algebras.
Consider a  {\em deformation-quantization}
of $\calo_M,$ i.e., a sheaf  $\QQ_M$ of
locally free complete
$\C[[\eps]]$-algebras
with a $G$-equivariant star-product
$a,b \mapsto a\star b,$ such that $a\star b-
b\star a \equiv \eps\cdot\{a,b\}_M\,({\mathsf{mod}}\,\eps^2).$
Let $\QQ_M[\frac{1}{\eps}]$
be the localization of  $\QQ_M$ with respect to $\eps$,
a  sheaf  of $\C((\eps))$-algebras on $M$.
We form the  cross-product
algebra $\QQ_M[\frac{1}{\eps}]\smsh G$.
Let $(\QQ_M[\frac{1}{\eps}]\smsh G)\mmod$ be the abelian category
of coherent $\QQ_M[\frac{1}{\eps}]\smsh G$-modules,
that is, the category of
$G$-equivariant coherent $\QQ_M[\frac{1}{\eps}]$-modules.

On the other hand, given a smooth symplectic resolution
$X\to M/G$, we consider similarly a sheaf $\QQ_X$ 
of $\C[[\eps]]$-algebras, which is a deformation-quantization
of the structure sheaf $(\calo_X,\,\{-,-\}_X), $ viewed
as a sheaf of Poisson algebras on
$X$. Let $\QQ_X[\frac{1}{\eps}]$
be  its localization with respect to $\eps$, 
a sheaf of $\C((\eps))$-algebras on $X$. Write
$\QQ_X[\frac{1}{\eps}]\mmod$ for the abelian category of
sheaves of coherent $\QQ_X[\frac{1}{\eps}]$-modules.

We propose the following conjecture that may be thought of as a
`quantum analogue' of Bridgeland-King-Reid theorem \cite{BKR}.

\begin{conjecture}\label{quant} For appropriate choices of
deformation-quantizations of symplectic manifolds $X$
and $M$, respectively,
there is a category equivalence
$$\QQ_X[\mbox{$\frac{1}{\eps}$}]\mmod\,\simeq\,
(\QQ_M[\mbox{$\frac{1}{\eps}$}]\smsh G)\mmod.$$
\end{conjecture}

\begin{remark} Notice that, unlike the situation
considered in \cite{BKR}, the conjecture above
involves no derived categories. This is 
 somewhat reminiscent of the
`$D$-affineness' property, proved 
by Beilinson-Bernstein \cite{BB} for flag manifolds.
\end{remark}

Now, given a space $U$ and a sheaf of algebras ${\mathcal{A}}_U$
on $U$, define the  Hochschild cohomology of  ${\mathcal{A}}_U$
by the formula $\HHH^\hdot({\mathcal{A}}_U):=
\Ext^\hdot_{\A_U\boxtimes\A_U^{\,op}}(\A_U,\A_U),$
where $\A_U$ is viewed as a sheaf on
the diagonal $U\subset U\times U$.
It was shown in \cite[\S15]{EG} that
Kontsevich's Formality theorem \cite{Kon}
yields a
 graded $\C((\eps))$-algebra isomorphism:
\begin{equation}\label{kon1}
\C((\eps))\otimes
H^\hdot(X)\cong\HHH^\hdot(\QQ_X[\mbox{$\frac{1}{\eps}$}]).
\end{equation}

Further, we expect that there is  a
 graded $\C((\eps))$-algebra isomorphism
\begin{equation}\label{kon2}\C((\eps))\otimes
H^\hdot_{\tt{orb}}(M;G)\cong \HHH^\hdot((\QQ_M[\mbox{$\frac{1}{\eps}$}])\smsh G),
\end{equation}
that may be viewed
as a `quantization' of the isomorphism of
 Proposition \ref{HHorb} (see \S6
below). Assuming this,
the equivalence of Conjecture \ref{quant} would yield
an isomorphism
between the Hochschild cohomology  algebras in the right-hand sides of
\eqref{kon1} and \eqref{kon2}. Hence, the corresponding
left-hand sides should also be isomorphic, and we would get
the desired  algebra isomorphism
$H^\hdot_{\tt{orb}}(M;G)\cong H^\hdot(X).$
\bigskip

\subsection{Poisson deformations.}\label{poidef.defn.sub}
The main idea of this paper, used in particular in the proofs of
Theorem~\ref{g2} and Theorem~\ref{vasserot}, is to relate smooth
symplectic {\em resolutions} of ${\yvg}$ to smooth symplectic {\em
deformations} of ${\yvg}$.
Specifically, a certain canonical  Poisson deformation of the Poisson
algebra $\C[V]^G$, which we propose to call {\em the
Calogero-Moser deformation}, has been introduced  in 
\cite{EG}.  Our
main  Theorem~\ref{versal} claims that -- under some
assumptions, and in an appropriate sense -- the Calogero-Moser
deformation is a versal deformation of $\yvg$ in the class of
Poisson algebras. In the course of proving the Theorem, we establish some
general basic results on Poisson  deformations,
which may be of independent interest.

We will now introduce the necessary definitions and state, one by
one, the technical results leading up to and including
Theorem~\ref{versal}. At the the end of this section, we will show how
the stated results imply Theorem~\ref{g2} and
Theorem~\ref{vasserot}.

Let $A$ be a commutative unital $\C$-algebra with product
$(a,b)\mapsto a\cdot b$, and $R\subset A$ a (unital)
subalgebra. Recall that $A$ is said to be a {\em Poisson algebra
over $R$}, or a {\em Poisson $R$-algebra}, if $A$ is equipped with
an $R$-linear skew-symmetric bracket $\{-,-\}$ that satisfies the
Leibniz rule
\begin{equation}\label{compt.eq}
\{a,(b\cdot c)\} = \{a,b\}\cdot c + \{a,c\}\cdot
b \qquad \text{ for all } a,b,c\in A,
\end{equation}
and the Jacobi identity
\begin{equation}\label{jacobi}
\{a,\{b, c\}\} + \{b,\{c,a\}\} + \{c,\{a,b\}\} = 0
\qquad \text{ for all } a,b,c\in A.
\end{equation}
Geometrically, the embedding $R\into A$ corresponds to a scheme
morphism $f: \Spec A \to \Spec R$, and the Poisson $R$-algebra
structure on $A$ makes each fiber of $f$ a Poisson scheme. In
particular, if $R$ is a local Artin algebra with maximal ideal $\m$,
then the fiber $A/\m$ over the special point $o \in S = \Spec R$ is
a Poisson algebra over $\C$.

\begin{defn}\label{poi.def.defn}
A {\em Poisson deformation} of a Poisson algebra $A$ over the
spectrum $S = \Spec R$ of a local Artin algebra $R$ with maximal
ideal $\m \subset R$ is a pair of a {\em flat} Poisson $R$-algebra
$A_R$ and a Poisson algebra isomophism $A_R/\m \cong A$.
\end{defn}

It is useful to introduce gradings
into the picture. Say that a unital associative algebra $\dis
A=\bigoplus_{n \geq 0}A_n$ is {\em positively graded} if for all
$n,m \geq 0$ we have $A_n \cdot A_m \subset A_{n+m}$, and $A_0 = \C
\cdot 1$.  Let $A$ be a commutative positively graded
$\C$-algebra, and let $R\subset A$ be a graded subalgebra. Fix $l >
0$, a positive integer. We will say that $A$ is a {\em graded
Poisson $R$-algebra of degree $l$} if $A$ is a Poisson $R$-algebra
and, for any $n,m\geq 0,$ we have $\{A_n,A_m\} \subset
A_{n+m-l}$. Repeating literally Definition~\ref{poi.def.defn}, we
obtain the notion of a graded deformation of a graded Poisson
algebra $A$ over the spectrum of a graded local Artin algebra
$R$. Note that our assumption  implies in particular  that the maximal
ideal $\m \subset R$ coincides with the augmentation ideal
$\bigoplus_{n \geq 1}R_n \subset R$.

A  Poisson
$R$-algebra $A_R$ over a complete local $\C$-algebra $\langle R,\m_R
\rangle$ equipped with a Possion isomorphism $A_R/\m_R \cong A$ is
called a
{\em universal formal Poisson deformation} of the Poisson algebra
$A$ if for every Poisson deformation $A_S$ over a local Artin base
$S$ there exists a unique map $\tau:S \to \Spec R$ such that the
isomorphism $A \cong A_R/\m_R$ extends to a Poisson isomorphism
\begin{equation}\label{eq.uni}
A_S \cong \tau^*A_R.
\end{equation}
Analogously, say that a graded Poisson $R$-algebra $A_R$ over a
positively graded $\C$-algebra $R$ equipped with a graded Poisson
isomorphism $A_R/R^{>0} \cong A$ is a {\em universal graded Poisson
deformation} of the Poisson algebra $A$ if for every graded Poisson
deformation $A_S$ over a local Artin base $S$ with an action of
$\C^*$ there exists a unique $\C^*$-equivaraint map $\tau:S \to
\Spec R$ such that the isomorphism $A \cong A_R/R^{>0}$ extends to a
graded Poisson isomorphism \eqref{eq.uni}.

To study Poisson deformations of a given Poisson $\C$-algebra $A$,
we develop a cohomology theory, $HP^\hdot(A)$, the {\em Poisson
cohomology} of $A$. Very much like Hochshcild cohomology of an
associative algebra, first order deformations of a Poisson algebra
$A$ are controlled by the group $HP^2(A)$, while obstructions to
deformations are controlled by the group $HP^3(A)$. The Poisson
cohomology groups are sufficiently functorial; in particular, for a
graded Poisson algebra $A$, each space $HP^k(A)$ carries a natural
additional grading.

\begin{remark}
For Poisson algebras such that $\Spec A$ is {\em smooth}, the theory
of Poisson cohomology goes back to J.-L. Koszul \cite{Ko}, and
J.-L. Brylinski \cite{B}. In the singular case, the approach must be
quite different; it is based on the general operadic formalism (see
\cite{Fr}, and
Appendix below).
\end{remark}

Assume now  that $HP^1(A)=0$ and $\dim
HP^2(A)< \infty$. Then, the result below says that 
 there exists a formal {\em universal} Poisson
deformation of $A$; moreover, under additional
assumptions, in the graded setting the formal universal deformation
comes from a deformation over a base of finite type over $\C$. 

To
make a precise statement, let  $\wh{HP^2(A)}$ denote the completion of
the affine space $HP^2(A)$.

\begin{theorem}[Universal Poisson deformations]\label{poi.def.uni}
Let $A$ be a Poisson algebra. Assume that $HP^1(A)=0$ and that
$HP^2(A)$ is a finite-dimensional vector space over $\C$.
\begin{enumerate}
\item There exists a closed subscheme $S
\subset \wh{HP^2(A)}$ and a Poisson $\C[S]$-algebra $A_S$ which is
a universal formal Poisson deformation of the algebra $A$.
\item\label{graded.uni} Assume in addition that the Poisson algebra
$A$ is positively graded of some positive degree $l$, and that the
induced grading on $HP^2(A)$ is also positive. View $HP^2(A)$ as a
$\C^*$-variety 
via the grading. 

Then there exists a
$\C^*$-stable closed subvariety $S \subset HP^2(A)$
and a positively graded Poisson
$\C[S]$-algebra $A_S$ of degree $l$ which is a universal graded
Poisson deformation of the algebra $A$.
\end{enumerate}
\end{theorem}

The proof of Theorem~\ref{poi.def.uni} is contained in the Appendix,
in Subsection~\ref{def.sub}.

\begin{remark}
We would like to note that the isomorphism \eqref{eq.uni} is not
canonical, nor indeed is it unique. Thus what we obtain is more than
a ``coarse deformation space'' and less than a ``fine moduli
space''. In the language of the stack formalism, our deformation
problems are classified by the quotient stack of some variety $S$ by
a {\em trivial} action of some group $G$. The group in question is
the connected component of unity in the group of Poisson
automorphisms of the algebra $A$.
\end{remark}

\subsection{Deforming symplectic resolutions.}
Let $X\aff:=\Spec H^0(X, {\mathcal{O}}_X)$
denote  the  affinization of
an  algebraic variety $X$. This is an
 affine algebraic variety.
We first introduce
\begin{defn} An irreducible smooth  algebraic
variety $X$ will be called {\it convex} if the canonical
 morphism $X \to X\aff$
is  projective and birational.
The affine algebraic variety $X\aff$ is then
 a normal variety.
\end{defn} 
For a convex symplectic manifold $X$, the symplectic structure on $X$
induces a Poisson bracket
on the algebra $H^0(X, {\mathcal{O}}_X)$, hence, a Poisson structure on $X\aff$.

Our next result concerns deformations of \scx
manifolds. Related results on infinitesemal and formal deformations have
 been
studied in the paper \cite{KV}. 
In order to generalize the
main theorem of \cite{KV} to a global setting,
 we need to assume that the affinization
$Y=X\aff$ is equipped with an  {\em expanding} action of the multiplicative group
$\C^*$. Equivalently, we assume that the algebra
$\C[Y]=H^0(X,\calo_X)$ carries a grading
by  {\em non-negative} integers such that
 the Poisson
bracket on $\C[Y]$ is of some fixed  degree $l > 0,$ that is,
for any homogeneous functions $a_1,a_2 \in \C[Y],$ we have
$$
\deg \{a_1,a_2\} = \deg a_1 + \deg a_2 -l.
$$

With these assumptions, the result below says that any 
convex symplectic manifold can be deformed nicely\footnote{this is
similar to the deformation of the Springer resolution
$\tilde\N \to \N$ (of the nilpotent variety $\N$ in a semisimple
Lie algebra $\g$) provided by Grothendieck's simultaneous resolution
$\tilde\g\to \g$. Here $X=\tilde\N\,,\, X\aff=\N,$
and $\X_\base$ in the Theorem below plays the role of $\tilde\g$.}
 to a smooth affine symplectic manifold without
changing the rational cohomology algebra.
\begin{theorem}\label{vit1}
Let 
$X$ be a \scx variety, and  $X \to Y=X\aff$  the corresponding resolution.
 Assume that
the algebra $\C[Y]$ carries a grading by non-negative integers such
 that the Poisson bracket has degree $l > 0$.  Put $\base:=H^2(X,\C)$ and 
equip the affine space $\base$ with a $\C^*$-action by $z \cdot b =
z^{-l}b$, $z \in \C^*$, $b \in \base=H^2(X,\C)$.

Then there exists a smooth $\C^*$-variety $\X_\base$ and a smooth
$\C^*$-equivariant morphism $\pi:\X_\base\to \base$ such that
\begin{enumerate}
\item $\X_\base$ is a relative symplectic manifold over $\base$, i.e., we
have a relative $2$-form $\om \in H^0(\X_\base,\Omega^2(\X_\base/\base))$
which induces a symplectic structure on each fiber $\X_b\,,\,b\in
\base=H^2(X,\C)$.

\item\label{coho.const} The relative cohomology sheaves $R^k\pi_*\Q$
are constant sheaves on $\base$ for all $k \geq 0$, and the canonical
base change morphism
$$
H^k(\X_b,\Q) \to \left(R^k\pi_*\Q\right)_b
$$
is an isomorphism for every point $b \in \base$.

\item\label{gen.1-1} The affinization, $(\X_\base)\aff,$ is flat over
$\base$. The canonical map $\X_\base \to (\X_\base)\aff$ is projective, birational, and
it is an isomorphism over the generic point $b \in \base$.

\item The special fiber $\X_0$ over $0 \in \base = H^2(X,\C)$ is
isomorphic to $X$ as a symplectic algebraic variety. The special
fiber $(\X_0)\aff$ of the affinization $(\X_\base)\aff$ is isomorphic to 
 $X\aff=Y$ as a Poisson algebraic variety with 
$\C^*$-action.
\end{enumerate}
\end{theorem}

In other words, there exists
a Zariski open, dense subset $\base^{\tt{generic}}\subset B$,
such that
one has a commutative diagram
$$
\xymatrix{
\X_{_{\base^{\tt{generic}}}}\;\ar@{=}[d]^<>(.5){\pi}\ar@{^{(}->}[rr]&&\X_\base\ar[d]^<>(.5){\pi}&
\X_o\ar@{_{(}->}[l]\ar[d]^<>(.5){\pi}\ar@{=}[r]& X\ar[d]\\
(\X_{_{\base^{\tt{generic}}}})\aff\;\ar@{^{(}->}[rr]\ar[d]&&(\X_\base)\aff\ar[d]&
(\X_o)\aff\ar@{_{(}->}[l]\ar[d]\ar@{=}[r]& Y\ar[d]\\
\base^{\tt{generic}}\;\ar@{^{(}->}[rr]&&\base&\{o\}\ar@{_{(}->}[l]\ar@{=}[r]&\{o\}\\
}
$$
In this diagram, the subscript `$\base^{\tt{generic}}$' indicates
restricion to $\base^{\tt{generic}}$ of a scheme over $B$.
Note, in particular, that according to Theorem \ref{vit1}.
the map $ \X_{_{\base^{\tt{generic}}}}\to \base^{\tt{generic}}$
is affine, so the arrow $\X_{_{\base^{\tt{generic}}}}\to 
(\X_{_{\base^{\tt{generic}}}})\aff,$ on the left of the diagram
above is an isomorphism.

The proof of Theorem~\ref{vit1} is contained in
Subsection~\ref{glb.sub}.


\subsection{Comparison.}\label{comp.sub}
Let $X$ be a convex symplectic manifold whose
af\-fi\-nization $Y=X\aff$ is a positively-graded Poisson algebra of
degree $l$, so that the assumptions of Theorem~\ref{vit1} are
satisfied. Assume in addition that $HP^1(\C[Y])$ $=0$, $\dim
HP^2(\C[Y]) < \infty$, and the natural grading on $HP^2(\C[Y])$ is
positive, so that we can apply Theorem~\ref{poi.def.uni}
to get the universal graded Poisson
deformation (of $Y$) over a base $S$. On the other hand,
let
$\X_B/B$ be the deformation provided by Theorem~\ref{vit1}.
 Then the affinization
$(\X_B)\aff$ is by construction a positively-graded flat Poisson
$\calo_B$-algebra.
Hence, by Theorem~\ref{poi.def.uni}~\ref{graded.uni},
 we get the  classifying map
$$
\tau:B \to S := \text{base of universal deformation}.
$$
In
Subsection~\ref{glb.sub} we will prove
\begin{prop}\label{cmp}
The classifying map $\tau: B \to S$ is a finite map onto an
irreducible component of the variety $S$.
\end{prop}

Now, fix  a symplectic vector space $V$ and a finite
group $G\subset Sp(V)$. 

\begin{defn}[\cite{EG}]\label{symplectic_refl} 
An element $g\in G$ is called a {\sl symplectic
reflection} if $\;\;\rk(\id_V-g)=2.$
\end{defn}
Let $\Sigma$ denote the set of symplectic reflections in $G.$
The group $G$ acts on $\Sigma$ by conjugation, and we put
$$
\bn:= \text{\sl number of $G$-conjugacy classes in $\Sigma$}.
$$

Introduce a grading on the polynomial algebra
$\C[V]$ by assigning degree $1$ to all linear functions. This turns
$\C[V]$ into a positively graded Poisson algebra of degree
$l=2$. Set $Y=V/G$.
Both the grading and the Poisson bracket descend to the
algebra  $\C[Y] = \C[V]^G$.
In  Section~\ref{hp2.sec} we will prove
\begin{prop}\label{vit2} 
We have $HP^1(\C[V]^G)= 0$ and $\dim HP^2(\C[V]^G)= \bn$. Moreover,
the natural grading on $HP^2(\C[V]^G)$ is positive.
\end{prop}

Let $X\to \yvg$ be a symplectic resolution.
The Betti numbers of $X$ are known by
\cite{batyrev} (see \cite{K} for an alternative proof), in
particular, we have $\dim H^2(X,\C) = \bn$.
Applying Proposition~\ref{cmp} we obtain a
sequence of maps
$$
H^2(X,\C) \overset{\tau}{\to} S \hookrightarrow HP^2(\C[\yvg]),
$$
where the first map is a finite surjection onto an irreducible
component, and the second map is a closed embedding. Since $\dim
HP^2(\C[\yvg]) = \bn = \dim H^2(X,\C)$, we conclude that $S$ is the
whole $HP^2(\C[\yvg])$ and $\tau:H^2(X,\C) \to S$ is a finite
dominant map.

\begin{remark} The composite map $H^2(X,\C) \to HP^2(\C[\yvg])$ is 
not a
bijection. Although both the source and
the target of this map are affine spaces, the map is not linear; it
is a ramified covering, and its differential at $0 \in H^2(X,\C)$
usually vanishes.
\end{remark}

\subsection{Calogero-Moser deformation.} Let $V$ and $G$ be as above.
Let $\Sigma \subset G$ be the set of all
symplectic reflections in $G$, and let $C$ be the vector space of all
$G$-invariant functions $c: \Sigma\to \C$. Clearly, $\dim C=\bn$. We
regard $\C[C]$, the polynomial algebra, as a positively graded
algebra by assigning degree $2$ to linear polynomials.

In \cite{EG}, P. Etingof and the first author have constructed a
certain flat graded Poisson $\C[C]$-algebra $\B$ of degree $2$,
which gives a graded $\bn$-parameter deformation of $\C[V]^G$.  Let
$\B_c$ denote the specialization of $\B$ at a point $c\in C$, and
put $\M_c=\Spec\B_c$. Thus, $\{\M_c\}_{c\in C}$ is a flat family of
Poisson affine algebraic varieties and, by construction, we have
$\M_0=V/G$.  We will refer to $\M_c$ as a {\em Calogero-Moser}
variety with parameter $c$, and we call the projection
$\M:=\Spec\B\onto C$ the {\em Calogero-Moser deformation} of
$V/G$.\footnote{Our terminology is motivated by the special case
where $G=S_n$ is the symmetric group, acting diagonally on
$V=\C^n\oplus\C^n$ by permutation of coordinates.  In that case, the
deformation $\M_c$ of $V/S_n$ is known (see [EG] and references
therein) to be the usual Calogero-Moser space.}

According to Proposition \ref{vit2}, there exists a canonical
classifying map $\kappa: C \to S \subset HP^2(\C[V]^G)$ such that
the Calogero-Moser deformation is obtained by pull-back from the
universal deformation. We will see that the map $\kappa$ is never
bijective. The best we can prove is provided by the following

\begin{theorem}\label{zcsmooth}\label{versal}
If the Calogero-Moser space $\M_c$ is smooth for generic values of
the parameter $c \in C$ then:
\begin{enumerate}
\item The base $S$ of the universal Poisson deformation
of $V/G$  coincides
with the whole vector  space $HP^2(\C[V]^G)$.
\item The classifying map $\kappa: C \to HP^2(\C[V]^G)$ is
surjective and generically \'etale.
\end{enumerate}
\end{theorem}

The proof of this theorem is contained in
Section~\ref{zcsmooth.sec}. Part \thetag{ii} of the Theorem says,
roughly speaking, that generically, every fiber of the universal
deformation is isomorphic to $\B_c$ for a suitable value of the
parameter $c \in C$.  Moreover, for general $c$, there exists only
a finite number of other values $c' \in C$ such that $\B_c \cong
\B_{c'}$. This is what we mean by saying that the Calogero-Moser
deformation is ``versal''.

It is not unreasonable to conjecture that the claim of
Theorem~\ref{zcsmooth} holds in the general situation. Even
stronger, we propose the following.

\begin{conjecture}\label{1-1}
For an arbitrary symplectic quotient singularity $V/G$,
the classifying map $C \to S$ of the  Calogero-Moser  deformation
$\M/C$ is
finite.
\end{conjecture}

The conjecture is motivated by analogy with the case of a symplectic
resolution $X \to \yvg$ and the corresponding deformation
$(\X_B)\aff$. The parameter spaces of these two deformations are the
same vector spaces with the same grading. Moreover, we will see in
Section~\ref{dim2} that they actually coincide when $\dim \yvg = \dim V
= 2$.

In general, view $V/G$ as a $\C^*$-variety, with $\C^*$-action being
induced
from the natural
one on the vector space $V$.
Then, the relation between symplectic resolutions of
$V/G$ and the Calogero-Moser deformation is provided by the
following theorem.

\begin{theorem}[Symplectic  resolutions and Calogero-Moser]\label{res}
Let $X$
$ \to Y=\yvg$ be a symplectic resolution of a symplectic quotient
singularity ${\yvg}$.
\begin{itemize}
\item Let $\Y/S$ be the universal Poisson deformation
of the Poisson variety $Y$, 
\item Let $\X_{\base}$ be the deformation
over $\base = H^2(X,\C)$ provided by Theorem~\ref{vit1},
\item  Let $\M/C$ be
the Calogero-Moser deformation, and denote by $\kappa:C \to S$ its
classifying map.
\end{itemize}
Then for  $c \in C$ general enough, the image $\kappa(c) \in S$ is
generic in the sense of Theorem~\ref{vit1}~\ref{gen.1-1}. More
precisely, for every point $b \in \base$ lying over $\kappa(c) \in S$,
the canonical map $\X_{b} \to \Y_{\kappa(c)}$ is an isomorphism.
\end{theorem}

Put  $\wt{C} = C \times_S \base$ and let $\psi:
\widetilde{C}\to \base$ and $\phi: \widetilde{C}\to C$ be the natural
projections. Then the statement of Theorem~\ref{res} can be
summarized by the following
commutative diagram.

{\small
$${\small
\xymatrix{
\M_{_{\widetilde{C}}}:=\widetilde{C}\times_{_{C}}\M
\enspace\ar[d]_{\phi\times\id_\M}\ar[dr]\ar@{=}[rr]^{\sim}&&
(\X_{_{\widetilde{C}}})\aff\ar[d]\ar[dl]\ar@{=}[r]^<>(.5){\pi_{_{\widetilde{C}}}}&
\X_{_{\widetilde{C}}}:=\widetilde{C}\times_{_{\base}}(\X_{_{\base}})\ar[d]\\
\M\ar[d]_{\text{\parbox[c]{4em}{\sf Calogero-\\[-1.5mm]
Moser}}}&\widetilde{C}
\ar[dl]_{\phi}\ar[dr]^{\psi}&(\X_{\base})\aff\ar[d]&\X_{\base}
\ar[l]_<>(.5){\pi=\pi_{_\base}}\\
C\ar[r]_<>(.5){\kappa} &\enspace
S=HP^2(\C[V]^G)\enspace&\base\ar[l]^<>(.5){\tau}&\\ }
}
$$}

The
proof of Theorem~\ref{res} is contained in Subsection~\ref{res.sub}.

\begin{corr}\label{res=>smooth} 
The existence of a symplectic resolution $X \to V/G$ implies
 that the generic fiber $\M_c$ of the Calogero-Moser
deformation $\M/C$ is smooth.
\end{corr}

\proof{} By definition we have $\M_c \cong \Y_{\kappa(c)}$,
where $\Y_B=(\X_B)\aff$. The
latter is isomorphic, for  $b\in B$ such that $\kappa(c)=\tau(b)\in S$ is
general enough,
 to $\X_{b}$, which is  smooth.
\endproof 

\subsection{Proofs of Theorem~\ref{g2} and Theorem~\ref{vasserot}.}
Let $V$ be a symplectic vector space, let $G \subset Sp(V)$ be a
finite group, and let $X \to Y=V/G$ be a crepant, hence symplectic,
 resolution.

To prove Theorem~\ref{g2}, note that in the assumptions of the Theorem, the
existence of $X$ implies by Corollary \ref{res=>smooth} that the
generic fiber $\M_c$ of the Calogero-Moser deformation $\M/C$ is
smooth. This contradicts  \cite[Proposition 7.3]{Go}, see also
\cite[Proposition 16.4(ii)]{EG}.

To prove Theorem~\ref{vasserot}, consider the universal deformation
$\Y/S$ and the deformation $\X/B$ provided by Theorem~\ref{vit1}.
Consider a general fiber $\M_c$ of the Calogero-Moser deformation
$\M/C$. Choose a point $b \in B$ lying over $\kappa(c) \in S$.  By
definition, we have isomorphisms of fibers $\M_c \cong
\Y_{\kappa(c)} \cong \X_{b}$, which induce algebra isomorphisms of the
rational cohomology:
$$
H^\hdot(\M_c,\Q) \cong H^\hdot(\Y_{\kappa(c)},\Q) \cong
H^\hdot(\X_{b},\Q)\cong H^\hdot(X,\Q),
$$
where the last isomorphism is due to Theorem~\ref{vit1}~\ref{coho.const}.
By Corollary~\ref{res=>smooth}, the generic fiber $\B_c$ is
smooth. By \cite[Theorem 1.8(i)]{EG}, this implies that the
left-hand-side is isomorphic to $\gr_\hdot (\Zf{G})$.\endproof

\subsection*{Acknowledgments.}
{We are grateful to V. Baranovsky,
R. Bezrukavnikov,
P.~Etingof,  and V. Ostrik for many useful discussions. 
We also thank Y. Ruan for bringing the question of the validity of
Lemma \ref{easy} 
to our attention.
The second author was
partially supported by CRDF Award
RM1-2354-MO02.}

\section{Generalities on Poisson deformations.}\label{def}
In this section, we review
  basic results on  the deformation theory of Poisson algebras
that will be used later. We have been unable to find an
adequate reference in the literature, so in the Appendix to this paper
the reader may find some details of  proofs
and precise definitions.
\smallskip

\noindent
{\bf Convention.}\, Given a commutative $\C$-algebra $A$,
we write $\Hom_A\,,\,\otimes_A,$ and $\Lambda_A^k$ for Hom,
tensor product and $k$-th wedge product over $A$, respectively.
If no subscript $A$ is indicated, then the corresponding
functors are understood to be taken  over $\C$,
e.g. $\otimes=\otimes_\C$.

We let $\T_M$ denote the tangent sheaf (or tangent bundle)
of a smooth algebraic variety $M$.
Further, given a finitely generated commutative algebra $A$ we write 
$\Omega^1A$ and $\T(A)=\T(\Spec A)$ for the $A$-modules of (global)
K\"ahler differentials and vector fields on the scheme $\Spec A$,
respectively.
\smallskip

\subsection{Poisson cohomology.} For smooth Poisson manifolds, the
notion of Poisson cohomology (in the differential geometric setting)
is due to J.-L. Koszul \cite{Ko} and
J.-L. Brylinski \cite{B}. In the
algebraic setting, the   Poisson cohomology has been 
introduced by B. Fresse \cite{Fr}. We review and extend it below
(relations with deformation theory were not discussed in 
\cite{Fr}).
For a more general  formalism of deformation theory of an algebra
over an arbitrary operad the reader may consult \cite{KS}.

Let $A$ be a  finitely generated commutative algebra.
 The standard
Lie bracket of vector fields extends  to the
so-called {\em Schouten bracket} $\{-,-\}$ on
$\Lambda^\hdot_A\T(A)$, the space of polyvector fields.
It is well-known that if  $\Spec A$ is smooth
then any  Poisson structure
on $A$ defines (and is defined by) a bivector field $\Theta \in
\Lambda^2_A\T(A)$ which satisfies the integrability condition
$\{\Theta,\Theta\}=0$.

Given a smooth Poisson algebra $A$ one
defines a map
$$
d:\Lambda^\hdot_A\T(A) \to \Lambda^{\hdot+1}_A\T(A)\quad,\quad
a \mapsto da:=\{\Theta,a\}\,.
$$
We have
 $d \circ d = 0$. Taking $d$ as the differential makes
$\Lambda^\hdot_A\T(A)$ into a complex. The cohomology groups of this
complex are called the {\em Poisson cohomology groups} of the
algebra $A$ and denoted by $HP^\hdot(A)$.

\medskip

It turns out that this formalism can be extended to arbitrary, not
necessarily smooth Poisson algebras $A$. The precise constructions
are a little bit technical; we give them in full in the
Appendix. Here we only describe the end result.

Recall that for an arbitrary finite-type commutative algebra $A$,
one can define the so-called {\em Harrison complex} $\Har_\idot(A)$,
a certain canonical complex of free $A$-modules respresenting the
{\em cotangent complex} of the algebra $A$. We recall the precise
construction of $\Har_\idot(A)$ in Subsection~\ref{har.def}. Here we
only note that when the algebra $A$ is smooth, the complex
$\Har_\idot(A)$ has non-trivial cohomology only in degree $0$, and
this non-trivial cohomology module is isomorphic to the module
$\Omega^1A$ of K\"ahler differentials.

Let now $A$ be an arbitrary Poisson algebra, not necessarily smooth.
Heuristically, to extend the Brylinski construction to $A$, one
replaces everywhere the module $\Omega^1A$ of K\"ahler differentials
with the Harrison complex $\Har_\idot(A)$. More precisely,
one considers the exterior powers $\Lambda^k_A\Har_\idot(A)$ of the
complex $\Har_\idot(A)$ of flat $A$-modules and sets
\begin{equation}\label{gr}
DP^{\hdot,k}(A) \cong \Hom_A(\Lambda^k_A\Har_\idot(A),A), \qquad k
\geq 0.
\end{equation}
This defines a canonical bigraded vector space $DP^{\hdot,\hdot}(A)$
and a differential $d\colon DP^{\hdot,\hdot}(A) \to
DP^{\hdot+1,\hdot}(A)$. In particular, we have 
$$DP^{0,2}(A)\cong
\Hom_A(\Lambda^2A \otimes A, A) = \Hom_\C(\Lambda^2A,A).$$
There is  a natural  {\em Gerstenhaber bracket}
on
 $DP^{\hdot,\hdot}(A)$, see \eqref{gerst}:
$$
\{-,-\}:DP^{p,q}(A) \otimes DP^{p',q'}(A)
\to DP^{p+p',q+q'-1}(A), 
$$
that makes  $DP^{\hdot,\hdot}(A)$ a DG Lie algebra (with shifted grading). The
Poisson structure on $A$ defines (and is defined by) an element
$\Theta \in DP^{0,2}(A) = \Hom(\Lambda^2 A,A)$ satisfying $d\Theta=0$
and $\{\Theta,\Theta\}=0$. We call this element {\em the Poisson
cochain}. Given such an element, one defines the differential
$\delta:DP^{\hdot,\hdot}(A) \to DP^{\hdot,\hdot+1}(A)$ by setting
$$\delta:\;DP^{\hdot,\hdot}(A) \to DP^{\hdot,\hdot+1}(A)\quad,\quad
a\mapsto \{\Theta,a\}.
$$
Thus we have a bicomplex $DP^{\hdot,\hdot}(A)$ with differentials
$d$, $\delta$. We define  Poisson cohomology of $A$
to be the cohomology groups  $HP^\hdot(A)$  of
the total complex associated with this bicomplex,
with respect to the total differential
$d+\delta$.

If the algebra $A$ is smooth, then we have a quasiisomorphism
$\Har_\idot(A) \qis \Omega^1A$, so that
$$
DP^{\hdot,k}(A) \qis \Hom(\Lambda^k_A\Omega^1A,A) \cong
\Lambda^k_A\T(A),
$$
and the general Poisson cohomology complex $DP^\hdot(A)$ is
quasiisomorphic to the Brylinski complex $\langle
\Lambda^\hdot_A\T(A),d \rangle$.

If the algebra $A$ is a Poisson graded algebra in the sense of
Subsection~\ref{poidef.defn.sub}, then the Poisson cohomology
bicomplex $DP^{\hdot,\hdot}(A)$ acquires an additional grading,
called the {\em $A$-grading}. However, since the Poisson cochain
$\Theta$ is of degree $l$ with respect to the grading, the
differential $\delta:DP^{\hdot,\hdot}(A) \to DP^{\hdot+1,\hdot}(A)$
does not preserve the $A$-grading but rather shifts it by $l$. To
cure this, we redefine the $A$-grading by shifting it by $(k-1)l$ on
$DP^{k,\hdot}(A)$.

\subsection{Globalization.}
The Poisson cohomology complex $DP^\hdot(A)$ becomes very simple if
the scheme $A$ is not only smooth but also symplectic. In that case
the canonical isomorphism $\Omega^1A \cong \T(A)$ extends to an
isomorphism $\Omega^\hdot(A) \cong \Lambda^\hdot\T(A)$ between the
Brylinski complex $\Lambda^\hdot\T(A)$ and the  de Rham complex
$\Omega^\hdot(A)$. Thus the Poisson cohomology $HP^\hdot(A)$
coincides with the de Rham cohomology $H^\hdot(\Spec A)$ of the
scheme $\Spec A$. It turns out that several features of the de Rham
cohomology formalism extends to the general case.

Firstly, one can define Poisson cohomology with coefficients, which
is analogous to the de Rham cohomology with coefficients in a local
system (or, more generally, in a $D$-module). For this one defines a
{\em Poisson module} $M$ over an arbitrary Poisson algebra $A$ in an
natural way. Then to every Poisson module $M$ one associates a
canonical bicomplex $DP^{\hdot,\hdot}(A,M)$ called the {\em
cohomology complex with coefficients in $M$}. The algebra $A$ is a
Poisson module over itself, and we have $DP^{\hdot,\hdot}(A,A) =
DP^{\hdot,\hdot}(A)$.

Secondly, one generalizes the notion of Poisson cohomology to the
scheme case. One defines a Poisson scheme $X$ and a Poisson sheaf
$\F$ of $\calo_X$-modules in the natural way. To a Poisson scheme
$X$ with a Poisson sheaf $\F$ (or, more generally, to a complex
$\F^\hdot$ of Poisson sheaves) one associates a canonical complex
$\HP^\hdot(X,\F^\hdot)$ of Zariski sheaves on $X$ called the {\em
local Poisson cohomology complex with coefficients in
$\F^\hdot$}. In the particular case $\F = \calo_X$, one obtains the
{\em local Poisson cohomology complex} $\HP^\hdot(X) =
\HP^\hdot(X,\calo_X)$. When the Poisson scheme $X$ is smooth and
symplectic, a Poisson sheaf on $X$ is the same as a $D$-module, the
Poisson cohomology complex $\HP^\hdot(X)$ is the de Rham complex of
the scheme $X$, and the Poisson cohomology complex with coefficients
$\HP^\hdot(X,\F)$ is the de Rham complex of the $D$-module $\F$.

Taking the hyperhomology, one obtains the groups
$$
HP^\hdot(X) = \HH^\hdot(X,\HP^\hdot(X))
$$ 
of global Poisson cohomology of the scheme $X$ and the groups
$$
HP^\hdot(X,\F^\hdot) = \HH^\hdot(X,\HP^\hdot(X,\F^\hdot))
$$ 
of global Poisson cohomogoly of $X$ with coefficients in the complex
$\F^\hdot$. When the scheme $X = \Spec A$ is affine and the complex
$\F^\hdot$ comes from a complex $M^\hdot$ of Poisson $A$-modules, we
have $HP^\hdot(X,\F^\hdot) \cong HP^\hdot(A,M^\hdot)$. When the
scheme $X$ is smooth and symplectic, $HP^\hdot(X,\F^\hdot)$ is the
 singular cohomology $H^\hdot(X,\F^\hdot)$ with coefficients
in the $D$-module $\F^\hdot$. When $X = Y \times Z$ is a product of
two Poisson schemes $Y$ and $Z$, we have the K\"unneth formula
(Proposition~\ref{geom}\thetag{i})
\begin{equation}\label{kun.eq}
HP^k(X) \cong \bigoplus_{p+q=k}HP^p(Y) \otimes HP^q(Z), \qquad k \geq 0.
\end{equation}
The notion of a Poisson sheaf is sufficiently functorial; in
particular, for any morphism $f:X \to Y$ between Poisson schemes,
the direct image $f_*\F$ of a Poisson sheaf $\F$ on $X$ is a Poisson
sheaf on $Y$. Moreover, one can represent the direct image $R^\hdot
f_*\F$ in the derived category by a complex of Poisson sheaves. In
the particular case of an open embedding $j:U \hookrightarrow X$, we
obtain a Poisson structure on the direct image $R^\hdot
j_*\calo_U$. Moreover, in this particular case we have 
$$
\HP^\hdot(X,R^\hdot j_*\calo_U) \qis R^\hdot j_*(\HP^\hdot(U)).
$$
One can also obtain a Poisson structure on the third term
$i_!\calo_Z$ in the exact triangle
$$
\begin{CD}
i_!\calo_Z @>>> \calo_X @>>> R^\hdot j_*\calo_U @>>>
\end{CD}
$$
where $Z \subset X$ is the closed complement to $U \subset X$ and
$i:Z \hookrightarrow X$ is its embedding. We define the {\em Poisson
cohomology $HP^\hdot_Z(X)$ of the scheme $X$ with supports in $Z
\subset X$} by setting $HP^\hdot_Z(X) \cong
HP^\hdot(X,i_!\calo_Z)$. When $X$ is smooth and symplectic, this is
the  ordinary singular cohomology $H^\hdot_Z(X)$ with supports in
$Z$. By definition, for a general Poisson $X$ we have the canonical
exact triangle
\begin{equation}\label{lng}
\begin{CD}
HP^\hdot_Z(X) @>>> HP^\hdot(X) @>>> HP^\hdot(U) @>>>
\end{CD}
\end{equation}
Note that the scheme $Z$ actually enters into this construction
only through its open complement $U \subset X$. In particular, there
is no need to assume that $Z \subset X$ is a Poisson subscheme in
any sense.

The reader will find the precise definitions and statements on
Poisson cohomology in the Apenndix, with all the proofs. The only
statements that we will actually use in the main body of the paper
are contained in Proposition~\ref{geom},
Corollary~\ref{supp.geom}, and in Lemma~\ref{cm.i}.

The reason we are interested in Poisson cohomology is its  role
in the study of Poisson deformations of a Poisson scheme $X$. This
role is completely analogous to the role of the standard cotangent
complex $\Omega_\idot(X)$ and the groups
$\Ext^\hdot(\Omega_\idot,\calo_X)$ in the standard deformation theory
of a scheme $X$. The analogy can be pushed quite far. 
We will use Poisson deformation theory in  Subsection~\ref{def.sub}
to prove Theorem~\ref{poi.def.uni}.

\section{The case of $\dim V = 2$.}\label{dim2}

Before we proceed to the study of general symplectic quotient
singularities $V/G$, we need to consider the particular case
$\dim V = 2$. This is the case of the so-called {\em McKay
correspondence}. Starting with the paper \cite{McK}, it has been
studied excessively by many authors. We recall here some of the
results.

Let $V$ be a complex vector space of dimension $\dim V = 2$. To
every finite subgroup $G \subset Sp(V)$, one canonically associates
a simply-laced root system whose rank is equal to the number of
non-trivial conjugacy classes in $G$. Let $\g$ be the simple Lie algebra
associated to this root system. The Cartan
subalgebra $\h \subset \g$ is naturally dual, $\h \cong C^*$, to the
base $C$ of the  Calogero-Moser  deformation $\M/C$ of the quotient
variety $V/G$. The conjugacy classes of elements $g \in G$
define a basis in the vector space $C \cong \h^*$, which is in fact
a basis of simple roots. The dual space $\g^*$ is naturally a
Poisson scheme. The nilpotent cone $\N \subset \g^*$ is a Poisson
subscheme consisting of a finite number of coadjoint orbits, all of which
are symplectic. There exists a unique orbit $\N_{subreg} \subset \N$
of dimension $\dim \N_{subreg} = \dim \N - 2$, called the {\em
subregular nilpotent orbit}. Take an arbitrary element $n \in
\N_{subreg}$ and let $\p \subset \g^*$ be an affine space passing
through $n \subset \g^*$ and transversal to $\N_{subreg} \subset
\g^*$. The Poisson structure on $\g^*$ induces a Poisson structure
on the affine space $\p$. The intersection $\p \cap \N \subset \p$
is a Poisson subscheme.

It turns out that the Poisson scheme $\p \cap \N$ is naturally
isomorphic to the quotient $V/G$. Moreover, the Poisson algebra
of functions on $\p$ has a large center, so that in fact we have a
Poisson deformation $\p/S$ over a base $S$ of dimension $\dim S =
\rk \g$. The base $S$ of the deformation $\p/S$ is canonically
isomorphic to the quotient $S = \h^*/W$ of the dual Cartan algebra
$\h^*$ by the Weyl group $W$. In particular $S$ is smooth. The
subscheme $Y \cong \p \cap \N \subset \p$ is the fiber of $\p/S$
over the point $0 \in S \cong \h^*/W$. In other words, we obtain a
Poisson deformation $\p/S$ of the Poisson scheme $V/G$ over a
smooth base $S$ of dimension $\dim S =\rk \g$.

The natural $\C^*$-action on the vector space $\g^*$ by dilatations
induces the natural grading on the Poisson algebra $\C[V/G]$. The
deformation $\p/S$ is a graded deformation. The grading on $S \cong
\h^*/W$ is induced by the $\C^*$-action by dilatations on $\h^*$.

  From now on, set $Y=V/G$. It is known that
the deformation $\p/S$ is a miniversal deformation of  $Y$
in the category of affine schemes, with the Poisson structure forgotten.
In particular, the canonical map
$$
T_oS \to \Ext^1(\Omega_\idot(Y),\calo_Y)
$$
between the Zariski tangent space $T_oS$ and the group which
classifies infinitesemal deformations of the scheme $Y$ is an
isomorphism. Moreover, the total space $\p$ and the base $S$ of the
deformation $\p/S$ are smooth; therefore the scheme $Y$ is a
complete intersection, and the cotangent complex $\Omega_\idot(Y)$
only has non-trivial cohomology in degree $0$.

Using these facts, we can now compute the Poisson cohomology groups
$HP^1(Y)$ and $HP^2(Y)$.

\begin{lemma}\label{dm2}
We have $HP^1(Y) = 0$, and
$$
HP^2(Y) \cong \Ext^1(\Omega(Y),\calo_Y) \cong T_oS.
$$
Moreover, the degrees of the natural grading on the group $HP^2(Y)$
coincide with the exponents of the Weyl group $W$.
\end{lemma}

\proof{} Denote $A = \C[Y]$. By definition of the Poisson
cohomology bicomplex $DP^{\hdot,\hdot}(A)$ we have $DP^{0,0}(A)
\cong A$, $DP^{0,k}(A) = 0$ for $k \geq 1$ and
$$
DP^{1,\hdot}(A) \qis \RHom^\hdot(\Omega^1A,A).
$$
In particular, there exists a canonical map $\kappa:HP^2(A) \to
\Ext^1(\Omega^1A,A)$ induced by the natural projection
$DP^{\hdot,k}(A) \to DP^{1,k}(A)$, $k \geq 1$. On the level of
deformation theory, the map $\kappa$ corresponds to forgetting the
Poisson structure. Since the universal deformation $\p/S$ of the
scheme $Y$ does admit a Poisson structure, the map $\kappa:HP^2(A) \to
\Ext^1(\Omega^1A,A)$ is surjective. Thus to prove the Lemma, it
suffices to prove that $HP^1(A) = 0$ and that $\kappa$ is an injective
map.

Consider the spectral sequence associated to the bicomplex
$HP^{\hdot,\hdot}(A)$. We have
$$
E_1^{p,q} = \Ext^q(\Lambda^p_A\Omega_\idot(A),A),
\qquad p,q \geq 0.
$$
The only term which contributes to $HP^1(A)$ is the term
$E_\infty^{1,0}$. The only non-trivial terms which contribute to
$\HP^2(Y)$ are $E_1^{2,0}$, $E_1^{1,1}$ and $E_1^{0,2}$. The term
$E_1^{0,2} = \Ext^2(A,A)$ vanishes. The term $E_1^{1,1}$ is
precisely $\Ext^1(\Omega^1A,A)$. Moreover, the claim about the
gradings on this term follows from the identification $S \cong
\h^*/W$. To prove that the map $HP^2(Y) \to E_1^{1,1}$ is injective,
it suffices to prove that the term $E_\infty^{2,0}$ vanishes. Thus
is suffices to prove that $E_\infty^{2,0} = E_\infty^{1,0} = 0$. We
will prove that already $E_2^{p,0} = 0$ for every $p \geq 1$.

Indeed, denote by $j:U \hookrightarrow Y$ the embedding of the open
complement $U = Y \setminus \{0\} \subset Y$ to the origin $o \in
Y$. By definition we have
\begin{align*}
E_1^{p,0} &\cong \Hom(\Lambda^p\Omega_Y,\calo_Y) \cong
\Hom(\Lambda^p\Omega_Y,j_*\calo_Y)\\
& \cong \Hom_U(\Lambda^p\Omega_U,\calo_U) \cong H^0(U,\Lambda^p\T_U).
\end{align*}
The quotient
map $\pi:V \to Y = V/G$ is \'etale over $U$, so that 
$$
H^0(U,\Lambda^p\T_U) \cong H^0(\pi^{-1}(U), \Lambda^p\T_V)^G.
$$
Moreover, since $\pi^{-1}(U) \subset V$ is the complement to a point
in a smooth scheme, the right-hand side is isomorphic to
$H^0(V,\Lambda^p\T_V)$. The differential $d_1:E_1^{p,0} \to
E_1^{p+1,0}$ in the spectral sequence is induced by the Poisson
differential on the space $H^0(V,\Lambda^p\T_V)$. Hence,
$
E_2^{p,0} \cong HP^p(V)^G.
$
Since $V$ is smooth and symplectic, and $H^p(V,\C) = 0$ for $p \geq
1$, this implies that $E_2^{p,0} = 0$ for $p \geq 1$. This finishes
the proof. \endproof

In particular we see that $\p/S$ is the universal {\em Poisson}
deformation of the Poisson scheme $Y = V/G$.

Recall
 the  Calogero-Moser  deformation $\M/C$ introduced in \cite{EG}.
The  deformation $\M/C$ does not
coincide with the universal deformation $\p/S$ but there is a
Cartesian square
$$
\begin{CD}
\M @>>> \p\\
@VVV @VVV\\
C @>>> S
\end{CD}
$$
Here  $C \cong \h$, $S \cong \h/W$,
and the bottom row in the square may be identified with the canonical
projection $C =\h\onto \h/W=S$ (which 
is not an isomorphism). Thus,
Lemma~\ref{dm2} provides a canonical identification $HP^2(Y) \cong
C/W$ (this is an identification of algebraic varieties, not of
vector spaces). 

The identification $HP^2(Y) \cong C/W$ yields:
$
\dim HP^2(Y) = \dim C = \bn
$,
 the
equality claimed in Propositon~\ref{vit2}. 
To study the
higher-dimensional case we will need  the
following twisted version of this equality.

\begin{lemma}\label{twist}
Let $G \subset G' \subset Sp(V)$ be two finite groups acting on the
two-dimensional symplectic space $V$. Assume that $G \subset G'$ is
normal, and let $H=G'/G$ be the quotient group, which act naturally
on the quotient variety ${\yvg}$ and on the vector space $C$.

Then we have:
$
\dim HP^2(Y)^H = \dim C^H.
$
\end{lemma}

\proof{} The group $H$ preserves the root system $\Delta \in C$
corresponding to $G \subset Sp(V)$ and the base of simple roots
defined by the conjugacy classes in $G$. Moreover, $H$ commutes with
the action of the Weyl group $w$. We have $HP^2(Y) \cong C/W$. The
space $HP^2(Y)^H$ of $H$-invariant vectors is the subvariety
$(C/W)_H \subset C/W$ of $H$-fixed points in $C/W$. Therefore it
suffices to use the following standard result.

\begin{lemma}
Let $\Delta \in C_\R$ be a root system with Weyl group $W$, and let
$H$ be a finite group of automorphisms of the root system $\Delta$
which preserves a Weyl chamber $C^+_\R \subset C_\R$. Let $C = C_\R
\otimes_\R \C$ be the complexification of the real vector space
$C_\R$. Then the quotient map $C \to C/W$ induces a surjective
finite map
$$
C^H \to (C/W)_H
$$
onto the set of $H$-fixed points in the quotient variety
$C/W$.
\end{lemma}

\proof{} It is well-known that the set $C^+ = C^+_\R +
\sqrt{-1}C_\R$ is a fundamental domain for the $W$-action on the
vector space $C$, so that we have an isomorphism $C^+ \cong C/W$.
Since $H$ preserves $C^+$, it induces an isomorphism
$$
C^+_H \cong C^+ \cap C^H \cong (C/W)_H,
$$
which gives a section of the map $C^H \to (C/W)_H$.
\endproof

To describe the Calogero-Moser deformation $\M/C$ more functorially,
we recall that the quotient ${\yvg}$ admits a canonical smooth
symplectic resolution $X \to Y$  (coming from the Springer
resolution $\tilde{\N} \to \N$ of the nilpotent cone  in
$\g^*$). Moreover, the scheme $X$ has a universal deformation
$\X/C$, and the total space $\X$ is symplectic over $C$. Therefore
the algebra $H^0(\X,\calo_{\X})$ of global functions on $\X$ has a
natural Poisson structure. This algebra coincides with the algebra
of functions on the Calogero-Moser deformation
$\M/C$. Geometrically, we have a natural birational projective map
$\pi:\X \to \M$ and an isomorphism $\pi_*\calo_{\X} \cong
\calo_{\M}$. The cohomology group $H^2(X,\C)$ of the resolution $X$
is naturally identified with the base $C$ of the deformation
$\X/C$. The map $\pi$ is compatible with the projections to the base
$C$. Moreover, it is an isomorphism over a generic point $c \in C$.

These facts taken together immediately imply Propositon~\ref{vit2},
Theorem~\ref{zcsmooth} and Theorem~\ref{res} in the case $\dim V =
2$.

\section{The computation of $HP^2(V/G)$.}\label{hp2.sec}

Let $V$ be an arbitrary finite-dimensional symplectic vector space,
and let $G \subset Sp(V)$ be a finite subgroup. In this section we will
compute the Poisson cohomology groups $HP^1(V/G)$ and $HP^2(V/G)$ of
the quotient $Y = V/G$. In particular, we will prove
Proposition~\ref{vit2}.

First, we introduce some notation. Notice that the vector space $V$
is naturally stratified by subspaces $V^H \subset V$ of
$H$-invariant vectors for various subgroups $H \subset G$. All these
subspaces are symplectic, hence even-dimensional. This induces a
stratification of the quotient variety $Y = V/G$.
The strata of the stratification are known to be 
smooth  symplectic locally-closed subvarieties in $V/G$,
in particular have even dimension.
Moreover, these strata turn out to be exactly the
{\em symplectic leaves} of the standard Poisson structure on $V/G$,
cf. e.g.~\cite{BG}.

In more detail, let $\Gamma\subset G$ be a subgroup
which is the isotropy group of an element of $V$.
Then $V^\Gamma\subset V$ is a nonzero vector subspace such that the
symplectic form on $V$ restricts to a nondegenerate 2-form
on $V^\Gamma$. Further, let $U^\Gamma\subset V$
be the set of  points of $V$ whose stabilizer
(in $G$) is equal to $\Gamma$. It is clear that
 $U^\Gamma \subset V^\Gamma$, moreover, it is
known from the theory of finite
group actions that $U^\Gamma$ is a non-empty Zariski open, hence dense,
subset in $V^\Gamma$. The strata of the stratification
of $V/G$ that has been mentioned in the previous paragraph
are defined to be 
the images  under the projection
$V\onto V/G$ of the   sets of the form $U^\Gamma$, as $\Gamma$
varies inside $G$.

It is straightforward to verify that
$N(\Gamma)$, the normalizer of $\Gamma$ in $G$,
preserves the set $U^\Gamma$. Furthermore, the
resulting $N(\Gamma)/\Gamma$-action on $U^\Gamma$ is free,
and the projection $V\onto V/G$ induces an isomorphism
of $U^\Gamma/(N(\Gamma)/\Gamma)$ with its image in $V/G$,
that is, with the corresponding stratum
of the  stratification.

Now, recall the notation $Y=V/G$ and write
 $U \subset Y$ for the complement to the union of all the
strata in $Y$ of codimension $\geq 4$. 
Thus we have $U=U_0\,\coprod \,
(\cup_{i\geq 1}\,U_i )$,
where $U_0$ is the unique open stratum in $Y$,
and $U_i\,,\,i=1,2,\ldots,$
are all the  codemension two strata 
of the  stratification. 
By the earlier discussion, each stratum is the image
 under the projection
$V\onto V/G$
of the set $U^{G_i}\subset V$
where $G_i \subset G$ is a subgroup
such that the fixed point set $V_i := V^{G_i} $
is a codimension 2
vector subspace in $V$ and, moreover, such that
the set $U^{G_i}=\{v\in V\;|\;\text{\em isotropy group of}\; v = G_i\}\,$
is a Zariski dense subset in $V^{G_i} $.
Two  subgroups as above give rise to  the same stratum $U_i$
 if and only if they  are conjugate within $G$.

Let $G'_i :=N(G_i)\subset G$ denote the normalizer of the subgroup $G_i$,
 and  $H_i := G_i'/G_i$ the quotient group. 
As we have explained above, the group
$G'_i$ preserves the set $U^{G_i}$,
the induced $H_i$-action   on $U^{G_i}$ is free;
furthermore, we have an isomorphism
$U^{G_i}/H_i\iso U_i$.

For each $i \geq 1$, let $W_i$ denote the annihilator of
$V_i$ with respect to the symplectic form.
Thus, $\dim W_i=2$, and there is a canonical $G'_i$-stable direct sum
decomposition $V=V_i\bigoplus W_i$.
Let $Y_i =
W_i/G_i$ be the quotient variety. The group $H_i$ acts on the
variety $Y_i$ and preserves the origin  $o \in Y_i$. The
 direct sum
decomposition $V \cong V_i \oplus W_i$ induces a map
$\eta_i:(V_i \times Y_i)/H_i \to Y$. 
We have a Zariski  open subset $U'_i\subset V_i$ such that
 $\eta_i$ 
maps $(U'_i\times \{o\})/H_i$ isomorphically
onto $U_i$ and, moreover, 
is \'etale in
a Zariski  open neighborhood of $U_i = (U'_i/H_i) \times \{o\}
\subset (V_i \times Y_i)/H_i$. We note for further record
that $V_i\smallsetminus U'_i$ is a finite union
of vector subspaces in $V_i$ of
(complex) codimension $\geq 2$. It follows that
$H^l(U'_i, \C)=0,$ for $l=1,2$.

As in Subsection~\ref{comp.sub}, let $C$ be the space of all
$G$-invariant functions on the set $\Sigma$ of symplectic
reflections in $G$. By definition we have $\dim C = \bn$. For every
$i$, denote by $C_i$ the vector space of $G_i$-invariant $\C$-valued
functions on the group $G_i$, and let $B_i \subset C_i$ be the
subspace of $G_i'$-invariant functions. The group $H_i$ acts on the
space $C_i$ and we have $B_i = C_i^{H_i} \subset C_i$.  Every
element $g \in G_i\smallsetminus\{1\}$
pointwise fixes the codimension two subspace $V_i$, hence is a symplectic
reflection.
Therefore we have a
natural restriction map $C \to C_i$. This map factors through a map
$C \to B_i = C_i^{H_i} \subset C_i$. The map $C \to B_i$ is
surjective. Indeed, it suffices to check that two elements $g_1,g_2
\subset G_i$ conjugate in $G$ are already conjugate in $G'_i$; and
if $g \cdot g_1 \cdot g^{-1} = g_2$, then $g$ preserves the subspace
$V_i = V^{g_1} \subset V$, hence lies in $G'_i$. Moreover, every
symplectic reflection $g \in G$ lies in one of the subgroups $G_i
\subset G$, -- namely, the stabilizer of the invariant subspace $V^g
\subset V$.  Therefore the projections $C \to B_i$ give a natural
splitting
\begin{equation}\label{splt}
C \cong \bigoplus\nolimits_i\,B_i = \bigoplus\nolimits_i\, C_i^{H_i}.
\end{equation}
We recall that the quotient singularity $Y = V/G$ carries a natural
Poisson structure, so that we have the Poisson cohomology groups
$HP^\hdot(Y)$. Moreover, $\C[Y]$ is a positively-graded Poisson
algebra of degree $2$, and this grading induces a grading on
$HP^\hdot(Y)$.

\proof[Proof of Proposition~\ref{vit2}.] The scheme $Y = V/G$ is
normal, and the smooth locus $U_0 \subset Y$ carries a
non-degenerate symplectic form.  Moreover, by \cite{Wat} the scheme
$Y$ is Gorenstein, hence Cohen-Macaulay. By
definition, the complement to the open subset $U \subset Y$ is of
codimension $\geq 4$. Therefore by Lemma~\ref{cm.i}, the natural map
$
HP^k(Y) \to HP^i(U)
$
is an isomorphism for $k = 1,2$.

Since the open subscheme $U_0 \subset U$ is smooth and symplectic,
we have $HP^\hdot(U_0) \cong H^\hdot(U_0,\C)$. But $U_0 = (V
\setminus Z)/G$ is the quotient of the complement in the vector
space $V$ to some closed subscheme $Z \subset V$ of $\codim \geq 2$.
Therefore the singular cohomology $H^\hdot(U,\C) = H^\hdot(V
\setminus Z)^G$ vanishes in low degrees,
$$
HP^k(U_0) \cong H^k(U_0,\C) \cong H^k(V \setminus Z,\C)^G \cong
H^k(V,\C)^G = 0, \qquad k=1,2,3.
$$
We conclude that for $k = 1,2$ the natural map $HP^k_{U \setminus
U_0}(U) \to HP^k(U)$ is an isomorphism, and the isomorphism $HP^k(Y)
\cong HP^k(U)$ factors through a canonical isomorphism
$$
HP^k(Y) \cong HP^k_{U \setminus U_0}(U).
$$
The complement $U \setminus U_0$ is the disjoint union of the closed
strata $U_i \subset U$, $i \geq 1$, and for every $i \geq 1$, we
have a map $\eta_i:(V_i \times Y_i)/H_i \to Y$. The map $\eta_i$
identifies $(U'_i/H_i) \subset (V_i \times Y_i)/H_i$ with the
stratum $U_i \subset Y$, and it is \'etale in an open neighborhood
of $U_i \cong (U_i'/H_i) \subset (V_i \times Y_i)/H_i$. By
Corollary~\ref{supp.geom}\thetag{i}, we have
$$
HP^\hdot_{U_i}(U) \cong HP^\hdot_{U_i}((U'_i \times Y_i)/H_i).
$$
By Proposition~\ref{geom}\thetag{iii}, the right-hand hand side
is naturally isomorphic to
$$
HP^\hdot_{U_i}((U'_i \times Y_i)/H_i) \cong
\left(HP^\hdot_{U'_i}(U'_i \times Y_i)\right)^{H_i}.
$$
Hence, 
we obtain a canonical direct sum decomposition
$$
HP^k(Y) \cong HP^k_{U \setminus U_0}(U) \cong \bigoplus\nolimits_i\,
\left(HP^k_{U'_i}(U'_i \times Y_i)\right)^{H_i}.
$$
The Kunneth
formula yields
$$HP^k_{U'_i}(U'_i \times Y_i)\cong
\bigoplus_{0\leq l\leq k}HP^{l}(U'_i)\otimes HP^{k-l}_{o}(Y_i).
$$
Further, recall that by construction  $U'_i$ is an open
subset in $V_i$ whose complement has complex codimension $\geq 2$.
In particular,  $U'_i$ is  smooth symplectic and
connected. Therefore, $HP^{l}(U'_i)=H^l(U'_i, \C)$;
moreover, this group is 1-dimensional (and has trivial
 action of the
group $H_i$) if $l=0$, and
is equal to zero if  $l=1,2$.
We conclude that the
product decomposition \eqref{kun.supp} reduces to an $H_i$-equivariant
isomorphism
$$
HP^k_{U'_i}(U'_i \times Y_i) \cong HP^k_o(Y_i) \otimes HP^0(U'_i)
\cong HP^k_o(Y_i),
\qquad k \leq 2.
$$
Therefore for $k = 1,2$ we have
$$
HP^k(Y) \cong \bigoplus\nolimits_i\, HP^k_o(Y_i)^{H_i}.
$$
But for every $i \geq 0$, the complement $Y_i \setminus \{o\} \cong
(W_i \setminus \{0\})/G_i$ is smooth, symplectic, and satisfies
$H^k(Y_i \setminus \{0\}) \cong H^k(W_i \setminus \{0\})^{G_i} =0$,
$k = 1,2,3$. Therefore $HP^k(Y_i \setminus \{o\}) = 0$ for $k =
1,2,3$, and we have an isomorphism $HP^k_o(Y_i) \cong
HP^k(Y_i)$.

Collecting together all of the above, we conclude that for $k = 1,2$
there exists a natural isomorphism
$$
HP^k(Y) \cong \bigoplus\nolimits_i\, \left(HP^k(Y_i)\right)^{H_i}.
$$
But we have already computed $HP^k(Y_i)$, $k = 1,2$ in
Section~\ref{dim2}. Lemma~\ref{dm2} immediately implies that
$HP^1(Y) = 0$, and Lemma~\ref{twist} shows that
$$
\dim HP^2(Y) = \sum_i \dim HP^2(Y_i)^{H_i} = \sum_i \dim C_i^{H_i}.
$$
By \eqref{splt} the right-hand side is equal to $\dim C$. Moreover,
the natural grading on each of the $HP^2(Y_i)$ has positive degrees.
\endproof

This Proposition of course applies to any pair $\langle V, G
\rangle$. In particular, one can take the trivial group $G=\{e\}$,
and obtain, as expected, the equality $HP^2(V)=0$. Taking the
quotient by some non-trivial $G \subset Sp(V)$ increases the second
Poisson cohomology group and creates new deformations. These
deformations are indeed new: they do not lift to $G$-equivariant
deformations of the symplectic vector space $V$. We do not formulate
this fact precisely because we will not need it. However, in
Section~\ref{res.sec} we will need the following claim.

\begin{lemma}\label{covers.trivial}
Let $\Y/S$ be a Poisson deformation of a symplectic quotient
singularity ${\yvg}$ over a local Artin base $S$. Assume that there
exists a Poisson deformation $\V/S$ of the vector space $V$ and a
map $\eta:\V/S \to \Y/S$ which extends the quotient map $\eta:V \to
Y$. Then the Poisson cocycle $\Theta_Y \in HP^2(\Y/S)$ is trivial.
\end{lemma}

\proof{} Assume that $\Theta_Y \neq 0$. Let $\m \subset \C[S]$ be
the maximal ideal, and let $k$ be the largest integer such that
$\Theta_Y =0 \mod \m^k$. Then $\Theta_Y  \mod \m^{k+1}$ is a
non-trivial element in the Poisson cohomology group
$$
HP^2(\Y/S,\m^{k}/\m^{k+1}) \cong H^2(Y) \otimes_\C
\left(\m^k/\m^{k+1}\right).
$$
But the group $HP^2(Y)$ admits a natural map $\eta:H^2(Y) \to
H^2(Y,\eta_*\calo_V)$ which identifies it with the direct summand
$H^2(Y,\eta_*\calo_V)^G \subset H^2(Y,\eta_*\calo_V)$. Therefore
$\eta(\Theta_Y)$ must be a non-trivial class in the group
$$
H^2(Y,\eta_*\calo_V) \otimes_\C \left(\m^k/\m^{k+1}\right).
$$
This is impossible. Indeed, by the functoriality of the Poisson
cohomology, the class $\eta(\Theta_Y)$ comes from the class
$\Theta_V \in HP^2(V)$ of the Poisson cocycle on the vector space
$V$, and the group $HP^2(V)$ is a trivial group.
\endproof

In truth, under the assumptions of the Lemma the whole deformation
$\Y/S$ must be trivial, not only its Poisson cocycle $\Theta_Y$. But
the proof of this fact is slightly harder. We settle for the weaker
version to save space.

One final remark is the following: since the Poisson algebra
$\C[Y]$ is positively graded, Lemma~\ref{compl.eq} immediately
shows that all the results of this Section are valid not only for
$Y$, but also for its completion $\wh{Y}$ at the origin  $o
\in Y$.

\section{Resolutions.}\label{res.sec}
In this section we will prove
 Theorem~\ref{vit1} and
Theorem~\ref{res}. 

Below, we use the notation $H^\hdot(M)=
H^\hdot(M,\C)$ for the singular cohomology of a variety $M$ with
{\em complex} coefficients. The homology of $M$ is usually taken
with  {\em rational} coefficients, i.e., $H_\idot(M)=
H_\idot(M,\Q).$
\subsection{Geometry of the resolution.}
We will need several general facts on the geometry of  symplectic
resolutions. Most of them are well-known. The reader can find the
proofs, for instance, in the papers \cite{K1}, \cite{K3}.

Let $X$ be an irreducible smooth variety over $\C$ equipped with a
closed nowhere-degenerate $2$-form $\om$, and let $Y = X\aff$ be
its affinization. Assume that the natural map $\pi:X \to Y$ is
projective and birational. By definition $Y$ is normal, and we have
$\pi_*\calo_X \cong \calo_Y$. The canonical bundle $K_X$ of the
manifold $X$ is trivial (it is trivialized by the top power of the
symplectic form $\om$). Therefore map $\pi:X \to Y$ is one-to-one
over the smooth locus $Y_o \subset Y$. By the
Grauert-Riemenschneider Vanishing Theorem, we have
\begin{equation}\label{vanish}
H^i(X,\calo_X) = H^0(Y, R^i\pi_*\calo_X) = H^0(Y,R^i\pi_*K_X) = 0
\end{equation}
for $i \geq 1$. Considering the exponential exact sequence on $X$,
one easily deduces that $R^1\pi_*\Z_X = 0$, where $\Z_X$ is the
constant sheaf on  $X$. 

The symplectic form $\om$ on the smooth variety $X$ induces a
Poisson structure on the algebra
$
\C[Y] \cong H^0(X,\calo_X).
$
Further, the resolution $X \to Y$ is known, see e.g. \cite{K},
 to be semismall (that is, $\dim X \times_Y X =
\dim X$). This implies that the Leray spectral sequence which
computes the group $H^2(X,\Q)$ degenerates, and we have a canonical
short exact sequence
of rational  cohomology groups
\begin{equation}\label{hyx}
\begin{CD}
0 @>>> H^2(Y,\Q) @>>> H^2(X,\Q) @>>> H^0(Y,R^2\pi_*\Q) @>>> 0
\end{CD}
\end{equation}
Set
$$
H_2(X/Y) := \Ker(H_2(X,\Q) \to H_2(Y,\Q))
$$
the subgroup in the homology group $H_2(X,\Q)$ dual to
$H^0(Y,R^2\pi_*\Q)$. The following fact will be very important.

\begin{lemma}\label{omega.lifts}
Let $[\om] \in H^2(X,\C)$ be the cohomology class of the
symplectic form $\om$. Then for every homology class $\alpha \in
H_2(X/Y)$, we have $\langle \alpha, [\om] \rangle = 0$, and
\begin{equation}\label{lft}
[\om] = \pi^*[\om]_Y,
\end{equation}
for some cohomology class $[\om]_Y \in H^2(Y,\C)$. Moreover, the
class $[\om]_Y$ depends only on the Poisson structure on the
variety $Y$, not on the resolution $X$.
\end{lemma}

\proof{} This is not new, see e.g. \cite{CF}, \cite[Lemma 2.2]{K3},
\cite{WW}. We only give a sketch of the proof. By \eqref{hyx}, to
establish \eqref{lft} it suffices to prove that $[\om]$ vanishes
on all the fibers of the map $Y \to X$. Taking a resolution of
singularities, it suffices to prove that $f^*[\om] = 0$ for every
map $f:Z \to X$ from a smooth projective manifold $Z$. By Hodge
theory, this is equivalent to proving that
$$
f^*[\overline{\om}] = 0 \in H^{0,2}(Z) = H^2(Z,\calo_Z),
$$ 
where $[\overline{\om}]$ is the complex-conjugate cohomology
class. But already $[\overline{\om}] = 0$, since $H^2(X,\calo_X)
= 0$ by \eqref{vanish}. To prove that the class $[\om]_Y$ is
canonical, note that its restriction to the non-singular part $U
\subset Y$ is represented by an explicit $2$-form $\om_Y$ which
is inverse to the non-degenerate Poisson bivector $\Theta$. Thus
$[\om]_Y|_U$ depends only on $\Theta$. Moreover, we know that the
form $\om_Y$ extends to some, hence to an arbitrary smooth
resolution $X \to Y$ and represents a cohomology class $[\om]_Y
\subset H^2(X,\C)$ coming from $H^2(Y,\C)$. 
The  class $[\om]_Y$ is
completely determined by the form $\om_Y$ since the map
$H^2(Y,\C) \to H^2(X,\C)$ is injective.  \endproof

In order to study deformations of the pair $\langle X,\om 
\rangle$, we follow the paper \cite{KV} and introduce the period
map. To do this, note that for every smooth deformation $\X/B$ over a
local Artin base $B$, the Gauss-Manin connection trivializes the
relative de Rham cohomology $H^2_{DR}(\X/B)$,
\begin{equation}\label{coho.constant}
H^2_{DR}(\X/B) \cong H^2(X,\C) \otimes_\C \C[B].
\end{equation}
Therefore the class $[\om]$ can be canonically considered as a
class in the group $H^2(X,\C) \otimes_\C \C[B]$. This induces a
map 
$$
P:B \to H^2(X,\C)
$$
called the {\em period map}. The fundamental theorem of \cite{KV}
claims that the period map completely determines the
deformation. More precisely, for every local Artin scheme $B$ and a
map $P:B \to H^2(X,\C)$ which sends the closed point $o \in B$ to
the class $[\om] \in H^2(X,\C)$, there exists a deformation
$\X/B$ of the pair $\langle X,\om \rangle$ with the period map
$P$. Moreover, such a deformation $\X/B$ is unique up to a
non-canonical isomorphism.

Set $B := H^2(X,\C)$, and let $\wh{B}$
denote the formal neighborhood of $[\om]
\in H^2(X,\C)$. This way, one obtains a {\em universal} deformation
${\wh\X}/\wh{B}$ of the pair $\langle X, \om \rangle$. Since one
has to pass to the limit, the universal deformation ${\wh\X}/\wh{B}$
is only a formal scheme.

\subsection{Globalizing the deformation.}\label{glb.sub}
We can now start proving Theorem~\ref{vit1}. Our method will be to
extend the formal deformation from \cite{KV} to an actual
deformation defined over a global base.

Throughout the subsection, let $X$ be a convex symplectic manifold with
symplectic form $\om$. We write $B := H^2(X,\C)$,
and equip the vector space $B$ with
$\C^*$-action so that $z\in\C^*$ acts via multiplication by
$z^{-l}$. Further, let $\wh{B}$
denote the formal neighborhood of $[\om]
\in B$, and
 let  $\m \subset \C[\wh{B}]$ be the
maximal ideal of the complete  algebra
$\C[\wh{B}]$. 
Write
$\wh{\X}/\wh{B}$ for the universal formal deformation of the pair
$\langle X,\om \rangle$.

Put $Y := X\aff$. Assume that
the Poisson algebra $\C[Y]=H^0(X,\calo_X)$
is a finitely-generated positively graded Poisson algebra of
degree $l > 0$, and consider the corresponding $\C^*$-action on the
scheme $Y$.  We split the proof of Theorem~\ref{vit1} into three
Propositions below.

\begin{prop}\label{omega.0}
The $\C^*$-action on $Y$ lifts uniquely to a $\C^*$-action on $X/Y$,
and the resulting action on $X$ extends to a $\C^*$-action on
$\wh{\X}/\wh{B}$. Moreover, we have $[\om] =0$ in $H^2(X)$.
\end{prop}
This result is known, and we refer to
\cite{Fu}, \cite{K1}, and \cite{V} for the proofs.
We also have the following  infinitesimal version of the  above proposition.

\begin{lemma}\label{vf.extends}
Let $Y$ be an irreducible algebraic variety, and let $X\to Y$ be
smooth projective semismall resolution of  $Y$. If the
canonical bundle  is trivial
then every vector field $\xi$ (a derivation of  ${\mathcal{O}}_Y$) on
$Y$ lifts
canonically to
a vector field on  $X$.
\end{lemma}

\proof{} (Compare \cite{CF}, where a similar statement is proved for
isolated singularities, but without the semismallness asumption.)
The vector field $\xi$ canonically lifts to a vector field defined
outside of the exceptional locus of the map $\pi:X \to Y$.  Bince
the variety $X$ is smooth, the sheaf $\T(X)$ of vector fields is
reflexive. Therefore it suffices to prove that $\xi$ to the generic
point of an arbitrary exceptional Weil divisor $E \subset X$. Since
$\pi:X \to Y$ is semismall, the image $\pi(E) \subset Y$ of the
divisor $E$ is a subvariety of codimension $2$.  Therefore near the
generic point of the subvariety $\pi(E) \subset Y$, the variety $Y$
is the product of a smooth variety and an isolated surface
singularity, and it suffices to prove the Lemma in the case $\dim Y =
2$. In this case, the triviality of the canonical bundle implies that $Y$ is a Du
Val point. Then Lemma easily follows from an explicit construction
of the minimal resolution as a transversal slice to a subregular
nilpotent orbit, see Section~\ref{dim2}.
\endproof

\begin{prop}
The formal scheme $\wh{\X}$ extends to an actual scheme $\X/B$
defined over the whole affine space $B=H^2(X)$. The map $\sigma:\X
\to H^2(X)$ is smooth and $\C^*$-equivariant, and the scheme $\X$ is
symplectic over $B$. Moreover, the relative cohomology sheaves
$R^\hdot\sigma_*\Q$ are constant sheaves on $H^2(X)$, and the
canonical base change morphism
$$
H^k(\X_b,\Q) \to \left(R^k\pi_*\Q\right)_b
$$
is an isomorphism for every point $b \in B$.
\end{prop}

\proof{} To extend the formal scheme $\wh{\X}/\wh{B}$ to a scheme
over the whole $B$, we repeat the argument of \cite[Lemma 4.2]{K3}.
Take an ample line bundle $L$ on $X$ and note that, since by
\eqref{vanish} $H^1(X,\calo_X) = H^2(X,\calo_X) = 0$, the bundle $L$
canonically extends to the formal scheme $\wh{\X}/\wh{B}$. Moreover,
using the same cohomology vanishing it has been shown
in \cite{K3} that $H^0(\wh{\X},\calo_{\wh{\X}})$,
the algebra of  global functions on $\wh{\X}$,
is a filtered algebra whose associated graded algebra
is isomorphic to $H^0(X,\calo_X)\otimes\C[t]$,
a Noetherian algebra. It follows that
$H^0(\wh{\X},\calo_{\wh{\X}})$ is  a
Noetherian algebra, which is also complete with respect to the $\m$-adic topology.
Thus,
$$ 
\wh{\Y} = \Spf(H^0(\wh{\X},\calo_{\wh{\X}})),
$$ 
the formal spectrum,
is an affine Noetherian formal scheme flat over $\wh{B}$, and
$\wh{\X}$ is projective over $\wh{\Y}$, with an ample line bundle
$L$. Therefore we can apply the Grothendieck {\it algebraization
theorem} \cite[III, Th\'eor\`eme 5.4.5]{EGA} to $\wh{\X}/\wh{\Y}$ and
conclude that the formal scheme $\wh{\X}$ extends to an actual
scheme $\X$ over $\wh{\Y}$ (and {\em a posteriori}, over
$\wh{B}$). Since the algebra $\C[B]$ of functions on the vector
space $B=H^2(X)$ is positively graded, we are done by
Lemma~\ref{compl.eq}. Moreover, by construction the scheme $\X/B$ is
symplectic and smooth, and the map $\X \to B$ is $\C^*$-equivariant.

Let us now analyze the relative de Rham cohomology sheaves
$H^\hdot_{DR}(\wh{\X}/\wh{B})$. For every $k \geq 0$, denote by $B_k
= \Spec\calo_{\wh{B}}/\m^k \subset \wh{B}$ the $k$-th infinitesimal
neighborhood of the special point in the local scheme $\wh{B}$, and
let $\X_k = \wh{\X} \times_{\wh{B}} B_k$. By \cite[III, Th\'eor\`eme
4.1.5]{EGA}, for every $p,q \geq 0$ we have canonical isomorphisms
\begin{equation}\label{limit}
H^q(\wh{\X},\Omega^p(\wh{\X}/\wh{B})) \cong \lim_{\leftarrow}
H^q(\X_k,\Omega^p(\X_k/B_k)),
\end{equation}
and the projective system on the right-hand side satisfies the
Mittag-Leffler condition \cite[0, 13.1]{EGA}. Moreover, for every
$k,n$ the de Rham cohomology module $H_{DR}^n(\X_k/B_k)$ is a
finitely-generated module over $\C[B_k]$, an Artinian algebra.
Therefore for every $n \geq 0$, the projective system
$H_{DR}^n(\X_k/B_k)$ also satisfies the Mittag-Leffler condition, and
by \cite[0, Proposition 13.2.3]{EGA} the maps \eqref{limit} induce
canonical isomorphisms
$
H_{DR}^n(\wh{\X}/\wh{B}) \cong \lim_{\leftarrow} H_{DR}^n(\X_k/B_k).
$ 

We conclude that for every $n \geq 0$, the $\C[\wh{B}]$-module
$H_{DR}^n(\wh{\X}/\wh{B})$ is finitely generated, complete and
separated with respect to the $\m$-adic topology. Since it carries a
flat connection -- namely, the Gauss-Manin connection -- it must be
a free $\C[\wh{B}]$-module.

By Lemma~\ref{compl.eq}, this implies that the relative de Rham
cohomology sheaves $H^\hdot_{DR}(\X/B)$ are also free sheaves of
$\C[B]$-modules. Therefore the $D$-modules
$R^\hdot\sigma_*\calo_{\X}$ are constant $D$-modules on $B$. This
means that the sheaves $R^\hdot\sigma_*\Q$ are constant sheaves.

Finally, let $b \in B$ be an arbitrary point with embedding $i_b:b
\to B$, and let $\X_b \subset \X$ be the fiber of $\X$ over $b$,
with embedding $\X_b \hookrightarrow \X$ denoted by the same letter
$i_b$. Since the scheme $\X/B$ is smooth, we can apply Poincare
duality together with the Proper Base Change Theorem and obtain a
canonical isomorphism
$$
H^\hdot(\X_b,\Q) \cong H^\hdot(X,i_b^!\Q)[\dim B] \cong
i_b^!R^\hdot\sigma_*\Q[\dim B].
$$
Since the sheaves $R^k\sigma_*\Q$ are locally constant on $B$, the
right-hand side is isomorphic to $(R^\hdot\sigma_*\Q)_b$.
\endproof

Below, by  `general enough' or `generic' point of $B$ we mean
a point from an appropriately chosen Zariski open dense subset of $B$.

To finish the proof of Theorem~\ref{vit1}, it remains to prove the
following. 

\begin{prop}\label{vit_Z}
The map $\pi:\X \to \Y$ is semismall. Moreover, for a  general enough point
$b \in B$, the induced map $\X_b \to \Y_b$ of  fibers over $b$ in the
varieties $\X/B$, resp., $\Y/B$, is an isomorphism.
\end{prop}

\proof{} (Compare \cite[Proposition 4.6]{K3}.) To prove that $\pi:\X
\to \Y$ is semismall, it suffices to prove that for every closed
subvariety $\ZZ \subset X$, the dimension of the generic fiber of
the map $\ZZ \to \pi(\ZZ)$ does not exceed the codimension $\codim
(\ZZ,\X)$. Let $Z = X \cap \ZZ$. Since $X \subset \X$ is the fiber
of a smooth map $\X \to H^2(X)$, we have $\codim (Z,X) \leq \codim
(\ZZ,\X)$. Since the map $\pi:X \to Y$ is semismall, the dimension
of the generic fiber of the map $\pi:Z \to \pi(Z)$ does not exceed
the codimension $\codim (Z,X)$. Therefore the same is true for the
generic fiber of the map $\pi:\ZZ \to \pi(\ZZ)$.

Let $b\in B$. Then,
 since $Y$ is normal, the
variety $\Y_b$ is also normal. By construction, the map $\X_b \to
\Y_b$ is projective and generically one-to-one. 
The set of points  $b\in B$ such that the map $\X_b \to
\Y_b$ has fibers
of dimension $>1$ is clearly a closed subset $B^{\tt{BAD}}\subset B$.
Hence, to complete the proof of the Proposition  it suffices to prove
that  $B^{\tt{BAD}}$ is a proper subset  in $B$.

To this end, we  choose $b\in B$ such that
the cohomology class $[\om]=b \in H^2(X) =
H^2(X,\C)$ does not annihilate any
{\it rational} homology class $[Z] \in
H_2(X/Y) \subset H_2(X,\Q)$ (the latter obviously form a countable
subset $H_2(X,\Q)\subset H_2(X,\C)$).
We claim that $b\not\in B^{\tt{BAD}}$, hence, $B^{\tt{BAD}}
\neq B$. Indeed, if not, then there is a fiber of
 the map $\X_b \to
\Y_b$ that contains a projective  curve $Z
\subset \X_b$.
The corresponding  homology class $[Z] \in
H_2(X/Y) \subset H_2(X,\Q)$ annihilates $\om$, i.e., we have
$
\left\langle [Z],[\om] \right\rangle = 0,
$ by Lemma~\ref{omega.lifts}.
But that would contradict the choice of $b$ which 
was made so that $[\om]$ does not  annihilate any
{\it rational} homology class $[Z] \in
H_2(X/Y)$. This completes the proof.
\endproof

\begin{remark} We included 
the claim that $\X \to \Y$ semismall 
 because it was used in the proof of  Lemma~\ref{vf.extends}
(if applied to the map  $\X \to \Y$); this case of the Lemma 
 will be used at one point of the argument  later, in
Subsection~\ref{on.c.i.sub}.
\end{remark}

We have completed the proof of Theorem~\ref{vit1}.

\proof[Proof of Proposition~\ref{cmp}.] Assume that the graded
Poisson algebra $\C[Y]$ satisfies the assumptions of
Theorem~\ref{poi.def.uni}, so that we have the universal graded
Poisson deformation $\Y_S/S$. Let $\tau:B \to S$ be the
classifying map of the deformation $\Y/B = \X\aff/B$. We have to
prove that the map $\tau$ is a finite map onto an irreducible
component of the variety $S$.

Let $b \in B=H^2(X)$ be a point sufficiently generic so that the map
$\X_b \to \Y_b$ is an isomorphism. Let $S' \subset S$ be the
irreducible component which contains the image $\tau(c)$ of the
point $b \in B$.

Since $\X_b \to \Y_b$ is one-to-one, the Poisson scheme $\Y_b$ is in
fact symplectic, and the Poisson deformation theory for $\Y_b$
reduces to the deformation theory of the pair $\langle \X_b =
\Y_b,\om \rangle$. Since the deformations $\X/B$, $\Y/S$ are
universal, this implies that the map $\tau:B \to S'$ is \'etale at
the point $b$.

We see that the map $\tau:B \to S'$ is generically \'etale.  To
prove that it is in fact finite, it suffices to prove that all of
its geometric fibers are finite. Moreover, since the map $\tau$ is
$\C^*$-equivariant, it suffices to consider the fiber over the
origin point $o \in S'$. We claim that this fiber in fact consists
of the point $0 \in H^2(X)$, and this follows from
Proposition~\ref{omega.0}. Indeed, for every point $b \in H^2(X)$
with $\tau(b) = o$, we apply Proposition~\ref{omega.0} to $\X_b/Y$
and deduce $[\om] = 0 \in H^2(\X_b,\C)$. This implies $b = 0$ by
the definition of the period map.
\endproof

\subsection{Comparison with the  Calogero-Moser
deformation.}\label{res.sub} 
We now turn to the proof of Theorem~\ref{res}. Let $Y = V/G$ be a
symplectic quotient singularity with a resolution $X$, let $\X/B$
be the deformation provided by Theorem~\ref{vit1}, and let $\M/C$ be
the Calogero-Moser deformation. Consider also the universal
deformation $\Y/S$ of the Poisson scheme $Y$.

By construction, we have a multi-valued map $\kappa:C \ratto B$. As
we have noted after stating Conjecture~\ref{1-1}, the domain and the
range of this map are in fact canonically isomorphic as vector
spaces.  If $\dim V = 2$, then the map $\kappa$ is single-valued and
provides an isomorphism between $C$ and $B$.

It is natural to conjecture that the same is true in the general
case: $\kappa:C \ratto B$ is single-valued and an isomorphism.  It
would immediately imply Theorem~\ref{res}. Unfortunately, we were
not able to prove it, and we have to settle for less. Namely, we
consider the direct sum decomposition $
C = \bigoplus\nolimits_i\, B_i,$ see \eqref{splt},
and the associated decomposition of the vector space $B = H^2(X,\C)
\cong C$. As the main technical step in the proof of
Theorem~\ref{res}, we show that the map $\kappa:C \ratto B$ is
compatible with these decompositions.

\begin{prop}\label{on.c.i}
Let $c \in B_i$ be a point in the subspace $B_i \subset C$ of the
base of the  Calogero-Moser  deformation $\M/C$. Let $\kappa:C \to S$
be the classifying map of the deformation $\M/C$, and let $b \in
B$ be an arbitrary point lying over $\kappa(c) \in S$.
\begin{enumerate}
\item The cohomology class $[\om]_{b} \subset B \cong
H^2(X,\C)$ of the symplectic form $\om \in \Omega^2(\X_{b})$
lies in the subspace $B_i \subset B$.
\item Assume that the point $c \in B_i$ is generic. Then a rational
homology class $\alpha \in H_2(X,\Q)$ satisfies $\langle
\alpha,[\om]_{b} \rangle = 0$ if and only if $\alpha$ is
orthogonal to the whole subspace $B_i \subset B$.
\end{enumerate}
\end{prop}

By definition of the period map, Proposition~\ref{on.c.i}\thetag{i}
means that the point $b \in B$ lies in
$B_i$. Proposition~\ref{on.c.i}\thetag{ii} means that for a generic
$c \in B_i$, the point $b \in B_i$ is generic.

The proof of Proposition~\ref{on.c.i} is rather technical. We will
give it in the next subsection. But first, we will deduce
Theorem~\ref{res} from the proposition.

\proof[Proof of Theorem~\ref{res}.] Let  $C^{\tt{GOOD}}\subset C$
be the set
of points  such that the claim of the Theorem holds for $c$.
This is
clearly a Zariski open subset in $C$. Hence, we must only show that 
$C^{\tt{GOOD}}$ is nonempty.

To this end, for every rational homology class $0\neq \phi
\in H_2(X,\Q)$, denote by $C_\phi \subset C$ the subset of elements
$c \in C$ such that
\begin{itemize}
\item For some point $b \in B$ lying over $\kappa(c) \in S$, the
class $\phi \in H_2(X,\Q) \cong H_2(\X_{b},\Q)=H_2(\X_{b})$ lies in the subgroup
$H_2(\X_{b}/\Y_{\kappa(c)}) \subset H_2(\X_{b},\Q)$.
\end{itemize}
All the subsets $C_\phi \subset C$ are closed algebraic
subvarieties.

Assume first that there exists $\phi \in H_2(X,\Q)\,,\, \phi\neq 0,$ such
that $C_\phi = C$. Then by definition,
for any $c\in C$, the class $\phi$ 
belongs to the subgroup
$H_2(\X_{b}/\Y_{\kappa(c)}) \subset H_2(\X_{b},\Q).$
It follows, by the first claim of Lemma~\ref{omega.lifts}, that  $\phi$
is orthogonal
to $[\om]_{b}$, for {\em all} elements $c \in C$. In particular,
this applies to generic elements in $B_i \subset C$ for each
$i\geq 1$. By Proposition~\ref{on.c.i}, we deduce that $\phi \in
H_2(X,\Q)$ is orthogonal to all the subspaces $B_i \subset C$. Thus,
 $\phi = 0$, contradicting the assumption $\phi\neq 0$. 

The contradiction implies that $C_\phi\neq C$, for
any nonzero class  $\phi\in  H_2(X,\Q).$
Since there is only a countable set of homology classes
$\phi \in H_2(X,\Q)$, it follows that there exists
$c\in C$ which is not contained in $C_\phi$, for any $\phi$.
It follows that,
for such a $c$, no nonzero element in $H_2(X,\Q)$
belongs to $H_2(\X_{b}/\Y_{\kappa(c)})$ for any $b$ lying
over $\kappa(c)$. This means  $H_2(\X_{b}/\Y_{\kappa(c)})=0$.
We deduce, similarly to the argument
in the proof of Proposition \ref{vit_Z},
 that the fibers of
$\X_{b}\to\Y_{\kappa(c)}$ do not contain projective curves. As in the
proof of Theorem~\ref{vit1}, this implies that the map $\X_{b} \to
\Y_{\kappa(c)}$ is an isomorphism. Thus, $c\in C^{\tt{GOOD}}$,
and we are done.
 \endproof

\subsection{Restriction to strata and the end of the
proof.}\label{on.c.i.sub} 

It remains to prove Proposition~\ref{on.c.i}. To do this, we need
some information on the Calogero-Moser deformation $\M/C$ and on the
geometry of the resolution $X/Y$.

Somewhat surprisingly, we need to know very little about the
Calogero-Moser deformation -- it suffices to know how it behaves
with respect to the change of the subgroup $G \subset Sp(V)$. The
precise statement is as follows. Let $G_1 \subset G$ be a subgroup,
let $Y_1 = V/G_1$ be the associated symplectic quotient singularity,
and let $\M_1/C_1$ be its Calogero-Moser deformation with its base
$C_1$. We have a canonical restriction map $C \to C_1$. Denote by
$B_1 \subset C_1$ its image, and consider the splitting $\iota:B_1
\to C$ of the projection $C \to B_1$ defined by
$$
\iota(c)(g) = 0 \text{  unless  } g \in G \text{  is conjugate to an
element  } g_1 \in G_1.
$$

\begin{prop}\label{red}
The canonical projection $\eta:Y' \to Y$ extends to a commutative
diagram 
\begin{equation}\label{red.dia}
\begin{CD}
\M_1 \times_{C_1} B_1 @>{\eta}>> \M\\
@VVV @VVV\\
B_1 @>{\iota}>> C
\end{CD},
\end{equation}
where $\M_1 \times_{C_1} B_1$ is the restriction of the
Calogero-Moser deformation $\M_1/C_1$ to $B_1 \subset
C_1$.\endproof
\end{prop}

Proposition~\ref{red} immediately follows from the definition of the
Calogero-Moser deformation, see \cite{EG}. We will aply it to the
terms in the direct sum decomposition \eqref{splt}.

\begin{remark}{\em (Added on Feb. 16, 2010.)} As G. Bellamy kindly
  indicated to us, the above Proposition is not contained in
  \cite{EG}, nor does it follow directly from the definition, and in
  fact it is not clear if the fact is true or not. Fortunately, in
  his beautiful recent paper {\tt arXiv:1001.0239}, I. Losev has
  proved that if one chooses a point $v \in V$, and lets $G_1
  \subset G$ be the stabilizer of this point, then a slightly
  stronger statement becomes true after completing both sides near
  $v$. This is Theorem 1.2.1 of {\tt arXiv:1001.0239}. Losev's
  Theorem is sufficient to save our proof of
  Proposition~\ref{on.c.i}, the only place where
  Proposition~\ref{red} is used.
\end{remark}

Recall (see Section~\ref{hp2.sec}) that the subspaces $B_i \subset
C$ correspond to co\-di\-men\-si\-on\--$2$ strata $U_i$ in the
natural stratification of the quotient variety $Y = V/G$. Every
stratum $U_i$ is of the form $U_i = V_i^o/G'_i$, where $V_i \subset
V$ is a symplectic vector subspace of codimension $2$, $G_i' \subset
G$ is the subgroup of elements which preserve $V_i \subset V$, and
$V_i^o \subset V_i$ is the open subset of elements with minimal
possible stabilizer. This stabilizer is a subgroup $G_i \subset G_i'
\subset G$. All elements of the subgroup $G_i$ are symplectic
reflections. They act trivially on the subspace $V_i \subset
V$. This induces a natural action of the group $G_i$ on the
$2$-dimensional quotient $W_i = V/V_i$. The Calogero-Moser
deformation of the quotient $V/G_i$ is simply the product of the
Calogero-Moser deformation of the quotient $Y_i = W_i/G_i$ and the
symplectic vector space $V_i$. Its base $C_i$ is the space of
$G_i$-invariant $\C$-valued functions on the set of all non-trivial
elements $g \in G_i$. The image $B_i \subset C_i$ of the
restriction map $C \to C_i$ is the subspace
$$
B_i = C_i^{H_i} \subset C_i
$$
of functions invariant with respect to the natural action of the
group $H_i = G_i'/G_i$. By Proposition~\ref{red}, for every stratum
$U_i$ we have a commutative diagram \eqref{red.dia}, where the
scheme in the left top corner is the product $\M_i \times V_i$ of
the Calogero-Moser deformation of the quotient $Y_i$ and the
symplectic vector space $V_i$. The Calogero-Moser deformation $\M_i$
is the standard deformation described in Section~\ref{dim2}.

These are all the properties of the Calogero-Moser deformation that
we will need. Proposition~\ref{on.c.i}, hence also Theorem~\ref{res}
will hold for every Poisson deformation $\M/C$ of the quotient
singularity $Y = V/G$ which admits a diagram \eqref{red.dia} for
every $G_i \subset G$, with $\M' \cong \M_i \times V_i$.

Aside from Proposition~\ref{red}, the decomposition \eqref{splt} has
a very clear geometric interpretation in terms of the isomorphism
$H^2(X,\C) \cong B$ provided by the generalized McKay
correspondence \cite{K}. Namely, let $y \in U_i$ be an arbitrary
point in the stratum $U_i \subset Y$. The map $Y_i \times U_i \to Y$
induced by the projection $V/G_i \to V/G$ is \'etale in
$y$. Therefore the formal neighborhood $\wh{\Y}_y$ of the point $y
\in Y$ naturally decomposes
\begin{equation}\label{decomp}
\wh{\Y}_y = \wh{Y_i} \times \wh{V_i},\qquad\text{product in the sense of
formal schemes},
\end{equation}
where $\wh{Y_i}$ is the completion of the quotient $Y_i = W_i/G_i$
at the origin $o \in Y_i$, and $\wh{V_i}$ is the completion of the
vector space $V_i$ at $0 \subset V_i$. Let $\wh{X}_y$ be the
completion of the variety $X$ at the fiber $\pi^{-1}(y) \subset
X$. Then by \cite[Proposition 5.2]{K1} the decomposition
\eqref{decomp} lifts to the formal scheme $\wh{X}_y$. Namely, we
have a canonical decomposition
$$
\wh{X}_y = \wh{X_i} \times \wh{V_i},
$$
where $\wh{X_i}$ is the completion of a smooth symplectic resolution
$X_i/Y_i$ of the $2$-dimensional quotient $Y_i$ at the exceptional
divisor $E \subset X_i$ of the map $X_i \to Y_i$. The cohomology
space
$
H^2(\wh{X_i}) \cong H^2(X_i) \cong H^2(E)
$
is naturally identified with the space $C_i$, and the natural map
$B \to C_i$ is given by the restriction map $H^2(X,\C) \to
H^2(E,\C)$.

\proof[Proof of Proposition~\ref{on.c.i}.]  We begin with
\thetag{i}. Consider the restriction $\M/B_i$ of the Calogero-Moser
deformation $\M/C$ to the subspace $B_i \subset C$. Let $\eta:\M_i
\times V_i \to \M/B_i$ be the projection provided by
Proposition~\ref{red}. It suffices to prove that $[\om]_{b}$
maps to $0$ under every projection $B \to B_j$ with $j \neq i$.

Choose an arbitrary point $y \in U_j$ and consider the completion
$\wh{\Y}_y$ together with the induced deformation
$\wh{\M}_y/\wh{B_i}$ (here $\wh{B_i}$ is the completion of the
space $B_i$ at $0 \subset B_i$.  Since $j \neq i$, there exists a
point $y' \in Y_i \times V_i$ such that $\eta(y') = y$, the map
$\eta$ is \'etale at $y'$, and the scheme $Y_i \times V_i$ is smooth
at $y' \in Y_i \times V_i$. Then the completion $\wh{(Y_i \times
V_i)}_{y'}$ of $Y_i \times V_i$ at $y'$ is isomorphic to the
completion $\wh{V}$ of the vector space $V$ at $0$, and the
projection
$$
\wh{V} \cong \wh{(Y_i \times V_i)}_{y'}\,
 \longrightarrow\, \wh{\Y}_y \cong \wh{Y_j}
\times \wh{V_j}
$$
is induced by the quotient map $V \to Y_j \times V_j \cong
V/G_j$. By Lemma~\ref{covers.trivial}, this implies that the Poisson
cocycle $\Theta \in HP^2(\wh{\M}_y/\wh{B_i})$ of the Poisson scheme
$\wh{\M}_y/\wh{B_i}$ is a coboundary. In other words, we have
$\Theta = d\xi$ for some vertical vector field $\xi$ on
$\wh{\M}_y/\wh{B_i}$.

Now replace, if necessary, the base $B_i$ of the deformation
$\M/B_i$ with a finite cover $B'_i$, and consider the resolution
$\pi:\X \to \M/B'_i$ provided by Theorem~\ref{vit1}. Let
$\wh{\X}_y$ be its completion at the fiber $\pi^{-1}(y) \subset
\X$. The map $B \to B_j \subset C_j$ is given by the restriction
from $B \cong H^2(X,\C)$ to $H^2(X_j,\C) \cong B_j$ and factors
through the restriction to $H^2(\wh{\X}_y,\C)$. By
Lemma~\ref{vf.extends}, the vector field $\xi$ extends to a vertical
vector field on $\wh{\X}_y/B'_i$. Since $\Theta = d\xi$, the Cartan
homotopy formula gives
$
\om = d(\om \cntrct \theta),
$
where $\om \in \Omega^2(\wh{\X}_y/B'_i)$ is the relative
symplectic form. This implies that $[\om]=0$ tautologically on
$\wh{\X}_y$. Therefore the same is true on the whole variety
$\X/B'_j$ and for every fiber $\X_{b}$.

This proves \thetag{i}. To prove \thetag{ii}, denote by $C_i^\perp
\subset H_2(X,\Q)$ the orthogonal to the subspace $B_i \subset B
\cong H^2(X,\C)$, and denote by $\om^\perp \subset H_2(X,\Q)$ the
orthogonal to the cohomology class $[\om]_{b} \in H^2(X,\C)$. We
have to prove that the embedding $C_i^\perp \subset \om^\perp$ is
in fact an equality. To prove this, it suffices to show that 
$$
\dim \om^\perp \leq \dim C_i^\perp = \dim B - \dim B_i.
$$
To prove this inequality, it suffices to exhibit a $\Q$-vector space
subspace $P$ of dimension $\dim P = \dim B_i$ and a map $f:P \to
H_2(X,\Q)$ such that the pairing with $[\om]_{b}$ induces an
embedding $P \to \C$. We claim that we can take
$$
P = \left(H_2(X_i,\Q)\right)^{H_i} = \left(H_2(E_i,\Q)\right)^{H_i}
$$ 
and $f = \eta_*:H_2(E_i,\Q) \to H_2(X,\Q)$.

Indeed, by Lemma~\ref{omega.lifts} the class $[\om]_{b} \in
H^2(\X_{b},\C)$ comes from a class $[\om] \in
H^2(\M_c,\C)$. Moreover, the class $[\om] \in H^2(\M_c,\C)$ is
canonical. In particular, the restriction $\eta^*[\om] \in
H^2((\M_i)_c \times V_i,\C)$ coincides with the class $[\om]_c
\otimes 1$, where $[\om]_c \in H^2((\M_i)_c,\C)$ is the class
corresponding to the Calogero-Moser deformation $\M_i$. But the
variety $Y_i$ is $2$-dimensional. Therefore the Calogero-Moser
deformation $\M_i$ comes from the universal symplectic deformation
of the resolution $X_i$, and the class $[\om]_c \in C_i \cong
H^2(X_i,\C)$ coincides with $c \in B_i = C_i^{H_i} \subset
C_i$. Thus we have to show that for a generic element $c \in
C_i^{H_i} \cong P^* \otimes \C$, the pairing with $c$ induces an
embedding $P \to \C$. This is clear.  \endproof

\section{Applications of Hochschild cohomology.}\label{zcsmooth.sec}
The primary goal of this section is to
 prove Theorem~\ref{zcsmooth}. We begin however with some
general results that relate Hochschild cohomology
to orbifold cohomology.

Let $A$ be an associative algebra. Recall
the Gerstenhaber bracket $[-,-]: \HHH^p(A) \times
\HHH^q(A) \to \HHH^{p+q-1}(A)$ on  Hochschild
cohomology (see \cite{Lo} and \S7 below).
Given  a Hochschild cocyle
$\Theta \in \HHH^2(A)$ such that $[\Theta,\Theta] = 0,$
we 
introduce a {\em
twisted Hochschild cohomology algebra}
 of $A$ as follows.

\begin{defn}\label{twisted_Hoch}
For
$\Theta \in \HHH^2(A)$ such that $[\Theta,\Theta] = 0$,
define
 $\HHH^\hdot_\Theta(A)$ to be 
the cohomology of the complex $\bigl(\HHH^\hdot(A),\,d_\Theta\bigr)$ where
the differential $d_\Theta:\HHH^\hdot(A) \to \HHH^{\hdot+1}(A)$ is
given by $d_\Theta(a) =
[\Theta,a]$.
\end{defn}

Let $A=\C[M]$ be the  coordinate ring
of  a smooth affine algebraic variety $M$.
Then we have $\HHH^\hdot(A)= \Gamma(M, \Lambda^\hdot{\mathcal{T}}_M)$,
by  Hochschild-Kostant-Rosenberg theorem, cf. \cite{Lo}.
Furthermore,  the Gerstenhaber bracket on
$\HHH^\hdot(A)$ reduces, in this case,
to the Schouten bracket on
polyvector fields, cf. \S\ref{smooth_case}. Thus, the
cocycle $\Theta\in \HHH^2(A)$ may be viewed as a bivector.
The equation  $[\Theta,\Theta] = 0$ says that this bivector
gives a Poisson structure on $A$. Moreover,
 according
to \eqref{PBW_smooth} of the Appendix, there is a natural isomorphism
$\HHH_\Theta^\hdot(A) \cong HP^\hdot(A)$. 
Note that if $M$ is
not smooth, both sides in the isomorphism are
still well-defined, but they are not 
necessarily isomorphic
 any more. 

Assume  that 
the smooth affine  variety $M$ is a symplectic
manifold.
We let the 2-cocycle $\Theta$ be the bivector on $M$
corresponding to (the inverse of) the symplectic 2-form.
In such a case, the Poisson cohomology
reduces to  De Rham cohomology, see \cite{B}
and section \ref{smooth_case} below,
hence we obtain graded algebra isomorphisms
\begin{equation}\label{twist_derham}
\HHH_\Theta^\hdot(\C[M]) \cong HP^\hdot(\C[M]) \cong H^\hdot(M,\C).
\end{equation}

Assume next that 
a finite group $G$ acts on $M$, a  smooth affine symplectic
variety, by symplectic automorphisms. The action on
$M$ induces one on $\C[M]$, and we form the cross-product
algebra
$$A:=\C[M]\smsh G.
$$

The bivector  corresponding to (the inverse of) the symplectic form
on $M$ is $G$-invariant, and the class 
$\Theta\in \HHH^2(\C[M])= \Gamma(M, \Lambda^\hdot{\mathcal{T}}_M)$
gives rise, by the standard deformation theory,
to 
a $G$-equivariant 
first-order infinitesimal
deformation of the associative algebra $\C[M]$.
Moreover, the equation $[\Theta,\Theta]=0$
insures that this first-order deformation may be
 extended (not uniquely) to a  $G$-equivariant 
second order deformation.
The latter gives rise, via the  cross-product construction,
 to a  second order
deformation of $\C[M]\smsh G$,
hence, to a class $\widetilde\Theta\in \HHH^2(A)=\HHH^2(\C[M]\smsh G)$, 
such that $[\widetilde\Theta, \widetilde\Theta]=0$.
We form  the corresponding twisted Hochschild cohomology algebra
$\HHH_{\widetilde\Theta}^\hdot(\C[M]\smsh G)$.

On the other hand, associated with the $G$-action on
$M$ one has the orbifold  cohomology
algebra 
$H^\hdot_{\tt{orb}}(M;G)$, see \eqref{Horb}.
We have

\begin{prop}\label{HHorb} Given a  symplectic $G$-action 
on an affine symplectic manifold, there is a natural
graded algebra isomorphism
$$
\HHH_{\widetilde\Theta}^\hdot(\C[M]\smsh G)\,\cong\,
H^\hdot_{\tt{orb}}(M;G)\,\left(=
\Bigl(\bigoplus\nolimits_{g\in G}\,\,
H^{\hdot-\dim M^g}(M^g)\Bigr)^G\right).
$$
\end{prop}

\begin{proof}\;
 The action of the finite group $G$ on $\C[M]$
being semisimple, the Hochschild cohomology
of the cross-product may be expressed in terms of
Ext-groups of $\C[M]$-bimodules as follows, cf. \cite{Al}: 
\begin{align}\label{bimod}
\HHH^\hdot(\C[M]\rtimes G)&=\Ext^\hdot_{\C[M]{\tt{-bimod}}}
\bigl(\C[M]\,,\,\C[M]\# G\bigr)^G\nonumber\\
&=
\Bigl(\bigoplus\nolimits_{g\in G} \Ext^\hdot_{\C[M\times M]{\tt{-mod}}}
\bigl(\C[M]\,,\,\C[M]\cdot g\bigr)\Bigr)^G.
\end{align}
Here $G$-invariants are taken with respect to the adjoint 
$G$-action, and in
the rightmost term of the formula  we identify 
$\C[M]$-bimodules with $\C[M\times M]$-modules.
Thus
$\C[M]\cdot g$ stands for the $\C[M\times M]$-module arising
from 
the coordinate ring of the graph-subvariety
${\tt{Graph}}(g: M \to M)\, \subset\, M\times M.$

In general, let $W$ be a smooth affine variety containing
 two smooth (closed) subvarieties $E,F\subset W,$  such that
the intersection $E\cap F$ is {\em clean}. The latter means
that  $E\cap F$ is smooth and that
the equation $\T_{E\cap F}={\T_E|}_{_{E\cap F}}\cap
{\T_F|}_{_{E\cap F}}$ holds for the corresponding tangent bundles.
In such a case, we put
$$d:=\dim W-\dim E - \dim F +\dim (E\cap F).
$$

A standard argument based on Koszul complexes shows that
the group $\Ext^i_{\C[W]{\tt{-mod}}}(\C[E]\,,\,\C[F])$
vanishes for all $i< d$, and for $i\geq d$ we have
\begin{align}\label{ext_gen}
\Ext^i_{\C[W]{\tt{-mod}}}&(\C[E]\,,\,\C[F])=\nonumber\\
&=
\Gamma\left(E\cap F,\,
\Lambda^{i-d}\bigl(\frac{{\T_W|}_{_{E\cap F}}}{{\T_E|}_{_{E\cap F}}+
{\T_F|}_{_{E\cap F}}}\bigr)\otimes
\det\bigl(\frac{{\T_E|}_{_{E\cap F}}+
{\T_F|}_{_{E\cap F}}}{{\T_E|}_{_{E\cap F}}}\bigr)\right),
\end{align}
where $\det(\ldots)$ denotes the top wedge power.

For any $g\in G$, there is a natural $g$-action
on  ${\T_M|}_{_{M^g}}$ by  vector bundle endomorphisms, and we have a
canonical
direct sum decomposition 
$${\T_M|}_{_{M^g}}={\op{Image}}(\id - g) \bigoplus
\T_{M^g},
$$
where ${\op{Image}}(\id - g)$ denotes the
subbundle (of locally constant rank) formed by
the  images of the fiberwise action of the operator $\id-g$.

Writing $E \subset M\times M$ for the diagonal, we get
$E\cap {\tt{Graph}}(g) \cong M^g$.
We deduce canonical isomorphisms
\begin{align*}
&\frac{{\T_E|}_{_{M^g}}+{\T_{{\tt{Graph}}(g)}|}_{_{M^g}}}{{\T_E|}_{_{M^g}}}
\,\cong\, {\op{Image}}(\id - g),\quad\text{and},\\
&\frac{{\T_{M\times M}|}_{_{M^g}}}{{\T_E|}_{_{M^g}}+{\T_{{\tt{Graph}}(g)}|}_{_{M^g}}}
\,\cong\, ({\T_M|}_{_{M^g}})\big/{\op{Image}}(\id - g)\, \cong\,
\T_{M^g}.
\end{align*}
Thus, $M^g$ is a symplectic submanifold and $d=
\dim M^g$. Further, the symplectic form restricts to a
non-degenerate 2-form on the fibers of the vector
bundle ${\op{Image}}(\id - g)$,
and this gives canonical trivializations
$\det\bigl({\op{Image}}(\id - g)\bigr)\cong\C$. Being canonical, the
trivializations
are
compatible with $G$-action (that permutes the fixed point sets $M^g$ for various 
elements $g$).
Thus,
from \eqref{bimod} and \eqref{ext_gen} we obtain 
$$\HHH^i(\C[M]\rtimes G) \cong\Bigl(\bigoplus\nolimits_{g\in G} \,\,
\Gamma(M^g,\, \Lambda^{i-\dim M^g}\T_{M^g})\Bigr)^G.
$$
Further, the symplectic form on $M^g$ induces an vector bundle
isomorphism
$\T_{M^g} \cong (\T_{M^g})^*$, hence a graded algebra isomorphism
$\Gamma(M^g, \Lambda^\hdot \T_{M^g})\cong
\Omega^\hdot(M^g)$.
Combining the formulas above, we obtain
\begin{equation}\label{HH1}
\HHH^\hdot(\C[M]\rtimes G) \cong\Bigl(\bigoplus\nolimits_{g\in G}\,\,
\Omega^{\hdot-\dim M^g}(M^g)\Bigr)^G.
\end{equation}

Observe next that 
the cochain ${\widetilde\Theta}\in \HHH^2(\C[M]\rtimes G)$ is given, essentially,
by the bivector corresponding to the symplectic 2-form. 
Therefore, it is easy to see (as e.g. in \cite{B}),
that the differential $[{\widetilde\Theta},-]$ on Hochschild cohomology
gets transported under the isomorphism \eqref{HH1} to
the direct sum of the de Rham differentials
$d: \Omega^{\hdot-\dim M^g}(M^g)\to \Omega^{\hdot+1-\dim M^g}(M^g)$.
The cohomology of the latter is nothing but
$H^{\hdot-\dim M^g}(M^g, \C)$, the singular cohomology
of $M^g$ (up to shift). 

This establishes the isomorphism of the Proposition.
Compatibility of the isomorphism with the algebra
structures is more difficult (it involves 
Kontsevich's theorem, see
\cite[\S8.4]{Ko}, on the {\em cup-product on
tangent cohomology}), and it will not be given here.
 Below, see \eqref{HH2}, we will only use
 a very special case of Proposition \ref{HHorb}
where such a compatibility is immediate from definitions.
\end{proof}

For the rest of this section, assume  $M=V$ a symplectic vector space,
and $G\subset Sp(V)$, a finite group, so that
the de Rham cohomology of each fixed point set $V^g$
 is trivial in all degrees but zero.
Thus, 
 writing $\C\langle k\rangle$ for a 1-dimensional  graded
vector space concentrated in degree $k$,
the  orbifold cohomology algebra
in the RHS of the isomorphism of Proposition
\ref{HHorb} reads
\begin{equation}\label{HH2}
H^\hdot_{\tt{orb}}(V;G)\cong\Bigl(\bigoplus\nolimits_{g\in G} \,
\C\langle \dim V^g\rangle \Bigr)^G \cong\bigl(\gr^F_\idot\C[G]\bigr)^G 
\cong\gr^F_\idot(\Zf{G}),
\end{equation}
where the associated graded algebra is taken with respect
to the filtration $F_\idot(\C[G])$ considered in \S1.

In \cite{EG}, the authors construct a certain deformation $\Hf_{t,c}$
of the cross-product algebra $\Hf_{0,0}:= \C[V] \smsh G,
$ which is parametrized by an affine line with
coordinate $t$ and the space $C \cong\gr^F_2(\Zf{G})$ of $G$-invariant functions on the
set on symplectic reflections in $G$. When $t=0$, the algebra
$\Hf_{0,c}$ has a large center $\Zf_c$. The Calogero-Moser space
$\M_c$ is obtained by taking its spectrum, $\M_c = \Spec \Zf_c$. The
deformation in the $t$-direction induces a Poisson structure on the
algebra $\Zf_c$ and on the variety $\M_c$. When, on the other hand,
$t$ is generic, the center of the algebra $\Hf_{t,c}$ is the
one-dimensional $\C$-vector space spanned by the unit element.

In \cite{EG}, the authors consider the
Hochschild cohomology  $\HHH^\hdot(\Hf_{t,c})$.
They
construct a canonical map $\chi:\gr^F_\idot(\Zf{G}) \to
\HHH^\hdot(\Hf_{t,c})$, which is
shown to be an isomorphism for
generic $t$. Moreover, the map $\chi$ is compatible with the
deformation $\Hf_{t,c}$, i.e., for every pair $\langle t_0,c_0 \rangle$,
the composite map:
$$
C\iso \gr^F_2(\Zf{G}) \stackrel{\chi}{\longrightarrow} \HHH^2(\Hf_{t_0,c_0})
$$
is the Kodaira-Spencer map for the family $\Hf_{t,c}$ near the point
$\langle t_0,c_0 \rangle$. In paritucular, if one fixes a generic
enough $t=t_0$, then the family $\Hf_{t_0,c}/C$ is the universal
deformation of the associative algebra $\Hf_{t_0,c_0}$.

When $t=0$, the map $\chi$ is still defined, but its image no longer
generates the Hochschild cohomology groups $\HHH^\hdot(\Hf_{0,c})$. In
fact, these groups become infinite-dimensional as $\C$-vector spaces.

Now, fix a $c \in C$, and consider $\Hf_{t,c}$ as a
family depending on $t$. Applying the Kodaira-Spencer map at the
point $t=0$, we obtain a cohomology class $\Theta_c \in
\HHH^2(\Hf_{0,c})$ which satisfies $[\Theta_c,\Theta_c]=0$.
Therefore, we may apply  Definition \ref{twisted_Hoch} and form
twisted Hochschild cohomology  $\HHH^\hdot_{\Theta_c}(\Hf_{0,c})$.

\begin{lemma}\label{twisted}
The canonical map $\chi:\gr^F_\idot(\Zf{G}) \to \HHH^\hdot(\Hf_{0,c})$
descends to a map $\chi:\gr^F_\idot(\Zf{G}) \to
\HHH^\hdot_{\Theta_c}(\Hf_{0,c})$, and the latter map is an isomorphism.
\end{lemma}

\proof{} The map
$\chi:\gr^F_\idot(\Zf{G}) \to \HHH^\hdot(\Hf_{t,c})$
is defined for all $(t,c)$.
It follows that, for $t=0$,
the image of this map commutes with $\Theta_0$, hence,
the map $\chi$ induces a well-defined map $\gr^F_\idot(\Zf{G}) \to
\HHH^\hdot_{\Theta_c}(\Hf_{0,c})$.

First, consider the case $c=0$ where $\Hf_{0,0}=\C[V]\# G$.
Applying Proposition \ref{HHorb} and using \eqref{HH2},
we obtain a graded algebra isomorphism
$\HHH^\hdot_{\Theta_0}(\Hf_{0,0})\cong\gr(\Zf{G}),$
in particular, all odd twisted Hochschild cohomology groups
of $\Hf_{0,0}$ vanish. Moreover, a
 calculation carried out in \cite{arg} for a Weyl algebra instead
of the polynomial algebra $\C[V]$ shows that the isomorphisms above
is the inverse of the map 
$\chi:\gr^F_\idot(\Zf{G}) \to \HHH^\hdot_{\Theta_0}(\Hf_{0,0})$ considered in
\cite{EG}. This completes the proof in the special case: $c=0$.

To complete the proof in the general case recall from \cite{EG}
that, for any $c$,
the algebra $\Hf_{0,c}$ comes equipped with a
canonical increasing filtration such that the associated graded
algebra $\gr\Hf_{0,c}$ is isomorphic to the algebra $\Hf_{0,0}$, see \cite{EG}. This
filtration induces a filtration on the
twisted  Hochschild complex
which is compatible  with the differential $[\Theta,-]$, and therefore gives rise
to a
spectral sequence 
$$  \HHH^\hdot_{\Theta_0}(\Hf_{0,0})=E^{\hdot,\hdot}_1\quad
\Longrightarrow\quad
E^{\hdot,\hdot}_\infty=\gr^\hdot\!\bigl(\HHH^\hdot_{\Theta_c}(\Hf_{0,c})\bigr)\,.
$$
Since  $\Hf_{0,0}$
has no odd twisted Hochschild
cohomology,
all the
differentials in this spectral sequence vanish,
 and we have: $E^{\hdot,\hdot}_1=E^{\hdot,\hdot}_\infty$.
Thus, it follows from the case $c=0$ of the Lemma that the
map $\gr(\chi):\,\gr^F_\idot(\Zf{G}) \to
\gr\bigl(\HHH^\hdot_{\Theta_c}(\Hf_{0,c})\bigr)$
is a bijection. Therefore, for any $c$, the map
$\chi: \gr^F_\idot(\Zf{G}) \to\HHH^\hdot_{\Theta_c}(\Hf_{0,c})$
is also a bijection, and the Lemma is proved.
\endproof

We now apply the above result assuming in addition that  $c \in C$ 
is such that the
Calogero-Moser space $\M_{c}$ is smooth.

\begin{lemma}\label{sm=>uni}
Let $c \in C$ be such that the Calogero-Moser space $\M_{c}=
\Spec\B_c$ is
smooth. Then the Calogero-Moser family $\M/C$ considered over the
formal neighborhood of the point $c$ gives a universal Poisson
deformation of the Poisson variety $\M_{c}$.
\end{lemma}

\proof{} 
It is proved in \cite{EG} that whenever $\M_{c} =
\Spec\B_{c}$ is smooth, the algebra $\Hf_{0,c}$ is
Morita-equivalent to  $\B_{c}$. Therefore we have
$\HHH^\hdot(\Hf_{0,c}) \cong \HHH^\hdot(\B_{c})$ and
$\HHH^\hdot_{\Theta_{c}}(\Hf_{0,c}) \cong
\HHH^\hdot_{\Theta_{c}}(\B_{c})$. 
Using the smoothness of $\M_{c}$ again,
 we conclude, see \eqref{PBW_smooth}, that the right-hand side is
isomorphic to $HP^\hdot(\B_{c})$. Combining this with
Lemma~\ref{twisted}, we see that the Kodaira-Spencer map for the
family $\M/C$ computed at the point $c \in C$ induces an
isomorphism $C \cong HP^2(\B_{c})$.
\endproof

We can now prove Theorem~\ref{zcsmooth} using the
gradings and the dimension estimate obtained in
Section~\ref{hp2.sec}.

\proof[Proof of Theorem~\ref{zcsmooth}.]
Denote $S' = HP^2(\C[Y])$, and let $S \subset S'$ be the base of the
universal graded Poisson deformation of the graded Poisson algebra
$\C[Y]=\C[V]^G$. By Lemma~\ref{sm=>uni}, the classifying map
$\kappa:C \to S$ of the Calogero-Moser family $\M/C$ is \'etale at a
generic point $c \in C$. In particular, its differential
$$
d\kappa:T_cC \to T_{\kappa(c)}S'
$$
is injective. Since, $\dim S' \leq \bn = \dim C$, this differential
is also surjective, and $S=S'$. 

It remains to prove that the map $\kappa:C \to S=S'$ is surjective.
The Calogero-Moser deformation $\M/C$ is graded, with grading given
by assigning degree $2$ to every element in the space $C$. Therefore
the map $\kappa:C \to S \subset S'$ is compatible with the
gradings. Taking quotients with respect to $\C^*$, we obtain a
rational map $\overline{\kappa}:\Pp(C) \ratto \Pp(S')$, where
$\Pp(C) \cong {\mathbb{P}}^{\bn-1}$ is the projectivization of $C$, and
$\Pp(S')$ is the projectivization of the vector space $S'$ (with its
grading, whatever it may be). Its image $\overline{\kappa}(\Pp(C))
\subset \Pp(S')$ is a closed subvariety. Since $\kappa:C \to S=S'$
is generically \'etale, the image $\overline{\kappa}(\Pp(C)) \subset
\Pp(S)$ coincides with the whole $\Pp(S)$.\endproof

\section{Appendix.}
\renewcommand{\thesection}{A}
\label{appendix}

\subsection{Harrison cohomology.}\label{har.def} 
Consider an arbitrary vector space $A$ over $\C$ (or, more
generally, over an arbitrary field of characteristic $0$). Let
$L_\idot A$ be the free graded Lie coalgebra generated by the vector
space $A$ placed in degree $-1$, so that $L_1 A = A$ and $L_2 A =
S^2A$, the symmetric square of $A$. Denote by $DL^\hdot(A)$ the
graded Lie algebra of coderivations of the Lie coalgebra $L_\idot
A$. Since $L_\idot$ is free, we have
$
DL^k(A) \cong \Hom(L_k A,A).
$
In particular, every map $m:S^2A \to A$ extends to a coderivation $d
\in DL^1$. The following is easily checked by a direct computation.

\begin{lemma}\label{har}
A commutative product $m:S^2A \to A$ is associative if and only if
the corresponding coderivation $m \in DL^1$ satisfies
$
\{m,m\} = 0: L_3 A \to A. 
$\hfill$\Box$
\end{lemma}

Assume that $A$ is a commutative associative algebra. Then we can
apply the Lemma and obtain a canonical coderivation $m \in DL^1(A)$
satisfying $\{m,m\}=0$. Setting
$
a \mapsto \{m,a\}
$
defines a differential $d:DL^\hdot(A) \to DL^{\hdot+1}(A)$ and turns
$DL^\hdot(A)$ into a DG Lie algebra. The differential $d:DL^0(A) \to
DL^1(A) \cong \Hom(S^2A,A)$ can be explicitly described in the
following way:
\begin{equation}\label{df}
(df)(a\otimes b) = f(ab) - af(b) - bf(a)\quad,\quad a\otimes b\in S^2A.
\end{equation}
The spaces $DL^\hdot(A) = \Hom(L_\idot A,A)$ carry natural
$A$-module structure -- $A$ acts on the target space $A$. It is easy
to check that $d:DL^\hdot(A) \to DL^{\hdot+1}(A)$ is an $A$-module
map. Moreover, we have
$$
DL^k(A) \cong \Hom(L_k A,A) \cong \Hom_A(L_k A \otimes A,A),
$$
and the the differential $d$ is dual to an $A$-module map
$d:L_{\idot+1} A \otimes A \to L_\idot A \otimes A$. The complex
$\langle L_\idot \otimes A, d \rangle$ is denoted by $\Har_\idot(A)$
and called the {\em Harrison complex} of the commutative
associative algebra $A$. Its first two terms are
$$
\begin{CD}
@>>> S^2A \otimes A @>{d}>> A \otimes A @>>> 0 @>>>
\end{CD}
$$
with the differential given by $d: a  b \otimes c\longmapsto a
\otimes b  c
+ b \otimes a  c - a  b \otimes c$. It is well-known that the Harrison
homology complex is quasiisomorphic to the so-called {\em cotangent
complex} $\Omega_\idot(A)$ of the algebra $A$\footnote{In degree
$0$, this is essentially the standard representation of the K\"ahler
differentials module $\Omega^1A$ as the quotient of the diagonal
ideal $I \subset A \otimes A$ by its square.}.

\subsection{Poisson cohomology.} Let $A$ be a unital commutative
$\C$-algebra with a Poisson bracket
$\{-,-\}:\Lambda^2_\C(A) \to A$ 
that satisfies the Leibniz formula \eqref{compt.eq},
the Jacobi identity
\eqref{jacobi}, and
such  that  $\{1,a\}=0,\,\forall a$.

Let $P_\idot(A)$ be the free graded Poisson coalgebra generated by
the vector space $A$ placed in degree $-1$ (here `free' means
that our coalgebra is a universal object in the category
of graded Poisson coalgebras). It is easy to see, cf. e.g. \cite{Fr},
 that the coalgebra
$P_\idot(A)$ may be constructed as the free super-symmetric (co)algebra on
the vector space $L_\idot(A)$(= free graded Lie coalgebra
generated by $A$). In more details, for each $m\geq 0$, one has a canonical
decomposition
\begin{equation}\label{pbw.P}
P_m(A) = \bigoplus_{p+q=m}P_{p,q}\quad\text{where}\quad
P_{\idot,k} = \Lambda^k(L_\idot(A)), \qquad k =0,1,\ldots .
\end{equation}
This gives a bigrading $P(A)= \bigoplus P_{p,q}\,,\,P_{p,q}=P_{p,q}(A)$
such that  $P_{\idot,1}(A) = L_\idot(A)$.
The comultiplication in $P_\idot(A)$ preserves the bigrading, while
the cobracket is of bidegree $(0,1)$.

Denote by $DP^{\hdot,\hdot}(A)$ the graded Lie algebra of
coderivations of the Poisson coalgebra $P_{\idot,\idot}$. For
reasons of convenience, we will shift the bigrading on
$DP^{\hdot,\hdot}(A)$ by $(1,0)$, so that the first non-trivial term
is $DP^{0,0}(A)$. Since $P_\idot$ is free, we have
$$
DP^k(A) \cong \Hom(P_k(A),A), \qquad k \geq 0.
$$
Moreover, \eqref{pbw.P} gives an identification
\begin{equation}\label{pbw.DP}
DP^{k,\hdot}(A) \cong \Hom(P_{k,\idot}(A),A) \cong
\Hom(\Lambda^k(DL^\hdot(A)),A).
\end{equation}
In particular, we have $DP^{0,2}(A) = \Hom(\Lambda^2(A),A)$, the
space of all skew-commutative binary operations on $A$.

Assume now that $A$ is an associative commutative algebra. Then by
definition, we have the Harrison complex $\Har_\idot(A)$ of
flat (in fact, free) $A$-modules. The identification \eqref{pbw.DP}
can be re-written as
\begin{equation}\label{hp.def.eq}
DP^{\idot,k}(A) \cong \Hom_A(\Lambda^k_A(\Har_\idot(A)),A), \qquad k
\geq 0,
\end{equation}
and the differential in the Harrison complex
$\Har_\idot(A)$ induces a differential $d:DP^{\hdot,\hdot}(A) \to
DP^{\hdot,\hdot+1}(A)$. This turns the graded Lie algebra
$DP^{\hdot}(A)$ into a DG Lie algebra.

\begin{lemma}\label{compt}\mbox{}
\begin{enumerate}
\item A skew-linear operation $\{-,-\}:\Lambda^2(A) \to A$ satisfies
the Leibnitz rule \eqref{compt.eq} if and only if for the corresponding
element
$$
\Theta \in DP^{0,2}(A) \cong \Hom(\Lambda^2(A),A)\quad\text{we
have}\quad
d\Theta = 0.
$$
\item The operation $\{-,-\}:\Lambda^2(A) \to A$ satisfies the
Jacobi identity \eqref{jacobi} if and only if $[\Theta,\Theta]=0$.
\item Let $1 \in A \cong DP^{0,0}(A)$ be the unit element. We have
$\{-,1\} = 0$ if and only if $[\Theta,1]=0$.
\end{enumerate}
\end{lemma}

\proof{} The first claim immediately follows from \eqref{df}. To
prove \thetag{ii}, note that $DP^{\hdot,0}(A) =
\Hom(\Lambda^\hdot(A),A)$ is in fact the graded Lie algebra of
coderivations of the free skew-commutative coalgebra $\Lambda^k(A)$
generated by $A$. Therefore $[\Theta,\Theta] = 0$ if and only if the
associated map $\{-,-\}:\Lambda^2(A) \to A$ extends to a
coerivation $\delta:\Lambda^{\hdot+1}(A) \to \Lambda^\hdot(A)$
satisfying $\delta \circ \delta = 0$. It is well-known (and easily checked)
that the latter is equivalent to the Jacobi identity \eqref{jacobi}
on the operation $\{-,-\}:\Lambda^2(A) \to A$. Finally,
\thetag{iii} is clear.  \endproof

By Lemma~\ref{compt}, giving a Poisson structure on a commutative
associative algebra $A$ is equivalent to giving an element $\Theta
\in DP^{0,2}(A)$ such that
$$
d\Theta = 0, \qquad\qquad [\Theta,\Theta] = 0, \qquad\qquad [\Theta,1]
=0. 
$$
We will call the element $\Theta \in DP^2(A)$ the {\em Poisson
cochain} corresponding to the Poisson structure on $A$. For every
Poisson algebra $A$, the map 
\begin{equation}\label{hp.def.eq2}
\delta:DP^{\hdot,\hdot}(A) \to
DP^{\hdot+1,\hdot}(A)\quad,\quad a\mapsto \delta(a) := [\Theta,a]
\end{equation}
satisfies $\delta \circ \delta = 0$.

\begin{defn}
The {\em Poisson cohomology} $HP^\hdot(A)$ of
a Poisson algebra $A$ is the cohomology of the total complex
$(DP^\hdot, d+\delta)$ associated
to the {\em Poisson
(bi)complex} $DP^{\hdot,\hdot}(A)$ with differentials
$d:DP^{\hdot,\hdot}(A) \to DP^{\hdot,\hdot+1}(A)$ and
$\delta:DP^{\hdot,\hdot}(A) \to DP^{\hdot+1,\hdot}(A)$.
\end{defn}

Here are the first few terms in the bicomplex $DP^{\hdot,\hdot}(A)$.
{\small
\begin{equation}\label{bicomplex}
\begin{CD}
@. @A{d}AA                    @A{d}AA\\
@. \Hom(L_3A,A) @>{\delta}>>  \Hom(\Lambda^2(S^2A),A)
\oplus \Hom(L_3A\otimes A,A) @>{\delta}>>\\
@. @A{d}AA                    @A{d}AA\\
@. \Hom(S^2A,A) @>{\delta}>>  \Hom(S^2A \otimes A,A) @>{\delta}>>\\
@. @A{d}AA                    @A{d}AA\\
A @>{\delta}>> \Hom(A,A) @>{\delta}>> \Hom(\Lambda^2(A),A)
   @>{\delta}>> 
\end{CD}
\end{equation}}
Here the leftmost vertical column is the Harrison complex
$\Hom(L_\idot{A},A)$, the 2-d vertical column is formed
by graded pieces of $\Lambda^2\bigl(\Hom(L_\idot{A},A)\bigr)$,
etc., and the bottom row is the standard cochain complex
for $A$ viewed as a Lie algebra with respect to the Poisson
bracket.
In particular we get
\smallskip

\begin{enumerate}
\item $HP^0(A)= \big\{z\in A\enspace\big|\enspace
\{z,a\}=0\,,\,\forall a\in A\big\}= \text{ {\it{Poisson center of }}} A$.

\item 
$HP^1(A)=ZP^1/BP^1\enspace$ where
$ZP^1,BP^1\subset\Hom(A,A)$
are defined by
\end{enumerate}
$$ZP^1:=\Big\{f\;\big|\enspace
f(ab)=af(b)+f(a)b\enspace,\enspace f(\{a,b\})=
\{f(a),b\}+\{a,f(b)\}\;,\,\forall a,b\Big\}$$
and $BP^1:=\Big\{f=f_c : a\mapsto \{c,a\}\enspace\big|\enspace c\in A\Big\}.$

Next, let $\phi\in \Hom(S^2A,A)$ and $\psi\in \Hom(\Lambda^2A,A)$.
Define two $\C$-bilinear maps $A\otimes A \to A[\eps]/(\eps^2)$
 by the formulas
\begin{equation}\label{hp.def.def}
a,b \longmapsto a\cdot_\eps b = a\cdot b + \eps\cdot\phi(a,b)\quad,\quad
a,b \longmapsto \{a,b\}_\eps =\{a,b\} + \eps\cdot\psi(a,b)
\end{equation}
It is straightforward to verify that 

\noindent
(iii) $\enspace$\parbox[t]{113mm}{
 $\phi\oplus\psi$ is a Poisson 2-cocycle if and only if
formulas
\eqref{hp.def.def} give a 1-st order infinitesimal deformation of $A$
as  a Poisson algebra. Furthermore, the group $HP^2(A)$
classifies 1-st order infinitesimal deformations up to equivalence.}

\begin{remark}
The  Poisson bracket gives a canonical cocycle  $\Theta \in DP^2(A)$, which may
or may not be trivial in  Poisson cohomology, depending on the algebra $A$.
\end{remark}

\subsection{Relation to Hochschild cohomology and the K\"unneth formula.} 
It is well-known that the  Poisson operad may be viewed as
a `degeneration' of the Associative operad. In this way,
 the  Poisson cohomology bicomplex of a Poisson algebra
may  be viewed as
a `degeneration' of the Hochschild cochain complex.

To explain the analogy, recall that, given a  vector space $A$,
one defines  three graded Lie algebras $DL^\hdot(A)$, 
$DP^\hdot(A)$, and  $DT^\hdot(A)$ as the Lie algebras of coderivations of
the free Lie colagebra $L_\idot(A)$,  free Poisson colagebra
$P_\idot(A)$
or  free associative coalgebra (with counit)
$T_\idot(A)$, respectively. As in the case
of $DP^\hdot(A)$, we shift the grading on $DT^\hdot(A)$ by one, so
that
$$
DT^k(A) \cong \Hom(T_k(A),A) = \Hom(A^{\otimes k},A), \qquad k \geq 0.
$$
Under this identification, the Lie bracket in $DT^\hdot(A)$
becomes the Gerstenhaber bracket given, for every $f \in DT^k(A)$, $g
\in DT^l(A)$,
 by the standard formula
\begin{multline}\label{gerst}
[f,g](a_1 \otimes \dots \otimes a_{k+l-1}) = \\
\begin{split}
&\sum_{1 \leq i \leq l}(-1)^ig(a_1 \otimes \dots \otimes f(a_i \otimes
\dots \otimes a_{i+k-1}) \otimes \dots \otimes a_{k+l-1}) \\ &\quad
- \sum_{1 \leq i \leq k}(-1)^if(a_1 \otimes \dots \otimes g(a_i \otimes
\dots \otimes a_{i+l-1}) \otimes \dots \otimes a_{k+l-1}).
\end{split}
\end{multline}

Recall that the free associative coalgebra
$T_\idot(A)$ is known to be isomorphic to the
universal enveloping coalgebra of the
free Lie coalgebra $L_\idot(A)$.  This means that there is
 a canonical decreasing filtration on the associative  coalgebra
$T_\idot(A)$,
such that   the corresponding 
associated graded $\gr_\idot T_\idot(A)$
is isomorphic, by
Poincar\'e-Birkhoff-Witt theorem, 
 to the free super-symmetric coalgebra on the vector space
$L_\idot(A)$. As we have mentioned earlier, this  free super-symmetric
coalgebra
is nothing but  the free
Poisson coalgebra $P_{\idot,\idot}(A)$ generated by $A$.
In other words, there is a canonical Poisson coalgebra
(bigraded) isomorphism
$$\gr_\idot T_\idot(A)\cong P_{\idot,\idot}(A)\,.$$
Further, the  decreasing filtration on $T_\idot(A)$
induces an increasing
filtration on the graded Lie algebra $DT^\hdot(A)$ which we call the
{\em PBW filtration}. The associated graded  $\gr^\hdot
DT^\hdot(A)$ acts on the associated graded  $\gr_\idot
T_\idot(A) \cong P_{\idot,\idot}(A)$ by Poisson coderivations;
therefore we have a Lie algebra map $\gr^\idot DT^\hdot(A) \to
DP^{\hdot,\hdot}(A)$. It is easy to check that this map is an
isomorphism. A similar relation between $DT^\hdot(A)$ and
$DP^{\hdot,\hdot}(A)$, formulated in terms of eulerian idempotents,
was considered in \cite{Fr}.

For every vector space $A$ we have $DT^1(A) \cong \Hom(A \otimes
A,A)$. An element $m \in \Hom(A \otimes A,A)$ satisfies $[m,m]=0$ if
and only if the corresponding binary operation $A \otimes A \to A$
is associative. In this case, one defines a differential
$DT^\hdot(A) \to DT^{\hdot+1}(A)$ and obtains the  {\em
Hochschild cochain complex} $DT^\hdot(A)$. Its cohomology groups
are called the {\em Hochschild cohomology} and denoted by
$\HHH^\hdot(A)$. If
the associative algebra $A$ is commutative, then the Hochshchild
differential $DT^\hdot(A) \to DT^{\hdot+1}(A)$ preserves the PBW
filtration. The complex $\gr^\idot DT^\hdot(A)$ with induced differential
is then isomorphic to the Poisson cohomology bicomplex
$HP^\hdot(\overline{A})$ of the Poisson algebra $\overline{A}$, which is the
commutative algebra $A$ equipped with trivial Poisson bracket
$\{-,-\}=0:\Lambda^2(A) \to A$. Note that the PBW filtration on
the complex $DT^\hdot(A)$ gives rise to a  PBW filtration
$F_\idot^{\tt{PBW}}\HHH^\hdot(A)$
on the Hochschild cohomology.

We will use the relationship between Hochschild and Poisson
complexes to prove the K\"unneth formula for Poisson cohomology. To
do this, we recall that for every two associative unitary algebras  $A$,
$B$,  there exists a
canonical quasiisomorphism
\begin{equation}
\sh:DT^\hdot(A \otimes B) \qis DT^\hdot(A) \otimes DT^\hdot(B),
\end{equation}
given by the  {\em shuffle  product}, see \cite[\S4.2]{Lo}. In more detail,
consider the category $\Delta$ of finite linearly ordered sets, and
denote by $[n] \in \Delta$, $n \geq 1$ the set of cardinality
$n$. By an {\em $(l,n)$-shuffle} we will understand an
order-preserving embedding $\Phi:[l] \to [n]$. Given an
$(l,k+l)$-shuffle $\Phi$, one defines the comlementary
$(k,k+l)$-shuffle $\overline{\Phi}:[k] \to [k+l]$. Taken together,
$\Phi$ and $\overline{\Phi}$ define a permutation $\sigma_\Phi:[l]
\cup [k] \to [k+l]$ of the set of $k+l$ elements. For every
shuffle $\Phi:[l] \to [k+l]$, we define a map $\sh_\Phi:A^{\otimes
l} \otimes B^{\otimes k} \to (A \otimes B)^{\otimes k + l}$ by
\begin{align*}
&\sh_\Phi(a_1 \otimes \dots \otimes a_l \otimes b_1 \otimes \dots
\otimes b_k) = c_1 \otimes \dots \otimes c_{k+l},\\ 
&c_i =
\begin{cases}
a_p \otimes 1, \quad &i = \Phi(p), p \in [l],\\
1 \otimes b_q, \quad &i = \overline{\Phi}(q), p \in [k],
\end{cases}
\end{align*}
The map $\sh:DT^{k+l}(A \otimes B) \to DT^l(A) \otimes DT^k(B) \cong
\Hom(A^{\otimes l} \otimes B^{\otimes k}, A \otimes B)$ is given by
\begin{multline}\label{shuffle}
\sh(f)(a_1 \otimes \dots \otimes a_k \otimes b_1 \otimes \dots
\otimes b_l) \\
= \sum_{\Phi:[l] \to
[k+l]}\sgn(\sigma_\Phi)f(\sh_\Phi(a_1 \otimes \dots \otimes a_k
\otimes b_1 \otimes \dots \otimes b_l))
\end{multline}
for every $f \in DT^{k+l}(A \otimes B) = \Hom((A \otimes B)^{\otimes
k+l},A \otimes B)$. 

In addition to the quasiisomorphism $\sh$, we also have an embedding
$
\kappa_A:DT^\hdot(A) \to DT^\hdot(A \otimes B)
$
given by
$$
\kappa_A(f)(a_1 \otimes b_1 \otimes \dots \otimes a_k \otimes b_k) =
f(a_1 \otimes \dots \otimes a_k) \otimes b_1\cdot \ldots\cdot b_k\;,
\quad f \in DT^k(A),
$$
and an analogously defined embedding
$\kappa_B:DT^\hdot(B) \to DT^\hdot(A \otimes B)$. Looking at the
formula \eqref{gerst} for the Gerstenhaber bracket, we immediately
see that both $\kappa_A$ and $\kappa_B$ are Lie algebra maps. Moreover,
say that a cochain $\om \in DT^k(A)$ is {\em reduced} if $f(a_1
\otimes \dots \otimes a_k) = 0$ whenever at least one of the
elements $a_1,\dots,a_k \in $ is equal to $1$. Then one can easily
derive from \eqref{gerst} and \eqref{shuffle} that for every reduced
cochain $\om \in DT^\hdot(A)$ and an arbitrary cochain $f \in
DT^\hdot(A \otimes B)$ we have
\begin{equation}\label{tau}
\sh([\kappa_A\om,f]) = [\om,\sh(f)]
\end{equation}
(here the bracket on the right-hand side acts on the first factor in
$DP^\hdot(A) \otimes DP^\hdot(B)$). The same statement holds for
$A$ replaced by $B$.

If the algebras $A$ and $B$ are commutative, then the
quasiisomorphism \eqref{shuffle} is compatible with the PBW
filtrations. The associated graded  map
\begin{equation}\label{knth}
\sh:DP^\hdot(\overline{A}) \otimes DP^\hdot(\overline{B}) \to
DP^\hdot(\overline{A} 
\otimes \overline{B})
\end{equation}
is also a quasiisomorphism; it is induced by the standard
quasiisomorphism
$$
\left(\Omega_\idot(A) \otimes B\right) \oplus \left(A \otimes
\Omega_\idot(B)\right) \qis \Omega_\idot(A \otimes B)
$$
between cotangent complexes.

\begin{lemma}\label{kun.alg}
For every two Poisson algebras $A$, $B$, the shuffle map
\eqref{shuffle} induces a natural
quasiisomorphism
$$
\sh:DP^\hdot(A \otimes B) \qis DP^\hdot(A) \otimes DP^\hdot(B).
$$
\end{lemma}

\proof{} It suffices to prove that the map $\sh$ is compatible with
the differentials; since the associated graded map \eqref{knth} is a
quasiisomorphism, $\sh:DP^\hdot(A \otimes B) \qis DP^\hdot(A)
\otimes DP^\hdot(B)$ will be a quasiisomorphism as well. By
definiton, the differential in the Poisson cohomology complex is the
sum of two differential, which we denoted by $d$ and
$\delta$. Compatibility with $d$ is contained in \eqref{knth}. Thus
we have to prove that
$$
\sh(\delta_{A \otimes B}f) = (\delta_A \otimes \id + \id \otimes
\delta_B)\sh(f)
$$
for every cochain $f \in DP^\hdot(A \otimes B)$. But by definition,
the Poisson cochain $\Theta_{A \otimes B}$ is of the form
$$
\Theta_{A \otimes B} = \kappa_A\Theta_A + \kappa_B\Theta_B.
$$
Since both Poisson cochains $\Theta_A$, $\Theta_B$ are clearly
reduced, we can apply \eqref{tau} and obtain
$$
\sh([\Theta_{A \otimes B},f]) = [\Theta_A,\sh(f)] +
[\Theta_B,\sh(f)].
$$
By \eqref{hp.def.eq2}, this is precisely what we had to prove.
\endproof

\subsection{Graded algebras.} As in
Subsection~\ref{poidef.defn.sub}, by a {\em graded} Poisson algebra
of degree $l$ we will understand a Poisson algebra $A$ equipped with
a grading such that the multiplication in $A$ is compatible with the
grading, and the Poisson bracket is of degree $-l$.  For a graded
Poisson algebra $A$ of degree $l \neq 0$, the Poisson cocycle
$\Theta \in DP^2(A)$ is canonically a coundary. Indeed, the grading
$A = \bigoplus A^\hdot$ induces a canonical derivation $\xi:A \to A$
by setting
$
\xi = k\id \text{  on  } A^k.
$
The formula $c\mapsto [\xi,c]$ extends this
 derivation  to a derivation of the Lie algebra 
$DP^{\hdot,\hdot}(A)$.
Since the Poisson bracket is of degree $-l$ with respect to the
grading, we have $[\xi,\Theta] = -l\cdot\Theta$, which can be rewritten
as $\Theta = -\frac{1}{l}d\xi$. 

The grading on $A$ induces an additional
grading on the Poisson cohomology bicomplex $DP^{\hdot,\hdot}(A)$,
called {\em the $A$-grading}. The differential
$d:DP^{\hdot,\hdot}(A) \to DP^{\hdot,\hdot+1}(A)$ preserves the
$A$-grading, and the differential $\delta:DP^{\hdot,\hdot}(A) \to
DP^{\hdot+1,\hdot}(A)$ shifts it by $l$. It will be convenient to
redefine the $A$-grading by shifting it by $(k-1)l$ on
$DP^{k,\hdot}(A)$, so that it is preserved by both differentials
$d,\delta$. We obtain a Poisson cohomology graded bicomplex
$DP^{\hdot,\hdot}(A)$. Note that after the shift, the Poisson
cochain $\Theta \in DP^{2,\hdot}(A)$ has $A$-degree $0$, and the
$A$-grading becomes compatible with the Lie bracket on
$DP^{\hdot,\hdot}(A)$.

\subsection{Modules and \'etale descent.}
Let $A$ be a Poisson algebra. By a {\em Poisson module} $M$ over $A$
we will understand an $A$-module $M$ equipped with an additional
operation $\{-,-\}:A \otimes M \to M$ such that
\begin{align*}
&\{ab,m\} = a\{b,m\}+b\{a,m\}, \qquad
\{a,bm\} = \{a,b\}m + b\{a,m\}, \\
&\{\{a,b\},m\} = \{a,\{b,m\}\} - \{b,\{a,m\}\}.
\end{align*}
Given a Poisson module $M$ one defines, following \cite{Fr},
its  Poisson cohomology bicomplex  $DP^{\hdot,\hdot}(A,M)$
(with coefficients in
$M$) by
$$
DP^{p,q}(A,M) = \Hom(P_{p,q}(A),M), \qquad p,q \geq 0,
$$
To get the differentials, it is convenient to 
treat  $\wh{A} :=
A \oplus M$ as a Poisson algebra by letting both operations to be zero on $M \subset
\wh{A}$ (the trivial square-zero extension). We grade this algebra by
assigning $\deg A = 0$, $\deg M = 1$, so that  $\wh{A}$ becomes a graded Poisson
algebra of degree $l = 0$. The differentials  on  $DP^{\hdot,\hdot}(\wh{A})$
preserve the induced grading, and  $DP^{\hdot,\hdot}(A,M)$
is nothing but the graded component in  $DP^{\hdot,\hdot}(\wh{A})$
of minimal degree.

In the special case $M = A$, we get back the original Poisson cohomology
bicomplex $DP^{\hdot,\hdot}(A)$. Indeed, the algebra $\wh{A} \cong
A\langle \eps \rangle = A[\eps]/(\eps^2)$ is simply the truncated
polynomial algebra over $A$. The bicomplex $DP^{\hdot,\hdot}(A)$
naturally lies in the degree-$0$ part of $DP^{\hdot,\hdot}(\wh{A})$,
and it is easy to check that, for $M=A$,
multiplication by $\eps$ provides an
isomorphism
$
DP^{\hdot,\hdot}(A,M) \cong DP^{\hdot,\hdot}(A).
$

Cohomology with coefficients appear naturally if
one wants to consider the functoriality properties of Poisson
cohomology groups. Let $f:A \to B$ be a Poisson map between Poisson
algebras. Then there is no natural map betwen the groups
$HP^\hdot(A)$, $HP^\hdot(B)$ -- just as there is no
natural map between the spaces of derivations of the algebras $A$,
$B$. However, there exist obvious natural maps
$f:DP^{\hdot,\hdot}(A) \to DP^{\hdot,\hdot}(A,f_*B)$ and
$df:DP^{\hdot,\hdot}(B) \to DP^{\hdot,\hdot}(A,f_*B)$. If
$\Theta_A$, $\Theta_B$ are the Poisson cochains of the algebra $A$,
$B$, then we have $f(\Theta_A) = df(\Theta_B)$. More generally, for
every Poisson module $M$ over $B$, there exists a natural map
$f:DP^{\hdot,\hdot}(B,M) \to DP^{\hdot,\hdot}(A,f_*M)$.

The main reason we want to consider functoriality is the following
compatibility result.

\begin{lemma}\label{loc}
Let $A \to B$ be an \'etale map between Poisson algebras, and let
$M$ be a Poisson module over the algebra $B$. Then the natural map
\begin{equation}\label{lc}
f:DP^{k,\hdot}(B,M) \to DP^{k,\hdot}(A,f_*M)
\end{equation}
induced by the map $A \to B$ is a quasiisomorphism for every $k \geq
0$. Consequently, the induced map $DP^\hdot(B,M) \to
DP^\hdot(A,f_*M)$ of total complexes is a quasiisomorphism.
\end{lemma}

\proof{} By \eqref{hp.def.eq}, it suffices to prove that the
canonical map
$$
\Har_\idot(A) \otimes B \to \Har_\idot(B)
$$
of Harrison complexes is a quasiisomorphism -- in other
words, that the pullback $\Omega_\idot(A) \otimes B$ of the
cotangent complex $\Omega_\idot(A)$ of the algebra $A$ is naturally
quasiisomorphic to the cotangent complex $\Omega_\idot(B)$ of the
algebra $B$. This is very well-known\footnote{In degree $0$, this is
the definition of an \'etale map.}.
\endproof

Notice that for every Poisson algebra $A$ and every multiplicative
system $S \subset A$, the localization $A[S^{-1}]$ carries a natural
Poisson structure. The natural map $A \to A[S^{-1}]$ is an example
of a Poisson \'etale map. This example will be used later in
constructing Poisson cohomology complexes of Poisson schemes.

\subsection{Poisson schemes.} We  define a {\em Poisson scheme}  as a
scheme $X$ equipped with a Poisson bracket on the structure sheaf
$\calo_X$. For any Poisson algebra $A$, the affine scheme $X = \Spec
A$ is a Poisson scheme. Poisson morphisms between Poisson schemes,
sheaves of Poisson $\calo_X$-modules and complexes of such sheaves
are defined in the natural way. For every Poisson module $M$ over a
Poisson algebra $A$, the corresponding coherent sheaf $\M$ on $X =
\Spec A$ is a sheaf of Poisson modules. Moreover, for every Poisson
morphism $f:X \to Y$ and a sheaf $\F$ of Poisson $\calo_X$-modules
on $X$, the direct image sheaf $f_*\F$ is a sheaf of Poisson
$\calo_Y$-modules. The same is true for the higher direct image
sheaves $R^pf_*\F$, $p \geq 1$, and for the whole complex $R^\hdot
f_*\F$.

The construction of the Poisson cohomology complex $DP^\hdot(A)$
generalizes straightforwardly to the scheme case. Namely, we notice
that the functors $P_k$, $k \geq 0$ from the category of vector
spaces into itself can be easily defined in an arbitrary symmetric
tensor category with the unit object. The terms $DP^k(A) =
\Hom(P_k(A),A)$ of the Poisson cohomology complex can be defined in
any symmetric monoidal category which admits internal
$\hhom$'s. Moreover, the differential in $DP^\hdot(A)$ for a Poisson
algebra $A$ is defined essentially by linear algebra, and  the
definition also works just as well in the general categorical
setting.  The same is true for the cohomology complex with
coefficients $HP^\hdot(A,M)$. We use this and define the {\em
Poisson cohomology complex} $\HP^\hdot(X,\F)$ of the Poisson scheme
$X$ with coefficients in an {\em injective} sheaf $\F$ of Poisson
$\calo_X$-modules by setting
$$
\HP^k(X) = \hhom(P_k(\calo_X),\F),
$$
where all tensor products and the $\hhom$ are taken in the category
of Zariski sheaves of groups on $X$. This complex carries the same
combinatorial structures as in the case of algebras, in particular,
the PBW filtration.

This definition makes sense for an arbitrary sheaf $\F$ of Poisson
modules -- in particular, for the structure sheaf $\calo_X$ -- but
it gives the wrong result: the resulting functor does not have
reasonable exactness properties, and it becomes impossible to prove
a crucial compatibility condition (Proposition~\ref{filtr.loc}
below). To get the correct definition, we use canonical injective
resolutions. Namely, let $\overline{X}$ be the union of all
scheme-theoretic points of $X$, and let $j:\overline{X} \to X$ be
the canonical map. Recall that for every sheaf $\F$ of abelian
groups on $X$, the canonical map
$
\F \to j_*j^*\F
$
is an embedding, and the sheaf $j_*j^*\F$ is injective. Iterating
this construction, one obtains the so-called {\em canonical Godement
resolution $\F^\hdot$} of the sheaf $\F$. If $\F$ is a sheaf of
Poisson $\calo_X$-modules, then all the terms in the Godement
resolution carry canonical structures of Poisson $\calo_X$-modules.

\begin{defn}\label{poi.coho.shf.defn}
Let $\F$ be a sheaf of Poisson $\calo_X$-modules on a Poisson scheme
$X$. The {\em local Poisson cohomology complex} $\HP^\hdot(X,\F)$ is
the total complex of the bicomplex
\begin{equation}\label{poi.coho.shf}
\HP^\hdot(X,\F^\hdot) = \hhom(P_k(\calo_X),\F^\hdot),
\end{equation}
where $\F^\hdot$ is the canonical Godement resolution of the sheaf
$\F$. The {\em Poisson cohomology}
$$
HP^\hdot(X,\F) = \HH^\hdot(X,\HP^\hdot(X,\F))
$$
with coefficients in $\F$ is the hypercohomology of the local
Poisson cohomology complex $\HP^\hdot(X,\F)$.
\end{defn}

We note that the functors $P_\idot(\calo_X) = P_{\idot,\idot}(X)$
carry an additional grading, so that in fact the local Poisson
cohomology complex $\HP^\hdot(X,\F) = \HP^{\hdot,\hdot}(X,\F)$ is a
bicomplex. For every $k \geq 0$, we have a quasiisomorphism
$$
\HP^{k,\hdot}(X,\F) \qis \RHom^\hdot(\Lambda^k\Omega_\idot(X),\F),
$$
where $\Omega_\idot(X)$ is the cotangent complex of the scheme $X$,
and $\Lambda^k$ is understood in the derived-category
sense. Moreover, the complex $\HP^{k,\hdot}(X,\F)$ is concentrated
in degrees $\geq 0$. Therefore we obtain a canonical embedding
$$
\Hom(\Lambda^k\Omega(X),\calo_X) \to \HP^{k,\hdot}(X,\F),
$$
where $\Omega(X)$ is the cotangent sheaf.

Since the Godement resolution is functorial,
Definition~\ref{poi.coho.shf.defn} immediately generalizes to the
case when $\F$ is itself not a single sheaf of Poisson
$\calo_X$-modules, but a complex of such sheaves. Then it descends
to a triangulated functor on the corresponding derived category. In
the case when $\F=\calo_X$ is the structure sheaf, we will denote
the local Poisson cohomology complex $\HP^\hdot(X,\calo_X)$ simply
by $\HP^\hdot(X)$, and we will denote the Poisson cohomology
$HP^\hdot(X,\calo_X)$ by $HP^\hdot(X)$.

Another important situation is the following. Let $Z \subset X$ be a
closed subscheme in a Poisson scheme $X$, and let $U \subset X$ be
the open complement. Let $i:Z \hookrightarrow X$, $j:U
\hookrightarrow X$ be the embeddings. Then there exists a canonical
map
$
\calo_X^\hdot \to j_*\calo_U^\hdot,
$ 
where $\calo_X^\hdot$ is the Godement resolution of the sheaf
$\calo_U$, and $\calo_U^\hdot$ is the Godement resolution of the
sheaf $\calo_U$. This map is compatible with the natural Poisson
structures on both sides. Therefore we also have a Poisson module
structure on its cone, denoted by $i_!\calo_Z^\hdot \subset
\calo_X^\hdot$ (this is a particular model of the canonical
quasicoherent complex $i_!\calo_Z$ used in the local duality
theory).

\begin{defn}
The Poisson cohomology 
$
HP^\hdot_Z(X) = HP^\hdot(X,i_!\calo_Z^\hdot)
$ 
is called the {\em Poisson cohomology of the scheme $X$ with
supports in $Z \subset X$}.
\end{defn}

Note that the closed subscheme $Z \subset X$ enters into this
construction only through its open complement $U \subset X$. In
particular, there is no need to assume that $Z \subset X$ is a
Poisson subscheme in any sense. We will need the following easy
vanishing result.

\begin{lemma}\label{cm.i}
Assume that the Poisson scheme $X$ is Cohen-Macaulay, and the closed
subscheme $Z \subset X$ is of codimension $\codim Z \geq k$ for some
integer $k \geq 2$. Then $HP^i_Z(X)=0$ for 
all $i\leq k-2$.
\end{lemma}

\proof{} Since $X$ is Cohen-Macaulay, the complex $i_!\calo_Z$ is
trivial in degrees $\leq k-1$.
\endproof

\begin{prop}\label{filtr.loc}
Let $X = \Spec A$ be a Poisson scheme with a Poisson sheaf $\F$ of
$\calo_X$-modules.
Then the  canonical map
$
HP^\hdot(X,\F) \to HP^\hdot(A,\RG^\hdot(X,\F))
$  is a quasiisomorphism.
\end{prop}

\proof{} The claim is functorial in $\F$. By
Definition~\ref{poi.coho.shf.defn}, it suffices to consider sheaves
of the form $\F = i_*M$, where $i:x \hookrightarrow X$ is the
embedding of a scheme-theoretic point $x \in X$, and $M \cong
\RG^\hdot(A,\F)$ is a Poisson module over the local ring
$\calo_{X,x}$. Let $\F = i_*M$ be such a sheaf. Then by definition
the local Poisson cohomology complex $\HP^\hdot(X,\F)$ is
quasiisomorphic
$$
\HP^\hdot(X,\F) \cong i_*HP^\hdot(\calo_{X,x},M),
$$ 
so that $HP^\hdot(X,\F) \cong HP^\hdot(\calo_{X,x},M)$. To prove the
claim, it suffices to note that the natural map $A \to \calo_{X,x}$
is \'etale, and apply Lemma~\ref{loc}.  \endproof

\begin{prop}\label{geom}\mbox{}
\begin{enumerate}
\item For every two Poisson schemes $X$, $Y$, there exists a
canonical quasiismorphism
\begin{equation}\label{kun.geom}
HP^\hdot(X \times Y) \cong HP^\hdot(X) \otimes HP^\hdot(Y).
\end{equation}
\item Let $f:X \to Y$ be an \'etale map of Poisson schemes. Then for
every complex $\F^\hdot$ of sheaves of Poisson modules on $X$, there
exists a canonical quasiisomprhism
$$
HP^\hdot(Y,f_*\F^\hdot) \cong HP^\hdot(X,\F^\hdot).
$$
\item Let $f:X \to Y$ be an \'etale Galois cover with Galois group
$G$. Then there exists a canonical quasiisomorphism
$
HP^\hdot(Y) \cong \left(HP^\hdot(X)\right)^G.
$
\end{enumerate}
\end{prop}

\proof{} By taking affine covers, we reduce \thetag{i} and
\thetag{ii} to the case of affine schemes $X$, $Y$. Then
Proposition~\ref{filtr.loc} reduces both statements to their
algebraic versions Lemma~\ref{kun.alg} and Lemma~\ref{loc}. Finally,
\thetag{iii} is a direct corollary of \thetag{ii}: we have
$$
\left(HP^\hdot(X)\right)^G \cong
\left(HP^\hdot(Y,f_*\calo_X)\right)^G \cong
HP^\hdot(Y,(f_*\calo_X)^G) \cong HP^\hdot(Y,\calo_Y). \quad\square
$$

\begin{corr}\label{supp.geom}\mbox{}
\begin{enumerate}
\item Let $f:X \to Y$ be a map between Poisson schemes, and assume
the restriction of $f$ to a closed subscheme $Z \subset X$ is
injective. Then, there exists a canonical quasiisomorphism
$$
HP^\hdot_Z(Y) \to HP^\hdot_Z(X).
$$
\item Let $Z \subset X$ be a closed subscheme in a Poisson scheme
$X$, and let $Y$ be another Poisson scheme. Then there exists a
canonical isomorphism
\begin{equation}\label{kun.supp} 
HP^\hdot_{Z \times Y}(X \times Y) \cong HP^\hdot_Z(X) \otimes
HP^\hdot(Y).
\end{equation}
\end{enumerate}
\end{corr}

\proof{} \thetag{ii} immediately follows from
Proposition~\ref{geom}\thetag{i}. To prove \thetag{i}, let $i^X:Z
\hookrightarrow X$, $i^Y = f \circ i^X:Z \hookrightarrow Y$ be the
embeddings. Notice that the map $f$ must be \'etale in a Zariski
neighborhood of the subscheme $Z \subset X$, and we have
$f_*i^X_!\calo_Z \cong i^Y_!\calo_Z$. The claim then follows from
Proposition~\ref{geom}\thetag{ii}.  \endproof

\subsection{The smooth case.}\label{smooth_case} Let now $X$ be a {\em smooth} Poisson
scheme. In this case, we have the original Koszul-Brylinski formalism for
local Poisson cohomology. Namely, one considers the graded Lie
algebra $\Lambda^\hdot\T(X)$ of polyvector fields on $X$ with the
so-called {\em Schouten bracket}. The Poisson structure is defined
by a bivector field $\Theta \in \Lambda^2\T(X)$ such that
$[\Theta,\Theta] = 0$. Taking commutator with $\Theta$ defines a
differential $d:\Lambda^\hdot\T(X) \to \Lambda^{\hdot+1}\T(X)$; we
will call the complex $\HP^\hdot_\Theta(X) \cong \langle
\Lambda^\hdot\T(X),d \rangle$ the {\em Koszul-Brylinski complex}.

For a general smooth Poisson scheme $X$, we have canonical
embeddings
$$
\Lambda^k\T(X) \to \HP^{k,\hdot}(X).
$$
These embeddings are compatible with the Koszul-Brylinski differential and
define a quasiisomorphism $\HP^\hdot_\Theta(X) \to \HP^\hdot(X)$. Indeed,
since the embeddings are canonical, it suffices to check this fact
locally, so that we can assume that the scheme $X = \Spec A$ is
affine. Then since $A$ is smooth, the Harrison homology
$\Har_\idot(A)$ is quasiisomorphic to the module $\Omega^1A$ of
K\"ahler differentials of the algebra $A$, and this module is flat
over $A$. Therefore the complex $DP^{k,\hdot}(A)$ only has
non-trivial cohomology in degree $0$, and  $HP^{k,0}(A) \cong
\Lambda^k\T(A)$, cf. \cite{Fr}. We have $DP^{k,l}(A) = 0$ for $l < 0$ and a
canonical embedding $\Lambda^k\T(A) \cong HP^{k,0}(A) \to
DP^{k,0}(A)$
 sends the
Poisson bivector $\Theta$ into the Poisson cochain, and it also
sends the Schouten bracket in $\Lambda^\hdot\T(A)$ into the
Gerstenhaber bracket in $DP^\hdot(A)$. Therefore it sends the
Koszul-Brylinski differential into the differential
$\delta:DP^{\hdot,\hdot}(A) \to DP^{\hdot+1,\hdot}(A)$, see also
\cite{Fr}. 

Recall that in case of a smooth affine Poisson variety $X=\Spec A$ there are
Hochschild-Kostant-Rosenberg isomorphisms
$$\HHH^\hdot(A)\cong \Lambda^\hdot_A(\Har_1(A))\cong \Lambda^\hdot_A\T(A)\,,
$$
and  the Koszul-Brylinski complex of $X$ may  be identified
via this isomorphism  with the `twisted'  Hochschild complex of $A$ considered
in \S6, cf. Definition \ref{twisted_Hoch}. Thus,
we have 
\begin{equation}\label{PBW_smooth}
HP^\hdot_\Theta(X)\cong\HHH^\hdot_\Theta(A)\,.
\end{equation}
Also,  in the smooth case, 
for the  PBW filtration on  Hochschild cohomology we get
$$
F_i^{\tt{PBW}}\HHH^j(A)=
\begin{cases}
\HHH^j(A) & \text{ if } i\geq j\\
0 & \text{ if } i <j .
\end{cases}
$$
Furthermore, 
the spectral sequence associated to the  bicomplex
\eqref{bicomplex}
collapses.

If the Poisson scheme $X$ is not only smooth but also symplectic,
then the symplectic form provides an identification $\T(X) \cong
\Omega(X)$. This identification extends to an isomorphism
$\HP^\hdot_\Theta(X) \cong \Omega^\hdot(X)$, the standard de Rham complex of
the smooth scheme $X$. Thus for symplectic schemes we have
$$
\HP^\hdot(X) \qis \HP^\hdot_\Theta(X) \cong \Omega^\hdot(X)
\quad\text{and}\quad
HP^\hdot(X) \cong H^\hdot(X,\C),
$$
the ordinary singular cohomology of the scheme $X$. We can also
obtain isomorphisms for cohomology with coefficients: Poisson
modules over $X$ are the same as $D$-modules, and the Poisson
cohomology $HP^\hdot(X,\F)$ with coefficients in a Poisson module is
isomorphic to the de Rham cohomology with coefficients in the
corresponding $D$-module. In particular, the Poisson cohomology with
supports is isomorphic to the singular cohomology with supports.

\subsection{Deformations.}\label{def.sub}
The shortest route to the classification of Poisson deformations
(in particular,  to the proof of Theorem~\ref{poi.def.uni}) is
through the formalism  of deformation
groupoids. This consists of obtaining deformed structures from 
solutions of the so-called Maurer-Cartan equation $d\gamma =
\{\gamma,\gamma\}$ in a certain differential-graded Lie algebra
$\LL^\hdot$.

More precisely, assume given a DG Lie algebra $\LL^\hdot$ with
differential $d:\LL^\hdot \to \LL^{\hdot+1}$ and commutator
$\{-,-\}$, and a local Artin algebra $S$ with maximal
ideal $\m \subset S$. Consider the DG Lie algebra $\LL^\hdot \otimes
S$ over $S$. Set $\g = \LL^0$. For any $g \in \g \otimes \m$, $a \in
\LL^1 \otimes S$ set 
$$
g \cdot a = dg + [g,a] \in \LL^1 \otimes S.
$$
Consider $g \cdot a$ as a tangent vector to the vector space $\LL^1
\otimes S$ at the point $a \in \LL^1 \otimes S$. Then the collection
$g \cdot a$, $a \in \LL^1 \otimes S$ defines a vector field on
$\LL^1 \otimes S$ for every $g \in \g \otimes \m$. These vector
fields glue together to an action of the Lie algebra $g \otimes \m$
on $\LL^1 \otimes S$. Since the Lie algebra $\g \otimes \m$ is
nilpotent, this action extends to an action of the corresponding
nilpotent Lie group $G_\m$.

\begin{defn}\label{MC}
The {\em Maurer-Cartan groupoid} $MC(\LL^\hdot,S)$ associated to the
pair $\langle \LL^\hdot,S\rangle$ is defined in the following way:
\begin{enumerate}
\item Objects of $MC(\LL^\hdot,S)$ are elements $\gamma \in \LL^1
\otimes \m$ satisfying the Maurer-Cartan equation
\begin{equation}\label{mc}
d\gamma = \{\gamma,\gamma\}.
\end{equation}
\item Morphisms between objects $\gamma_1,\gamma_2 \in
MC(\LL^\hdot,S)$ are elements $g \in G_\m$ in the group $G_\m$ such
that $g \cdot \gamma_1 = \gamma_2$.
\end{enumerate}
\end{defn}

To check that this definition is consistent, one has to prove that
the $G_\m$-action preserves the Maurer-Cartan equation, which is an
elementary computation.

The Maurer-Cartan formalism admits an obvious graded version. If the
DG Lie algebra $\LL^\hdot$ is equipped with a grading, -- or,
equivalently, witha $\C^*$-action, -- then one can consider local
Artin algebras $S$ equipped with a $\C^*$-action and form the
grouppoid $MC_{gr}(\LL^\hdot,S)$ of the $\C^*$-invariant solutions
to the Maurer-Cartain equation \eqref{mc}.

Application of the Maurer-Cartan formalism to the Poisson
deformation theory is an immediate corollary of Lemma~\ref{compt}
and Lemma~\ref{har}.

\begin{lemma}\label{def.poi}
Let $A$ be a Poisson algebra over $\C$. Consider the Poisson
cohomology complex $DP^\hdot(A)$, and let 
$$
\LL^k = \begin{cases}DP^{k+1}(A), &\quad k \geq 0,\\0, &\quad k <
0.\end{cases}
$$
Let $S$ be a local Artin algebra with maximal ideal $\m$. 

Then the category $\Def(A,S)$ of flat Poisson algebras $\wt{A}$ over
$S$ equipped with an isomorphism $\wt{A}/\m\cdot\wt{A} \cong A$ is
equivalent to the Maurer-Cartan groupoid $MC(\LL^1,S)$.
\end{lemma}

\proof{} By Lemma~\ref{compt} and Lemma~\ref{har}, defining a
Poisson algebra structure on a vector space $V$ is equivalent to
giving a cochain $\alpha = m+\Theta \in DP^2(V) = DL^1(V) \oplus
DP^{2,0}(V)$ which satifies $[\alpha,\alpha]=0$. To define a functor
$MC(\Der^\hdot(A),S) \to \Def(A,S)$, it suffices to notice that for
every solution $\gamma \in \LL^1 = DP^2(A) \otimes \m$ of the
Maurer-Cartan equation, the element
$$
\wt{\alpha} = \alpha - \mbox{$\frac{1}{2}$}\gamma \in \LL^1 \otimes S
$$
satisfies $\{\wt{\alpha},\wt{\alpha}\}=0$, hence defines a Poisson
algebra structure on $A \otimes S$. This functor is surjective on
objects -- indeed, by definition for every deformation $\wt{A}$ we
have $\wt{A} \cong A \otimes S$ as $S$-modules. To prove that it is
an equivalence, it suffices to prove that it is surjective on
morphisms. But the Lie algebra $\LL^0 = DP^1(A) \cong \Hom(A,A)$ by
definition consists of {\em all} linear endomorphisms of the vector
space $A$. Therefore morphisms in the Maurer-Cartan groupoid are all
$S$-module morphisms $A \otimes S \to A \otimes S$ which reduce to
identity modulo $\m \subset S$. This is the same as morphisms in
$\Def(A,S)$.  \endproof

For any DG Lie algebra $\LL^\hdot$, the Maurer-Cartan groupoids
$MC(\LL^\hdot,S)$ form a stack over the category of local Artin
algebras $S$. The general machinery shows that this stack is
(pro)-representable. The representing object is a certain
canonically defined subscheme in the vector space $H^1(\LL^\hdot)$
completed at zero, modulo the action of the (pro)-algebraic group
corresponding to the Lie algebra $\LL^0$. 

The machinery works in the case of graded DG Lie algebras without
any change: the (pro)-representing object carries a natural
$\C^*$-action compatible with the induced grading on
$H^1(\LL^\hdot)$, and it (pro)-represents the Maurer-Cartain stack
$MC_{gr}(\LL^\hdot,S)$ on the category of local Artin schemes
equipped with a $\C^*$-action. However, if the degrees of the
induced gradings on $H^1(\LL^\hdot)$, $H^0(\LL^\hdot)$ are strictly
positive, we can obtain more.

Say that a commutative algebra $A = \bigoplus A^\hdot$ is positively
graded if $\dim A^k = 0$ for $k < 0$, $\dim A^0 = 1$ and $\dim A^k <
\infty$ for $k > 0$. Further, we say that a complete local algebra
$\langle A,\m \rangle$ is equipped with a good $\C^*$-action if the
group $\C^*$ acts on $\wh{A}$, and this action induces a positive
grading on every Artin quotient $A/\m^k A$ and on the associated
graded quotient
$
\gr^\hdot\wh{A} = \oplus_k\,\m^k/\m^{k+1}.
$
For every positively graded algebra $A^\hdot$, denote by $\wh{A}$
the completion of the algebra $A^\hdot$ with respect to the maximal
ideal $\m=A^{\geq 1} \subset A^\hdot$. We note the following standard
fact (essentially, this is a toy version of the Grothendieck
algebraization, \cite[III, 5.4]{EGA}).

\begin{lemma}\label{compl.eq}
\begin{enumerate}
\item The correspondence $A^\hdot \mapsto \wh{A}$ is an equivalence
between the category of positively graded commutative algebras
$A^\hdot$ and the category of complete local algebras $\langle
\wh{A},\m \subset \wh{A}\rangle$ equipped with a good $\C^*$-action.

\item For any positively graded algebra $A^\hdot$, the completion
induces an equivalence between the category of finitely generated
graded $A^\hdot$-modules and the category of finitely generated
$\C^*$-equivariant $\wh{A}$-modules.

\item If the algebra $A^\hdot$ is finitely generated, then the
completion induces an equivalence between the category of projective
$\C^*$-equivariant schemes of finite type over $\Spec A^\hdot$ and
the category of projective $\C^*$-equivariant schemes of finite type
over $\Spec \wh{A}$.

\item Let $X$ be a projective $\C^*$-equivariant scheme of finite
type over $\Spec A^\hdot$, and let $\X/\Spec\wh{A}$ be its
completion. Then completion induces an equivalence between the
category of coherent $\C^*$-equivariant sheaves on $X$ and coherent
$\C^*$-equivariant sheaves on $\X$. A sheaf $\F$ on $X$ is flat if
and only if its completion $\wh{F}$ is flat. Moreover, for every
sheaf $\F$ on $X$, the completion $\wh{H^k(X,\F)}$ of the cohomology
group $H^k(X,\F)$ with respect to the $\m$-adic topology is
canonically isomorphic to the cohomology group $H^k(\X,\wh{\F})$ of
the completion $\wh{\F}$
\end{enumerate}
\end{lemma}

\proof[Sketch of proof.] To prove \thetag{i}, we will just give
the inverse to the equivalence of categories. It is given by
\begin{equation}\label{inverse}
\langle A,\m \rangle \mapsto \bigoplus\nolimits_p\, \lim_{\leftarrow}
(A/\m^k)^p. 
\end{equation}
By our definition of a good $\C^*$-action, for every integer $p$ the
projective system of the corresponding graded components
$(A/\m^k)^p$ actually stabilizes at some finite level $k \geq 0$,
and the limit is a finite-dimensional vector space. Thus the
right-hand side is indeed a positively graded commutative
algebra. To prove \thetag{ii}, we define the inverse equivalence by
the same formula \eqref{inverse}. 

To prove \thetag{iii}, we write $X = \Proj \wh{B}^\hdot$ for a
scheme $\wh{X}$ projective over $\Spec \wh{A}$, apply
\eqref{inverse} and obtain a scheme $X = \Proj B$ over $\Spec A$
whose completion is $\wh{X}$ (here $\wh{B}$ has two independent
gradings: the first one is related to the $\C^*$-action, and the
second one comes from the definition of $\Proj$). This defines the
inverse equivalence on the level of objects. To lift it to
morphisms, one identifies morphisms with their graphs.

Finally, to prove \thetag{iv} one uses Serre's Theorem and
identifies the category of coherent sheaves on a projective scheme
$B$ with the quotient of the category of finitely-generated graded
$B$-modules by the fat subcategory of finitely-generated $B$-modules
with finite number of non-trivial graded components. Again using
\eqref{inverse} -- and keeping track of the additional grading --
one checks easily that the completion is an equivalence both between
the categories of modules and between the fat subcategories one has
to mod out. Therefore it also induces an equivalence between the
quotient categories.

This equivalence is clearly compatible with the tensor products,
hence preserves flatness. To prove the statement on cohomology, it
suffices to identify the graded component $H^k(X,\F)^p$ of some
degree $p$ with the $\Ext$-group $\Ext^k(\calo_p,\F)$, where
$\calo_p$ is the structure sheaf on $X$ with grading shifted by $p$.
\endproof

If the degrees of the natural gradings on $H^1(\LL^\hdot)$,
$H^0(\LL^\hdot)$ are strictly positive, then the $\C^*$-action on the
universal object for the Maurer-Cartan stack is good in the sense of
Lemma~\ref{compl.eq}. Therefore this object falls within the scope
of the Lemma, and we automatically obtain a universal solution for
the Maurer-Cartan equation over a positively graded algebra $S$, not
only over its completion $\wh{S}$.

In the case of Poisson deformations of a Poisson algebra $A$, we
have $H^1(\LL^\hdot) \cong HP^2(A)$. The Lie algebra
$H^0(\LL^\hdot)$ is the algebra of all Poisson cocycles in $DP^1(A)
= \Hom(A,A)$ which coincides with the Lie algbera of all Poisson
derivations. This Lie algebra is rather big; in particular, it
contains all the Hamiltonian vector fields on $A$. However, we do
not need to know the precise stack structure on $MC(\LL^\hdot,S)$ --
our goal is the weaker Theorem~\ref{poi.def.uni}. To prove it, it
suffices to prove that the Lie algebra $H^0(\LL^\hdot)$ acts
trivially on the vector space $H^1(\LL^\hdot)$.

\proof[Proof of Theorem~\ref{poi.def.uni}.] Assume that
$HP^1(A)$. We have to prove that the Lie algebra $H^0(\LL^1)$ of
Poisson derivations of the algebra $A$ acts trivially on the
cohomology group $H^1(\LL^\hdot)$. Consider the map $\LL^\hdot \to
DP^{\hdot-1}(A)$. This is a Lie algebra map and a isomorphism in all
non-negative degrees. In particular, it induces an isomorphism
$H^1(\LL^\hdot) \cong HP^2(A)$, and the action of the algebra
$H^0(\LL^\hdot)$ on $HP^2(A) \cong H^1(\LL^\hdot)$ is induced by the
action of the Lie algebra $HP^1(A)$. By assumption, the Lie algebra
$HP^1(A)$ is trivial.  \endproof

\subsection*{A linear algebra result.}
Let $V$ be a finite dimensional
complex symplectic vector space, and $G\subset Sp(V)$ a finite subgroup.
Write $V^g$ for the fixed point space of
$g\in G$. 
\begin{lemma}\label{easy} For any elements $g,h\in G$ such that 
$V=V^g+V^h$, we have $V^{g\cdot h}= V^g\cap V^h$.
\end{lemma}
\noindent

\proof{} It is clear that $V^g\cap
V^h\subseteq V^{g\cdot h}$. To prove the opposite inclusion,
let $x\in V^{g\cdot h}$. The assumption
$V=V^g+V^h$, implies that $x=u+v$, for some $u\in V^g$ and $v\in V^h.$
Thus, $ gh(u+v)=u+v$, hence, $hu+v=u+g^{-1}v,$ or equivalently
$(h-1)u=(g^{-1}-1)v.$ But, for any $a\in Sp(V)$ of finite order, there is an orthogonal
direct sum decomposition
 $V={\tt{Image}}(\id-a) \oplus V^a$.
We deduce
that $(h-1)u$ in the LHS of the last equation is
orthogonal to $V^h$, resp.   $(g^{-1}-1)v$
in the RHS is
orthogonal to $V^g$, with respect to the symplectic form.
Thus, each side  is orthogonal to $V=V^g+V^h$, hence, vanishes.
It follows that $(h-1)u=(g^{-1}-1)v=0$. Therefore $u,v\in V^g\cap
V^h,$ hence $x\in V^g\cap
V^h,$ and we obtain $V^{g\cdot h}\subseteq V^g\cap
V^h.$~\endproof

\footnotesize{
{\bf V.G.}: Department of Mathematics, University of Chicago,
Chicago, IL
60637, USA;\\
\hphantom{x}\quad\, {\tt ginzburg@math.uchicago.edu}
\smallskip

{\bf D.K.}: Steklov Mathematical Institute,
Moscow, Russia\\
\hphantom{x}\quad\, {\tt kaledin@balthi.dnttm.ru}}

\end{document}